%% file: main.tex
\documentclass[]{article}
\pdfoutput=1
\usepackage{authblk}
\usepackage{setspace}
\usepackage[letterpaper, margin=1in]{geometry}

\usepackage[cmex10]{amsmath}
\usepackage{amsfonts, amssymb}
\usepackage{color}
\usepackage{bbm}
\usepackage{amsthm}
\usepackage{graphicx}
\usepackage{enumitem}
\graphicspath{{./figs/}}%{./orth-frame/figure/}{./bi-orth/figure/}}
\makeatletter
\newcommand*{\ov}[1]{\m@th\overline{\mbox{#1}\raisebox{.8em}{}}}
\newcommand{\V}[1]{\boldsymbol{#1}}% boldface for vector and matrix
\usepackage[toc,page]{appendix}
\usepackage{wrapfig}
\usepackage{color}
\usepackage{verbatim}
\usepackage{caption}

\newcommand{\mo}[1]{m_{#1}(\boldsymbol{\omega})}
\newcommand{\m}[1]{\widetilde{m_{#1}}(\boldsymbol{\omega})}
\renewcommand{\mp}[2]{\widetilde{m_{#1}}(\boldsymbol{\omega}+\boldsymbol{\pi}_{#2})}
\newcommand{\mc}[1]{\widetilde{m_{#1}}^C(\boldsymbol{\omega})}

\newcommand{\M}{\widetilde{\mathbf{M}}}
\newcommand{\Msub}{\widetilde{\mathbf{M}}^{\Box}}
\newcommand{\G}{\mathcal{G}}

\newcommand{\mhat}[1]{\widehat{\widetilde{m_{#1}}}(\boldsymbol{\omega})}
\newcommand{\xvec}{\mathbf{x}}
\newcommand{\wvec}{\mathbf{w}}
\newcommand{\mrow}[1]{\widetilde{\mathbf{m}}^{#1}}
\newtheorem{theorem}{Theorem}[section]
\newtheorem{lemma}[theorem]{Lemma}
\newtheorem{proposition}[theorem]{Proposition}

\renewcommand{\qed}{\hfill\ensuremath{\square}}
\newcommand{\mteven}{\overlinespace{\widetilde{\mathbf{m}_0}}^{\mathcal{E}}}
\newcommand{\meven}{\mathbf{m}_0^{\mathcal{E}}}
\newcommand{\mtodd}{\overlinespace{\widetilde{\mathbf{m}_0}}^{\mathcal{O}}}
\newcommand{\modd}{\mathbf{m}_0^{\mathcal{O}}}
\makeatletter

\newcommand{\mylabel}[2]{#2\def\@currentlabel{#2}\label{#1}}
\newcommand*{\overlinespace}[1]{\overline{\mbox{$#1$}\raisebox{3mm}{}}}
\newcommand{\sbarm}[1]{\overlinespace{\widetilde{m_{#1}}}(\boldsymbol{\omega})}
\newcommand{\sbarmp}[2]{\overlinespace{\widetilde{m_{#1}}}(\boldsymbol{\omega} + \boldsymbol{\pi}_{#2})}
\newcommand{\sbarmn}[2]{\overlinespace{\widetilde{m_{#1}}}(\boldsymbol{\omega} + \boldsymbol{\nu}_{#2})}
\makeatother

\begin{document}

\newtheorem{lem}{Lemma}
\newtheorem*{mydef}{Definition}
\newtheorem{thm}{Theorem}
\newtheorem{prop}{Proposition}
\newtheorem*{notat}{Notation}

\abovedisplayskip=2pt
\belowdisplayskip=2pt
\abovedisplayshortskip=2pt
\belowdisplayshortskip=2pt

\title{Directional Wavelet Bases Constructions \\with Dyadic Quincunx Subsampling}
\author[1]{Rujie Yin}
\author[2]{Ingrid Daubechies}
\affil[1]{Department of Mathematics, Duke University, rujie.yin@duke.edu}
\affil[2]{Department of Mathematics, Duke University, ingrid@math.duke.edu}
\date{}
\maketitle

\begin{abstract}
We construct directional wavelet systems that will enable building efficient signal representation schemes with good direction selectivity. In particular, we focus on wavelet bases with dyadic quincunx subsampling.
%; dilated quincunx downsampling is used to construct orthonormal and bi-orthogonal bases and standard dyadic downsampling for low-redundancy frames. 
In our previous work \cite{yin2014orthshear}, We show that the supports of orthonormal wavelets in our framework are discontinuous in the frequency domain, yet this irregularity constraint can be avoided in frames, even with redundancy factor less than 2. In this paper, we focus on the extension of orthonormal wavelets to biorthogonal wavelets and show that the same obstruction of regularity as in orthonormal schemes exists in biorthogonal schemes. In addition, we provide a numerical algorithm for biorthogonal wavelets construction where the dual wavelets can be optimized, though at the cost of deteriorating the primal wavelets due to the intrinsic irregularity of biorthogonal schemes.
\end{abstract}

\input{intro}

\input{setup}

\input{orthonormal}

%\include{orthonormal}

\input{frame}

\input{bi-orth}

\input{conclusion}
%\include{conclusion}
%\section{Conclusion}\label{sec: end}

\textit{•}\section{Acknowledgment}
This work is support by the NSF grant 1516988.

\bibliographystyle{IEEETrans}%bib/te}
\bibliography{ref}

\clearpage
\newpage
\begin{appendices}
\input{parseval-cond}
\input{lemmas}
\input{QCQP}
\input{m0-opt}

\end{appendices}

\newpage
\end{document}

%% file: intro.tex
\section{Introduction}

\iffalse
\begin{itemize}
\item review construction of directional wavelet, shearlet
\item what's new in our construction?
\item summary: framework, technique, reference(Durand, Cohen)
\item organization of the paper
\end{itemize}
\fi

% introducing directional wavelets
In image compression and analysis, 2D tensor wavelet schemes are widely used. Despite the time-frequency localization inherited from 1D wavelet, 2D tensor wavelets suffer from poor orientation selectivity: only horizontal or vertical edges are well represented by tensor wavelets. To obtain better representation of 2D images, several directional wavelet schemes have been proposed and applied to image processing, such as directional wavelet filterbanks (DFB) and various extensions.

% summary
Conventional DFB \cite{DFB92} divides the square frequency domain associated with a regular 2D lattice into eight equi-angular pairs of triangles; such schemes can be critically downsampled (maximally decimated) with perfect reconstruction (PR), but they typically do not have a multi-resolution structure. 
Different approaches have been proposed to generalize DFB to multi-resolution systems, including non-uniform DFB (nuDFB), contourlets, curvelets, shearlets and dual-tree wavelets.
nuDFB is introduced in \cite{nuDFB05} based on multi-resolution analysis (MRA), where at each level of decomposition the square frequency domain is divided into a high frequency outer ring and a central low frequency domain. For nuDFB, the high frequency ring is primarily divided further into six equi-angular pairs of trapezoids and the central low frequency square is kept intact for division in the next level of decomposition, see the left panel in Figure \ref{fig: partition}. The nuDFB filters are solved by optimization which provides non-unique near orthogonal or bi-orthogonal solutions depending on the initialization without stable convergence.
Contourlets \cite{do2005contourlet} combine the Laplacian pyramid scheme with DFB which has PR but with redundancy $4/3$ inherited from the Laplacian pyramid.
Shearlet \cite{shearlet12book,easley2008sparse} and curvelet \cite{candes2006fast} systems construct a multi-resolution partition of the frequency domain by applying shear or rotation operators to a generator function in each level of frequency decomposition. Available shearlet and curvelet implementations have redundancy at least 4; moreover, the factor may grow with the number of directions in the decomposition level.
% to be modified and add M-band version
Dual-tree wavelets \cite{selesnick2005dual} are linear combinations of 2D tensor wavelets (corresponding to multi-resolution systems) that constitute an approximate Hilbert transform pair, where the high frequency ring is divided into pairs of squares of different directional preference.

% what's the problem we will address, don't focusing too much on shearlet for the moment
However, none of these multi-resolution schemes is PR, critically downsampled and regularized (localized in both time and frequency). In the framework of nuDFB (\cite{nuDFB05}), it was shown by Durand \cite{durand2007} that it is impossible to construct orthonormal filters localize without discontinuity in the frequency domain, or -- equivalently -- regularized filters without aliasing. His construction of directional filters uses compositions of 2-band filters associated to quincunx lattice, similar to that of uniform DFB in \cite{nuDFB05}; as pointed out in \cite{nuDFB05} the overall composed filters are not alias-free. It is not clear whether Durand's argument also precludes the existence of a regularized wavelet system, if one slightly weakens the set of conditions.
%Among these,  shearlets and curvelets have optimal asymptotic rate of approximation for ``cartoon images''(piecewise smooth, with jumps occurring along piecewise $C^2$-curves), due to the parabolic scaling rule in the frequency domain \cite{guo2007optimally,candes2005curvelet}; they have been successfully applied to image denoising \cite{easley2009shearlet}, restoration \cite{candes2002new} and separation \cite{kutyniok2012image}. Despite their theoretical potential, the (often high) redundancy of curvelets and shearlets impedes their practical usage. 
%Redundancy is useful in image processing tasks such as denoising, restoration and reconstruction, but a non-redundant basis decomposition is preferred in tasks where computation cost is of concern.

To study this question, we consider multi-resolution directional wavelets corresponding to the same partition of frequency domain as nuDFB and build a framework to analyze the equivalent conditions of PR for critically downsampled as well as more general redundant schemes. In our previous work \cite{yin2014orthshear}, we show that in MRA on a dyadic quincunx lattice, PR is equivalent to an identity condition and a set of shift-cancellation conditions closely related to the frequency support of filters and their downsampling scheme. Based on these two conditions, we rederived Durand's discontinuity result of orthonormal schemes; we also show that a slight relaxation of conditions allows frames with redundancy less than 2 that circumvent the regularity limitation. Furthermore, we have an explicit approach to construct such regularized directional wavelet frames by smoothing the Fourier transform of the irregular directional wavelets.
The main contribution of this paper is that we extend our previous work and show that the same obstruction to regularity as in orthonormal schemes exists in biorthogonal schemes. Different from our previous approach in the orthonormal case, our analysis of bi-orthogonal schemes is inspired by Cohen et al's approach in \cite{cohen1993compactly} for numerical construction of compactly supported symmetric wavelet bases on a hexagonal lattice. We extend and adapt their numerical construction to our bi-orthogonal setting.

The paper is organized as follows. In Section \ref{sec: setup}, we set up the framework of an MRA with dyadic quincunx downsampling. In Section \ref{sec: orth}, we review the regularity analysis of orthonormal schemes and its extension to frames in \cite{yin2014orthshear}. In particular, we derive two conditions, {\it identity summation} and {\it shift cancellation}, equivalent to perfect reconstruction in this MRA with critical downsampling. These lead to the classification of {\it regular/singular} boundaries of the frequency partition %and the corresponding smoothing techniques to improve spatial localization. We compare our and Durand's directional wavelet constructions. 
%In section \ref{sec: frame}, 
and a {\it relaxed shift-cancellation} condition for low-redundancy MRA frame allows better regularity of the directional wavelets. 
In Section \ref{sec: bi-orth}, we extend the orthonormal schemes to biorthogonal schemes as well as the corresponding {\it identity summation} and {\it shift cancellation} conditions. We then introduce Cohen et al's approach in \cite{cohen1993compactly} and adapt it to the regularity analysis on our biorthogonal schemes due to these conditions. We show that the biorthogonal schemes have the same irregularity as in the orthonormal schemes.
In Section \ref{sec: solve-quincunx}, we  propose a numerical algorithm for the construction of biorthogonal schemes along with further analysis on the regularity constraints.
%the resulting linear system does not have feasible solution satisfying all regularity constraints, especially continuity of Fourier transforms of wavelet filters.  
Finally, we present and discuss numerical results of our algorithm in Section \ref{sec: numerics}, and conclude our current work in Section \ref{sec: end}.

%% file: setup.tex
\pdfoutput=1
\section{Framework Setup}\label{sec: setup}
We summarize 2D-MRA systems and the relation between frequency domain partition and sublattice of $\mathbb{Z}^2$ with critical downsampling following \cite{yin2014orthshear}.

\subsection{Notations and conventions}
Throughout this paper, we use upper case bold font for matrices $(e.g.\;\V{A},\V{B})$, lower case bold font for vectors $(e.g.\;\V{a},\V{b})$ and upper case italics for subsets $(e.g.\;C_1,\,C_2)$ of the frequency domain. We denote the conjugate transpose of a matrix $\V{A}$ by $\V{A}^*$. For $\V{a}$ in a $d$-dimensional vector space over $\mathbb{F}$, we use the convention $\V{a}\in\mathbb{F}^{d\times 1}$ and $\V{a}^*$ for its conjugate row vector. 

We adopt conventions in scientific computing programs and packages. For matrices and vectors, the indexing of rows and columns starts with zero. For the axes of the frequency plane, we denote the vertical axis as $\omega_1$-axis with values increasing from top to bottom and the horizontal axis as $\omega_2$-axis with values increasing from left to right, e.g. Figure \ref{fig: partition}.

\subsection{Multi-resolution analysis and sublattice sampling}
In an MRA, given a scaling function $\phi\in L^2(\mathbb{R}^2)$, s.t. $\Vert\phi\Vert_2=1$,
the base approximation space is defined as $V_0 = \overlinespace{span\{\phi_{0,\boldsymbol{k}}\}}_{\boldsymbol{k}\in\mathbb{Z}^2}$, where $\phi_{0,\boldsymbol{k}} = \phi(\boldsymbol{x}-\boldsymbol{k})$. If $\langle \phi_{0,\boldsymbol{k}},\phi_{0,\boldsymbol{k'}}\rangle = \delta_{\boldsymbol{k,k'}}$, then $\{\phi_{0,\boldsymbol{k}}\}$ is an orthonormal basis of $V_0$. In addition, $\phi$ is associated with a scaling matrix $\mathbf{D}\in\mathbb{Z}^{2\times 2}$, s.t. the dilated scaling function
 $\phi_1(\boldsymbol{x}) = |\mathbf{D}|^{-1/2}\phi(\mathbf{D}^{-1}\boldsymbol{x})$ is a linear combination of $\phi_{0,\boldsymbol{k}}$.
Equivalently, $\exists\,m_0(\boldsymbol{\omega}) = m_0(\omega_1,\omega_2)$, $2\pi-$periodic in $\omega_1,\omega_2$, s.t. in the frequency domain
\begin{align}\label{eq: m0}
\widehat{\phi}(\mathbf{D}^T\boldsymbol{\omega}) = m_0(\boldsymbol{\omega})\widehat{\phi}(\boldsymbol{\omega}).
\end{align}
The recursive expression \eqref{eq: m0} of $\widehat{\phi}(\V{\omega})$ implies that
\begin{align}\label{eq: phi-m0}
\textstyle \widehat{\phi}(\boldsymbol{\omega}) = (2\pi)^{-1}\prod_{k=1}^{\infty}m_0(\mathbf{D}^{-k} \boldsymbol{\omega}),
\end{align}
where we have implicitly assumed that $\phi\in L^1(\mathbb{R}^2)$ and $\int\phi\,dx = 1$ (which follows from the other constraints if $\phi$ has some decay at $\infty$).
\\[.2em]
Let $\phi_{l,\boldsymbol{k}} = \phi(\V{D}^{-l}\boldsymbol{x}-\boldsymbol{k})$
and $V_l = \overlinespace{span\{\phi_{l,\boldsymbol{k}};\boldsymbol{k}\in\mathbb{Z}^2\}},\,l\in\mathbb{Z}$ be the nested approximation spaces. Define $W_l$ as the orthogonal complement of $V_l$ with respect to $V_{l-1}$ in MRA. 
Suppose there are $J$ wavelet functions $\psi^j\in L^2(\mathbb{R}^2)$, {\small $1 \leq j \leq J$}, and $\mathbf{Q}\in\mathbb{Z}^{2\times2}$, s.t.
\begin{align*}
\hspace{-1em} W_l = \bigcup_{j=1}^J W_l^j = \bigcup_{j=1}^J \overline{span\{\psi^j_{l,\boldsymbol{k}};\boldsymbol{ k}\in \mathbf{Q}\mathbb{Z}^2\}}= \bigcup_{j=1}^J \overline{span\{\psi^j(\mathbf{D}^{-l}\boldsymbol{x-k});\boldsymbol{ k}\in \mathbf{Q}\mathbb{Z}^2\}},
\end{align*}
 an $L$-level multi-resolution system with base space $V_0$ is then spanned by
 \begin{align}\label{eq: MRA}
 V_L\oplus\,\bigoplus_{l=1}^L\Big(\bigcup_{j=1}^J\, W_l^j\Big) =
 \{\phi_{L,\boldsymbol{k}}\,,\psi^j_{l,\boldsymbol{k'}}\,, \, {\small 1\leq l \leq L,\, \boldsymbol{k}\in \mathbb{Z}^2,\,\boldsymbol{k'}\in \mathbf{Q}\mathbb{Z}^2,\,1\leq j \leq J\}}.
\end{align}  
%In particular, we set $\mathbf{D} = \mathbf{D_2}\doteq\bigl(\begin{smallmatrix} 2&0\\0&2\end{smallmatrix}\bigr)$ and $\mathbf{Q}\doteq\bigl(\begin{smallmatrix} 1&1\\-1&1\end{smallmatrix}\bigr)$.
As $W_1\subset V_0$, each rescaled wavelet $\psi^j(\mathbf{D}^{-1}\cdot)$ is also a linear combination of $\phi_{0,\boldsymbol{k}}$, so that $\exists\, m_j$ analogous to $m_0$
satisfying 
\begin{align}\label{eq: mj}
\widehat{\psi}^j(\mathbf{D}^T\boldsymbol{\omega}) = m_j(\boldsymbol{\omega})\widehat{\phi}(\boldsymbol{\omega}),\hspace{1cm} 1\leq j \leq J.
\end{align}
%In this specific construction of MRA,
% the scaling function $\phi$ and all the wavelet functions $\psi^j$ share the same scaling matrix $\mathbf{D}$, yet the family of shifted $\phi_{\boldsymbol{k}}$ is defined on $\mathbb{Z}^2$, whereas the family of shifted $\psi^j_{\boldsymbol{k}}$ is defined on a sub-integer lattice $\mathbf{Q}\mathbb{Z}^2$. Hence 
% the corresponding subsampling matrix of $\phi_{1,\boldsymbol{k}}$ is $\mathbf{D}$ and that of $\psi^j_{1,\boldsymbol{k}}$ is $\mathbf{QD}$, the dyadic quincunx subsample (see the right panel in Fig.\ref{fig: partition}), as in \cite{durand2007}. 
% We haven't yet imposed any condition on this MRA, or equivalently, on $m-$functions and the subsampling matrices $\mathbf{D}$ and $\mathbf{Q}$; this comes next.

\subsection{Frequency domain partition and critical downsampling}\label{subsec: frequency partition}

%critical downsampling thus depends only on the subsampling matrices $\mathbf{D}$ and $\mathbf{Q}$. The space decomposition structure $V_0 = V_1\bigoplus W_1$ in MRA and \eqref{eq: m0}, \eqref{eq: mj} require consistency between the $m-$functions and the subsampling matrices $\mathbf{D}$ and $\mathbf{Q}$. 

\begin{figure}[!t]
\centering
\begin{minipage}[c]{.35\textwidth}
\includegraphics[width=\textwidth]{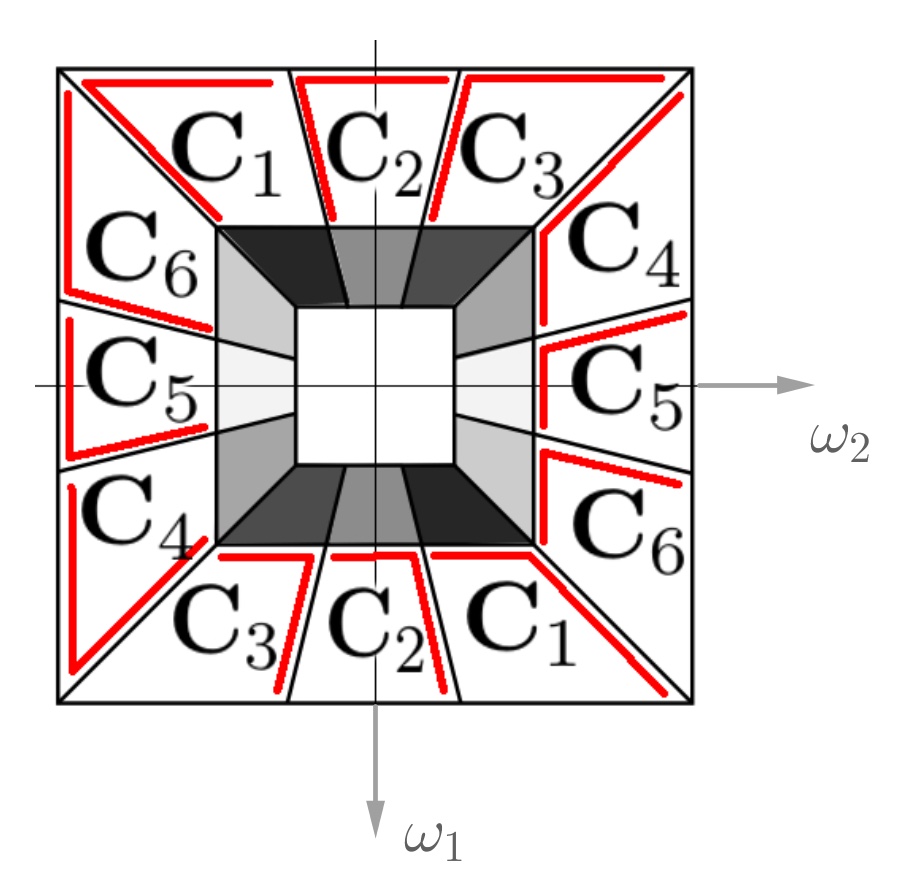}
\end{minipage}\hspace*{3em}
\begin{minipage}[c]{.25\textwidth}
\includegraphics[width=\textwidth]{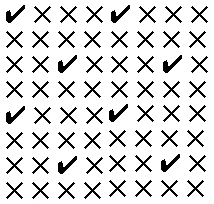}
\vspace*{1em}
\end{minipage}
\caption{Left: partition of $S_0$ and boundary assignment of $C_j$, $j = 1,\cdots,6$ (each $C_j$ has boundaries indicated by red line segments), Right: dyadic quincunx sublattice. Note that the $\omega_1$-axis is vertical and the $\omega_2$-axis is horizontal by our convention.}
\label{fig: partition}
\vspace*{-5mm}
\end{figure}

%\begin{mydef}
%If $\mathcal{L}$ is the lattice generated by $\boldsymbol{a}_1,\boldsymbol{a}_2$, i.e. $\mathcal{L} = \sum_{i=1,2}k_i\boldsymbol{a}_i,\,k_i\in\mathbb{Z}$,the {\bf reciprocal lattice} $\mathcal{L}^*$ of $\mathcal{L}$ is the lattice generated by the vectors $\boldsymbol{b}_1,\boldsymbol{b}_2$, s.t. $\boldsymbol{b}_i^T\boldsymbol{a}_j = 2\pi\delta_{ij}$. 
%\end{mydef}

%\begin{mydef}
%Given a lattice $\mathcal{L}$, a {\bf fundamental domain} $S$ in $\mathbb{R}^2$ with respect to $\mathcal{L}$, denoted as $S = \mathbb{R}^2/\mathcal{L}$, is a subset of $\mathbb{R}^2$, such that $\bigcup_{l\in\mathcal{L}}(S+l) = \mathbb{R}^2$ and $S\cap(S+l)=\varnothing,\,\forall l\in\mathcal{L}\setminus\{\mathbf{0}\}$.
%A set $S$ is a {\bf frequency support} of $\mathcal{L}$ if $S = \mathbb{R}^2/\mathcal{L}^*$.
%\end{mydef}
%Furthermore, each sub-lattice $\mathcal{L}$ of $\mathbb{Z}^2$ is associated with a set of shifts $\Gamma, \, s.t. \bigcup_{\V{\gamma}\in\Gamma}\mathcal{L}+\V{\gamma} = \mathbb{Z}^2$ and $|\Gamma| = |\mathbb{Z}^2/\mathcal{L}|$.

Consider the canonical frequency square, $S_0 = [-\pi,\pi)\times[-\pi,\pi)$ associated with the lattice $\mathcal{L} = \mathbb{Z}^2$. 
%Due to \eqref{eq: m0} and \eqref{eq: mj}, this is equivalent to $S_0=\bigcup_{0\leq j\leq J} supp(m_j\vert_{S_0})$.That is, 
For $L=1$, the 1-level decomposition \eqref{eq: MRA} together with \eqref{eq: m0} and \eqref{eq: mj} implies that the union of the support of $m_j,\,0\leq j\leq J$ covers $S_0$.
Furthermore, $\exists\, C_j\subset supp(m_j), 0\leq j\leq J,$ such that they form a partition of $S_0$; conversely, given a partition $C_j$ of $S_0$, we may construct an MRA where $m_j$ are ``mainly" supported on $C_j$ (this will become more explicit in Section \ref{subsec: discontinuity}).
%, and it is natural to take $C_j$ as the main support of $m_j$.
To build an orthonormal basis with good directional selectivity, we choose the partition of $S_0$ shown in the left of Figure \ref{fig: partition}, which is the same for Example B in \cite{durand2007} and the least redundant shearlet system \cite{kutyniok2012digital}. In this partition, $S_0$ is divided into a central square $C_0 = \bigl(\begin{smallmatrix} 2&0\\0&2\end{smallmatrix}\bigr)^{-1}S_0$ and a ring: the ring is further cut into six pairs of directional trapezoids $C_j$ by lines passing through the origin with slopes $\pm 1, \pm 3$ and $\pm \frac{1}{3}$. The central square $C_0$ can be further partitioned in the same way to obtain a two-level multi-resolution system, as shown in Figure \ref{fig: partition}.

%To build our first example, in which $\hat{\phi},\,\widehat{\psi}^j$ are indicator functions in $\mathbb{R}^2$, we consider the case where  and we pick $S_0=\mathbb{R}^2/(\mathbb{Z}^2)^*$, . 
%Since $\phi_1,\psi^j_1$ and their shifts span the space $V_0$, $supp(\widehat{\phi}_1)$ and $supp(\widehat{\psi}^j)$, together, should thus cover $S_0$. Due to \eqref{eq: m0} and \eqref{eq: mj}, this is equivalent to $S_0=\bigcup_{0\leq j\leq J} supp(m_j\vert_{S_0})$. That is, if $C_j,\,{\small 0\leq j\leq J}$ are the main support of $m_j,\,0\leq j\leq J$ respectively, then they form a partition of $S_0$. An non-uniform admissible partition is defined as follows,

%\begin{mydef}
%$C_j, 0\leq j\leq J$ is an {\bf admissible} partition of $S_0$ if and only if $\exists \mathbf{D}, \mathbf{Q}\in\mathbb{Z}^{2\times 2}$, s.t. the low frequency piece $C_0 = \mathbb{R}^2/(\mathbf{D}\mathbb{Z}^2)^*,$ and the high frequency pieces $C_j = \mathbb{R}^2/(\mathbf{QD}\mathbb{Z}^2)^*,\,j = 1,\cdots,J$.
%\end{mydef}
In the corresponding MRA generated by \eqref{eq: MRA}, $J=6$ and $\mathbf{D} =\bigl(\begin{smallmatrix} 2&0\\0&2\end{smallmatrix}\bigr)$, and we choose  $\mathbf{Q}$ specifically to be $\bigl(\begin{smallmatrix} 1&1\\-1&1\end{smallmatrix}\bigr)$. Because $ |\mathbf{D}|^{-1} + J|\mathbf{QD}|^{-1} = 1/4 + 6/ (2\cdot 4) = 1$, the corresponding MRA generated by \eqref{eq: MRA} achieves critical downsampling(\cite{durand2007}). The scaling matrix of $\psi^j$ is $\V{QD}=\bigl(\begin{smallmatrix} 2 & 2\\-2 & 2\end{smallmatrix}\bigr)$, which corresponds to downsampling on the dyadic quincunx sublattice $\V{QD}\mathbb{Z}^2$ (see the right panel in Figure \ref{fig: partition}), as in \cite{durand2007}. 

%\noindent{\it Remark.}
%The dyadic quincunx subsampling $\V{QD}$ is not ideal for $\psi^2$ with frequency support on  $C_2$ that characterizes horizontal features. However, we can adapt the subsampling lattice $\V{QD}\mathbb{Z}^2$ by shearing it. In particular, let $\V{S} = \bigl(\begin{smallmatrix} 1&0\\1&1\end{smallmatrix}\bigr)$ be the shearing matrix, then $\V{SQD}\mathbb{Z}^2 = \bigl(\begin{smallmatrix}2 & 2\\ 0 & 4\end{smallmatrix}\bigr) \mathbb{Z}^2= \bigl(\begin{smallmatrix}2 & 0\\0 & 4\end{smallmatrix}\bigr)\mathbb{Z}^2$, a sublattice where the horizontal subsampling factor is twice of the vertical one. Since $|\V{S}|=1$, the critical downsampling property still holds and all the analysis in the rest of this paper stays the same. For $\psi^5$, we can shear $\V{QD}\mathbb{Z}^2$ by $\V{S}' = \bigl(\begin{smallmatrix} 1 & -1\\ 0 & 1\end{smallmatrix}\bigr)$ to obtain $\V{S}'\V{QD}\mathbb{Z}^2 = \bigl(\begin{smallmatrix} 4 & 0\\0 & 2\end{smallmatrix}\bigr)\mathbb{Z}^2$.

This downsampling scheme is compatible with $C_j$. Consider two sets of shifts in the frequency domain $\Gamma_0= \{\V{\pi}_i,\, i = 0,2,4,6\}$ and $\Gamma_1 = \{\V{\pi}_i,\, i = 0,1,\cdots,7\}$, where {\small $\V{\pi}_0 = (0,0), \V{\pi}_1 = (\pi/2,\pi/2), \V{\pi}_2 = (\pi,0),\V{\pi}_3 = (-\pi/2,\pi/2), \V{\pi}_4 = (0,\pi), \V{\pi}_5 = (\pi/2,-\pi/2),\V{\pi}_6 = (\pi,\pi), \V{\pi}_7=(-\pi/2,-\pi/2)$}. $\Gamma_0$ and $\Gamma_1$ characterize the sublattices $\V{D}\mathbb{Z}^2$ and $\V{QD}\mathbb{Z}^2$ respectively by
$\sum_{\V{\pi}\in\Gamma_0}e^{i\V{\alpha}^\top\V{\pi}} = |\Gamma_0| \, \mathbbm{1}_{\V{D}\mathbb{Z}^2}(\V{\alpha})$ and $ \sum_{\V{\pi}\in\Gamma_1}e^{i\V{\alpha}^\top\V{\pi}} = |\Gamma_1| \, \mathbbm{1}_{\V{QD}\mathbb{Z}^2}(\V{\alpha})$, where $\mathbbm{1}$ is the indicator function. 
We observe that each $C_j$ forms a tiling of $S_0$ under the shifts associated with the sublattice where the coefficients of $\psi^j$ are downsampled:
% Each $C_j$ together with its shifts form a tiling of $S_0$, i.e.
\begin{align}\label{eq: tiling}
S_0 = \bigcup_{\V{\pi}\in\Gamma_1}\left(C_j + \V{\pi}\right) = \bigcup_{\V{\pi}\in\Gamma_0}\left(C_0 + \V{\pi}\right),\quad j = 1,\cdots,6.
%\mathbb{R}^2/(\mathbb{Z}^2)^*=\bigcup_{\V{\pi}\in\Gamma_0}\left(\mathbb{R}^2/(\V{D_2}\mathbb{Z}^2)^* %+ \V{\pi}\right)
%=\bigcup_{\V{\pi}'\in\Gamma_1}\left(\mathbb{R}^2/(\V{QD_2}\mathbb{Z}^2)^* + \V{\pi}'\right),
\end{align}
%where%then $\mathbf{D_2}\mathbb{Z}^2$ is associated with the set of shifts 
%$\Gamma_0= \{\V{\pi}_i,\, i = 0,2,4,6\}$ and %$\mathbf{QD_2}\mathbb{Z}^2$ is associated with 
%$\Gamma_1 = \{\V{\pi}_i,\, i = 0,1,\cdots,7\}$.
Alternatively, we say that $\{\,C_j,\, j = 0,\cdots,6\,\}$ is an {\it admissible} partition of $S_0$ with respect to the dyadic quincunx downsampling scheme.
%This partition is admissible and corresponds to $\mathbf{D} = \mathbf{D_2}\doteq\bigl(\begin{smallmatrix} 2&0\\0&2\end{smallmatrix}\bigr)$ and $\mathbf{Q}:=\bigl(\begin{smallmatrix} 1&1\\-1&1\end{smallmatrix}\bigr)$. The wavelet coefficients are taken on the dilated quincunx sub-lattice $\mathbf{QD_2}\mathbb{Z}^2$ (see the right panel in Fig.\ref{fig: partition}).
%In addition, $|\mathbf{D_2}| = 4,\,|\mathbf{Q}| = 2$ so that $1/4 + 6/ (2\cdot 4) = 1$, and the system is critical downsampling.
The admissible property guarantees the existence of orthonormal bases consisting of directional filters on the dyadic quincunx sublattice with frequency support in $C_j$.

%% file: orthonormal.tex
\pdfoutput=1
\section{Orthonormal Bases}\label{sec: orth}

In this section, we discuss the conditions on $m_j$ such that the corresponding MRA forms an orthonormal bases. 
%consider orthonormal bases with $m-$functions defined in \eqref{eq: m0} and \eqref{eq: mj} whose supports mainly corresponding to the partition of $S_0$ in Fig.\ref{fig: partition}.
%we always consider a multi-resolution system with scaling function $\phi$ and quasi-shearlets $\psi^j$, {\small$ j = 1,\dots,6$} defined by $(M_0, D_2)$ and $(M_j, Q)$, {\small$j = 1,\dots,6$} respectively. Furthermore, the essential support of $M_j$'s corresponds to the partition of $\mathbf{S}_0$ in a shearlet system.

%The construction of \eqref{eq: MRA} reduces to design $m_0$ in \eqref{eq: m0} and $m_j, j= 1,\cdots,6$ in \eqref{eq: mj}.
We begin with the two key conditions, i.e. {\it identity summation} and {\it shift cancellation}, on $m_j$ such that the system \eqref{eq: MRA} is perfect-reconstruction (PR) or equivalently a Parseval frame in MRA.% weaker than the orthonormal condition.

%\subsection{Identity summation and shift cancellation}
\subsection{orthonormal conditions on $m_j$}\label{subsec: northonormal cond}
In MRA, \eqref{eq: MRA} is PR if $\forall f\in L_2(\mathbb{R}^2)$,
\begin{equation}
%\textstyle 
\sum_{\V{k}\in\mathbb{Z}^2}\langle f,\phi_{0,\V{k}}\rangle\phi_{0,\V{k}} = \sum_{\V{k}\in\mathbb{Z}^2}\langle f,\phi_{1,\V{k}}\rangle\phi_{1,\V{k}} + \sum_{j=1}^J\sum_{\V{k}'\in\V{Q}\mathbb{Z}^2}\langle f,\psi^{j}_{1,\V{k}'}\rangle\psi^j_{1,\V{k}'}.\label{eq: PR}
\end{equation}
Using \eqref{eq: m0} and \eqref{eq: mj} together with the admissibility of the frequency partition \eqref{eq: tiling}, condition \eqref{eq: PR} on $\phi$ and $\psi^j$ yields:
\begin{thm}\label{thm: conds}
Let $J=6,\, \V{D} =\bigl(\begin{smallmatrix} 2&0\\0&2\end{smallmatrix}\bigr)$ and $\V{Q}=\bigl(\begin{smallmatrix} 1&1\\-1&1\end{smallmatrix}\bigr)$ in \eqref{eq: MRA}.
Then the perfect reconstruction condition holds for \eqref{eq: MRA} if and only if the following two conditions hold.
\begin{align}\label{eq: id-sum}
|m_0(\boldsymbol{\omega})|^2 + \sum_{j = 1}^6|m_j(\boldsymbol{\omega})|^2 = 1.
\end{align}
\begin{equation}\label{eq: shift-cancel}
 \begin{cases}
%M_0(\boldsymbol{omega})\overline{M_0(\boldsymbol{omega}+\boldsymbol{\gamma})} + 
\sum_{j = 0}^6m_j(\boldsymbol{\omega})\overlinespace{m_j}(\boldsymbol{\omega} + \boldsymbol{\pi}) = 0, & \boldsymbol{\pi}\in \Gamma_0\setminus\{\boldsymbol{0}\}.\\[.5em]
\sum_{j=1}^6m_j(\boldsymbol{\omega})\overlinespace{m_j}(\boldsymbol{\omega}+\boldsymbol{\pi}) = 0, & \boldsymbol{\pi}\in\Gamma_1\setminus\Gamma_0.
\end{cases}
\end{equation}
%where  $\Lambda = (QD\mathbb{Z}^2)^*/(\mathbb{Z}^2)^*,\,\Gamma = (D\mathbb{Z}^2)^*/(\mathbb{Z}^2)^*.$ %$\{(\frac{\pi}{2},\frac{\pi}{2}),(\frac{3\pi}{2},\frac{\pi}{2}),(\frac{\pi}{2},\frac{3\pi}{2}),(\frac{3\pi}{2},\frac{3\pi}{2}),$ $(0,0),(0,\pi),(\pi,0),(\pi,\pi)\}$
\end{thm}

Theorem \ref{thm: conds} is a corollary of Proposition 1 and Proposition 2 in \cite{durand2007}. We give an alternate proof in Appendix \ref{app: cond-thm}.
In Theorem \ref{thm: conds}, \eqref{eq: id-sum} is the {\it identity summation} condition, guaranteeing conservation of $l_2$ energy; \eqref{eq: shift-cancel} is the {\it shift cancellation} condition such that aliasing is canceled correctly in reconstruction from wavelet coefficients. %, such that downsampling of scaling and shearlet coefficients is valid. 
 %$\Lambda = \{(\frac{\pi}{2},\frac{\pi}{2}),(\frac{3\pi}{2},\frac{\pi}{2}),(\frac{\pi}{2},\frac{3\pi}{2}),(\frac{3\pi}{2},\frac{3\pi}{2}),$ $(0,0),(0,\pi),(\pi,0),(\pi,\pi)\},\Gamma = \{(0,0),(0,\pi),(\pi,0),(\pi,\pi)\}$ are the sets of shifts associated to quincunx-dyadic dilation $QD$ and dyadic dilation $D$ respectively.
%Each $M_j$ contributes a term $M_j(\boldsymbol{omega})\overline{M_j(\boldsymbol{omega} + \boldsymbol{\nu})}$ in the cancellation condition of any shift $\boldsymbol{\nu}$ corresponding to the downsampling scheme of $M_j$.
Because each $m_j$ is $(2\pi,2\pi)$ periodic, we only need to check these conditions $\forall \V{\omega}\in S_0$.

%\subsection{Extra condition for basis}
%By Theorem \ref{thm: conds}, the system \eqref{eq: MRA} is a Parseval frame ; 
Moreover, for  \eqref{eq: MRA} to be an orthonormal basis,  $\{\phi_{\boldsymbol{k}}\}_{\boldsymbol{k}\in\mathbb{Z}^2}$ need to be an orthonormal basis, which is determined by $m_0$ in \eqref{eq: phi-m0}. In 1D MRA, Cohen's theorem in \cite{cohen1992biorthogonal} provides a necessary and sufficient condition on $m_0$ such that \eqref{eq: MRA} is an orthonormal basis. %Such a condition in 1D wavelet MRA is given by Cohen \cite{cohen1992biorthogonal}. 
This theorem generalizes to 2D in  e.g. \cite{yin2014orthshear}, as follows.
\begin{thm}\label{thm: basis cond}
Assume that $m_0$ is a trigonometric polynomial with $m_0(\V{0})=1$, and define $\hat{\phi}(\boldsymbol{\omega})$ as in \eqref{eq: phi-m0}.\\
If $\phi(\cdot - \boldsymbol{k}),\boldsymbol{k}\in\mathbb{Z}^2$ are orthonormal, then $\exists K$ containing a neighborhood of 0, s.t. $\forall\boldsymbol{\omega}\in S_0,\,\boldsymbol{\omega}+2\pi\mathbf{n}\in K$ for some $\mathbf{n}\in\mathbb{Z}^2, $ and $\inf_{k>0,\,\boldsymbol{\omega}\in K}|m_0(\mathbf{D_2}^{-k}\boldsymbol{\omega})| >0$. 
 Further, if $\sum_{\boldsymbol{\V{\pi}}\in \Gamma_0} |m_0(\boldsymbol{\omega}+\boldsymbol{\pi})|^2 = 1$, then the inverse is true.
\end{thm}
%Theorem \ref{thm: basis cond} can be proved similarly to Cohen's theorem (\cite{cohen1992biorthogonal}).
%Because it is difficult to directly design $M_0$ that satisfies the conditions in Theorem \ref{thm: basis cond}, 
%Below, we construct $m-$functions imposing only \eqref{eq: id-sum} and \eqref{eq: shift-cancel} and then check if the resulting Parseval frame is an orthonormal basis by applying Theorem \ref{thm: basis cond} to $m_0$.

\subsection{Regularity of $m_j$ supported on the $C_j$}\label{sec: design}
%In this section, we define Shannon-type directional orthonormal basis same as in \cite{durand2007} and \cite{nuDFB05}. Then, we apply direct smoothing to its $m-$functions to improve spatial localization, this leads to a critical analysis of the boundary regularity of $S_1$ and the $C_j$'s.

%\subsection{Shannon-type wavelets and smoothing}
In this subsection, we consider $m_j$ supported on the $C_j$ introduced in Section \ref{subsec: frequency partition} that satisfy orthonormal conditions in Section \ref{subsec: northonormal cond}. We begin with the Shannon-type wavelet construction,
where $m_j$ are indicator functions $ m_j = \mathbbm{1}_{C_j},\,0\leq j \leq 6,\,$ and we use the boundary assignment of $C_j$ in Figure \ref{fig: partition}.
The identity summation follows from the partition of $S_0$ by the $C_j$, and the shift cancellation follows from the admissible property \eqref{eq: tiling}. % which implies $m_j(\boldsymbol{\omega})\overline{m_j(\boldsymbol{\omega} + \boldsymbol{\pi}_i)}\equiv 0,\,\forall j,\, i\neq 0.$ %Let $\partial\mathbf{C}_j = \overline{\mathbf{C}_j}-\mathbf{C}_j^{\circ}$ be the boundary of $\mathbf{C}_j$,  
Applying Theorem \ref{thm: basis cond} to $m_0$, we verify that the Shannon-type wavelets generated from these $m_j$ form an orthonormal basis.

\begin{figure}[!t]
\centering
%\hspace*{-5mm}
%\begin{minipage}{.6\textwidth}
%\includegraphics[width=\textwidth]{fig_overlap.png}\hspace*{-2mm}
%\end{minipage}
\begin{minipage}{.25\textwidth}
\includegraphics[width=\textwidth]{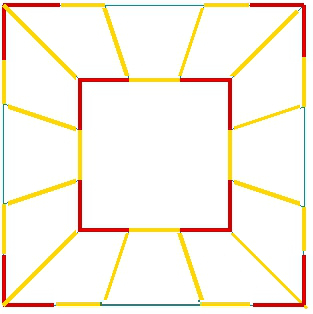}
\end{minipage}
\caption[caption]{
%\textcolor{red}{
%{\it Left top}: the supports of $m_j$ (green) and $m_j(\cdot+\V{\pi}_2$ (red) for $j = 5,6$ after smoothing, overlap on the vertical boundary at $\omega_1 =  \pm \pi/2$ of $C_5$ (green) and its shift (red) by $\boldsymbol{\pi}_2 = (\pi,0)$. Note that two copies of shifted $S_0$ (red) overlap the un-shifted $S_0$(green) due to the $(2\pi,2\pi)$ periodicity of $m_j$. Only $m_0$ and $m_5$ have overlapping smoothed boundaries by $\boldsymbol{\pi}_2$. 
%\\\hspace{\textwidth}
%{\it Left bottom}: intersection of $\mathcal{B}(0,\boldsymbol{\pi}_2)$ and $\mathcal{B}(5,\boldsymbol{\pi}_2)$ in yellow and $\mathcal{C}(0,\boldsymbol{\pi}_2) = \mathcal{B}(0,\boldsymbol{\pi}_2)\setminus\mathcal{B}(5,\boldsymbol{\pi}_2)$ in red. Smoothing $m_0$ in the red (singular) region is impossible without violating \eqref{eq: shift-cancel}. %This leads to the distinction between regular (yellow) and singular (red) boundaries at $omega_1 = \pm\pi/2$.
%\\\hspace{\textwidth}
%{\it Right}: 
Boundary classification, singular (red) and regular (yellow) }%after similar arguments for all shifts $\boldsymbol{\pi}_i$.}
 %}
\label{fig: boundary}
\end{figure}

%However, because of the discontinuity of $m_j$ across the boundary of its support, the corresponding wavelet has slow decay in the time domain. 
%Such $m_j$'s are not regularized due to 
Because of the discontinuity at $\partial C_j$, the boundaries of the $C_j$, these $m_j$ are not smooth, and hence the corresponding wavelets are not spatially localized. The $m_j$ can be regularized by smoothing at the $\partial C_j$.
%We take a different regularization approach from Durand's \cite{durand2007}, where three regular quincunx filter banks are constructed and then composed to obtain the desired regular quincunx dyadic filter banks. Here, we smooth the discontinuous boundaries of $m-$functions directly. 
However, as shown in Proposition 3 in \cite{durand2007}, it is not possible to smooth the behavior of the $m_j$ at  {\it all} the boundaries with discontinuity if the $m_j$ have to satisfy the perfect reconstruction condition.
In \cite{yin2014orthshear}, the $\partial C_j$ are segmented into {\it singular} and {\it regular} pieces with respect to the shift cancellation condition \eqref{eq: shift-cancel} in Theorem \ref{thm: conds}. On regular boundaries, pairs of $(m_j,\, m_{j'})$ share a boundary and can both be smoothed in a coherent way such that both \eqref{eq: id-sum} and \eqref{eq: shift-cancel} remain satisfied. The singular pieces are boundaries for just one $m_j$, which can then not be smoothed without violating the shift cancellation condition. Figure \ref{fig: boundary} shows the boundary classification, where the corners of $S_0$ and $C_0$ are singular, hence $m_0$ and the $m_j$'s in two diagonal directions of an orthonormal bases are discontinuous there. A mechanism of constructing orthonormal bases by smoothing Shannon-type $m_j$ on regular boundaries is provided in \cite{yin2014orthshear}.%However, it remains unclear if some of the discontinuity can be removed by direct smoothing.

%Next, we analyze the limitation of direct smoothing in detail and show that there are regular boundaries which can be smoothed without violating \eqref{eq: id-sum} and \eqref{eq: shift-cancel}, and singular boundaries which cannot be smoothed.

%\textcolor{red}{In the remainder of this paper, we explore how this can be done for the first dyadic layer (without cutting further at higher frequencies); we call the corresponding basis {\it quasi-Shearlets}.}

%% file: frame.tex
%\section{Low-redundancy frame construction}\label{sec: frame}
\subsection{Extension to low-redundancy tight frame}\label{sec: frame}
%Consider the $L-$level directional wavelet MRA system
The irregularity of orthonormal bases can be overcome in the following low-redundancy tight frame construction,
 \begin{align}\label{eq: MRA-frame}
 \{\phi_{L,\boldsymbol{k}}\,,\psi^j_{l,\boldsymbol{k'}}\,, \, {\small 1\leq l \leq L,\, \boldsymbol{k},\,\boldsymbol{k'}\in \mathbb{Z}^2,\,1\leq j \leq 6\}}.
\end{align}  
In \eqref{eq: MRA-frame}, all wavelet coefficients are downsampled on the dyadic sublattice and the redundancy of any such $L-$level frame does not exceed $\frac{J/|D|}{1-1/|D|} = \frac{6/4}{1-1/4} = 2$.
%where $\phi,\psi^j$ satisfy \eqref{eq: m0} and \eqref{eq: mj} as before. Instead of taking the dilated quincunx subsampling of directional wavelet coefficients of  \eqref{eq: MRA}, a dyadic subsampling is taken instead. A 1-level MRA frame \eqref{eq: MRA-frame} has redundancy $\frac{1}{|D|} + \frac{J}{|D|} = 1/4 + 6/4 = 7/4$, and the redundancy for any $L-$level MRA frame doesn't exceed $\frac{J/|D|}{1-1/|D|} = \frac{6/4}{1-1/4} = 2$. 
Similar to Theorem \ref{thm: conds}, we have
\begin{thm}\label{thm: frame-conds}
%Set $\Gamma = (D\mathbb{Z}^2)^*/(\mathbb{Z}^2)^*.$ 
\eqref{eq: MRA-frame} has PR if and only if the following two conditions both hold.
\begin{align}
\textstyle |m_0(\boldsymbol{\omega})|^2 + \sum_{j = 1}^6|m_j(\boldsymbol{\omega})|^2 &= 1. \\
\textstyle\sum_{j = 0}^6\,m_j(\boldsymbol{\omega})\overlinespace{m_j}(\boldsymbol{\omega} + \boldsymbol{\pi}) &= 0,\quad  \boldsymbol{\pi}\in \Gamma_0\setminus\{\boldsymbol{0}\}. \label{eq: reduced-shift-cancel}
\end{align}
\end{thm}
Theorem \ref{thm: frame-conds} can be proved analogously to Theorem \ref{thm: conds}, 
but with fewer shift cancellation constraints. Following the same analysis of boundary regularity as before, we show in \cite{yin2014orthshear} that all boundaries are regular with respect to \eqref{eq: reduced-shift-cancel} and can be smoothed properly. Hence, we were able to obtain directional wavelets with much better spatial and frequency localization than those constructed by Durand in \cite{durand2007}. 
%We can define {\it singular} boundaries as before, %and Lemma \ref{lem: singular-bdy} holds for \eqref{eq: reduced-shift-cancel} as well.
%but only $\{\mathcal{B}(j,\boldsymbol{\pi})\}_{\boldsymbol{\pi}\in\Gamma_0\setminus\{\boldsymbol{0}\}}$ need to be considered, which results in fewer singular boundaries $\{\mathcal{C}_s(j,\boldsymbol{\pi})\}_{\boldsymbol{\pi}\in\Gamma\setminus\{\boldsymbol{0}\}}$; 
% In particular, we construct a directional wavelet tight frame with redundancy of 2 by using the classical dyadic downsampling $D_2$, with shift cancellaiton constraints \eqref{eq: shift-cancel} only on set $\Gamma\setminus\{\boldsymbol{0}\}$. 
%We check that within these singular boundaries, 
%and no "double" singular boundaries now.

%This means that even though $supp(m_0)$ still cannot be extended outside of the four corners of $S_1$ due to $\mathcal{C}_s(0,(\pi,0))$ and $\mathcal{C}_s(0,(0,\pi))$, $m_1$ can penetrate into the inside of $S_1$ because $\mathcal{C}_s(1,(\pi/2,3\pi/2))$ is not a singular boundary in \eqref{eq: MRA-frame}. The same is true for $m_3,m_4$ and $m_6$. This makes smoothing the boundaries of $m_0$ inwards possible without violating \eqref{eq: id-sum}, see Fig. \ref{fig: many-squares}(c). At the price of double redundancy, we obtain directional wavelets with much better spatial localization; see Fig. \ref{fig: many-squares}(d)(e):
%the discontinuities of a directional wavelets basis in the frequency domain around the singular boundaries can be removed in a low redundant directional wavelet tight frame.

So far, we have considered two directional wavelet MRA systems \eqref{eq: MRA} and \eqref{eq: MRA-frame} such that the directional wavelets characterize 2D signals in six equi-angled directions. 
%The orthonormal basis we construct has better frequency localization than the one constructed by Durand in \cite{durand2007} ( see Fig. \ref{fig: design} and \ref{fig: many-squares}(b)(c)), but has long tails in certain spatial directions, unavoidable because of "double" singular boundaries. 
%By doubling the redundancy we obtain spatially well localized directional wavelets.
Furthermore, these wavelets are well localized in the frequency domain such that $supp(m_j)$ is convex and $\exists\,\epsilon\; s.t.$
\begin{align}\label{eq: no-alians}
 \sup_{\boldsymbol{\omega}'\in supp(m_j)}\inf_{\boldsymbol{\omega}\in C_j}\Vert\boldsymbol{\omega'} - \boldsymbol{\omega}\Vert < \epsilon,\quad  0\leq j\leq 6.
\end{align}
This desirable condition is hard to obtain by multi-directional filter bank assembly of several elementary filter banks.

In the next section, we analyze the more general case of directional bi-orthorgonal filters constructed with respect to the same frequency partition. 

%% file: bi-orth.tex
\section{Biorthogonal Bases}\label{sec: bi-orth}
In this section, we analyze biorthogonal bases in the following form of MRA,
\begin{align}\label{eq: bi-orth MRA}
\{\phi_{L,\V{k}},\widetilde{\phi}_{L,\V{k}}, \psi_{l,\V{k}'}^j,\widetilde{\psi}_{l,\V{k}'}^j,\, 1\leq l\leq L,\,\V{k}\in\mathbb{Z}^2,\, \V{k}'\in\mathbf{Q}\mathbb{Z}^2,\,1\leq j\leq 6 \},
\end{align}
where $\phi$ and $\psi^j$ satisfy \eqref{eq: m0} and \eqref{eq: mj} respectively, and likewise for $\widetilde{\phi}$ and $\widetilde{\psi^j}$,
\begin{align}\label{eq: mj_dual}
\widehat{\widetilde{\phi}}(\V{D}^T\V{\omega}) = \widetilde{m_0}(\V{\omega})\widehat{\widetilde{\phi}}(\V{\omega}),\quad \widehat{\widetilde{\psi^j}}(\V{D}^T\V{\omega}) = \widetilde{m_j}(\V{\omega})\widehat{\widetilde{\phi}}(\V{\omega}).
\end{align}
For such biorthogonal bases, we have the similar identity summation and shift cancellation conditions to those in Theorem \ref{thm: conds}.
\begin{thm}\label{thm: bi-orth conds}
\eqref{eq: bi-orth MRA} has PR if and only if the following two conditions hold
\begin{align}\label{eq: id-sum 2}
m_0(\boldsymbol{\omega})\sbarm{0} + \sum_{j = 1}^6 m_j(\boldsymbol{\omega})\sbarm{j} = 1,
\end{align}
\begin{equation}\label{eq: shift-cancel 2}
\begin{cases}
\sum_{j = 0}^6m_j(\boldsymbol{\omega})\overlinespace{\widetilde{m_j}}(\boldsymbol{\omega} + \boldsymbol{\pi}) = 0, & \boldsymbol{\pi}\in \Gamma_0\setminus\{\boldsymbol{0}\}.\\[.5em]
\sum_{j=1}^6m_j(\boldsymbol{\omega})\overlinespace{\widetilde{m_j}}(\boldsymbol{\omega}+\boldsymbol{\pi}) = 0, & \boldsymbol{\pi}\in\Gamma_1\setminus\Gamma_0.
\end{cases}
\end{equation}
\end{thm}
\vspace{1em}
%where $\M\in\mathbb{C}^{8\times 7}$ and $\mathbf{m}_0\in\mathbb{C}^7$.
We also have the following analogue of Theorem \ref{thm: basis cond}.
\begin{thm}\label{thm: basis cond 2}
Assume that $m_0, \widetilde{m_0}$ are trigonometric polynomials with $m_0(\V{0})=\widetilde{m_0}(\V{0}) = 1$, which generate $\phi,\widetilde{\phi}$ respectively.\\
If $\phi(\cdot - \boldsymbol{k}),\widetilde{\phi}(\cdot - \boldsymbol{k}),\,\boldsymbol{k}\in\mathbb{Z}^2$ are biorthogonal, then $\exists K$ containing a neighborhood of 0, s.t. $\forall\boldsymbol{\omega}\in S_0,\,\boldsymbol{\omega}+2\pi\mathbf{n}\in K$ for some $\mathbf{n}\in\mathbb{Z}^2, $ and $\inf_{k>0,\,\boldsymbol{\omega}\in K}|m_0(\mathbf{D_2}^{-k}\boldsymbol{\omega})| >0$, $\inf_{k>0,\,\boldsymbol{\omega}\in K}|\widetilde{m_0}(\mathbf{D_2}^{-k}\boldsymbol{\omega})| >0$. 
 Furthermore, if  $\sum_{\boldsymbol{\V{\pi}}\in \Gamma_0} m_0(\boldsymbol{\omega}+\boldsymbol{\pi})\sbarmp{0}{} = 1,$ then the inverse is true.
\end{thm}
By Theorem \ref{thm: basis cond 2}, $m_0$ and $\widetilde{m_0}$ need to satisfy the following identity constraint for the MRA \eqref{eq: bi-orth MRA} to be biorthogonal,
\begin{align}\label{eq: identity-cond}
m_0\sbarm{0} + m_0\sbarmp{0}{2} + m_0\sbarmp{0}{4} + m_0\sbarmp{0}{6} = 1.
\end{align}
Furthermore, the identity summation and shift cancellation conditions \eqref{eq: id-sum 2} and \eqref{eq: shift-cancel 2} from Theorem \ref{thm: bi-orth conds}
 can be combined into a linear system with respect to $m_j$ as follows,
\begin{align}\label{eq: LS-new}
%\overline{\M}(\V{\omega})\mathbf{m}_0(\V{\omega})=
\begin{bmatrix}
    \,\sbarm{0} & \sbarm{1} & \hdots & \sbarm{6}\;  \\
    \;0 & \sbarmp{1}{1}  & \hdots  & \sbarmp{6}{1}\; \\
    \,\sbarmp{0}{2} & \sbarmp{1}{2} & \hdots & \sbarmp{6}{2}\;\\
    \;\vdots & \vdots & \vdots & \vdots \; \\
    \;0 & \sbarmp{1}{7} & \hdots & \sbarmp{6}{7}\;
\end{bmatrix}
\begin{bmatrix}
\;\mo{0}\; \\
\;\mo{1}\; \\
\;\mo{2}\; \\
\; \vdots\; \\
\;\mo{6}\; 
\end{bmatrix} 
=
\begin{bmatrix}
1\\
0\\
0\\
\vdots \\
0
\end{bmatrix}
\end{align}
In summary, the construction of a biorthogonal basis \eqref{eq: bi-orth MRA} is equivalent to find feasible solutions of \eqref{eq: LS-new} with constraint \eqref{eq: identity-cond}.\footnote{It can be shown that as long as \eqref{eq: LS-new} has a unique solution for $m_j$ given fixed $\widetilde{m_j}, \, j = 0,\cdots,6,$ \eqref{eq: identity-cond} always holds. See Section \ref{subsec: adapt-cohen}.}
Our approach to this is inspired by the approach in \cite{cohen1993compactly} for constructing compactly supported symmetric biorthogonal filters on a hexagon lattice. We next review the main scheme in \cite{cohen1993compactly} and adapt it to our setup of biorthogonal bases on the dyadic quincunx lattice.

\subsection{Summary of Cohen et al's construction}\label{subsec: cohen-summary}
We summerize the main setup and the approach in \cite{cohen1993compactly}. Consider a biorthogonal scheme consisting of three high-pass filters $m_1,m_2$ and $m_3$ and a low-pass filter $m_0$ together with their biorthogonal duals $\widetilde{m_j}$, s.t.
$m_0$ and $\widetilde{m_0}$ are $\frac{2\pi}{3}$-rotation invariant and $m_1,\, m_2,\, m_3$ and their duals are $\frac{2\pi}{3}$-rotation co-variant on a hexagon lattice.

This biorthogonal scheme satisfies the following linear system (
Lemma 2.2.2 in \cite{cohen1993compactly} )
\begin{align}\label{eq: LS}
\begin{bmatrix}
    \, \sbarm{0} &  \sbarm{1} &  \sbarm{2} & \sbarm{3}\; \\
    \;\sbarmn{0}{1} & \sbarmn{1}{1}  & \sbarmn{2}{1}  & \sbarmn{3}{1}\; \\
    \;\sbarmn{0}{2} & \sbarmn{1}{2}  & \sbarmn{2}{2}  & \sbarmn{3}{2}\; \\
    \;\sbarmn{0}{3} & \sbarmn{1}{3} & \sbarmn{2}{3} & \sbarmn{3}{3}\;
\end{bmatrix}
\begin{bmatrix}
\;\mo{0}\; \\
\;\mo{1}\; \\
\;\mo{2}\; \\
\;\mo{3}\; 
\end{bmatrix} 
=
\begin{bmatrix}
1\\
0\\
0\\
0
\end{bmatrix}.
\end{align}
 where $\V{\nu}_i = \V{\pi}_{2i},\, i = 1,2,3.$ %$\V{\nu}_1 = (\pi,0),\V{\nu}_2 = (0,\pi),\V{\nu}_3=(\pi,\pi)$.
 Let $\widetilde{\mathbf{M}}(\V{\omega})\in\mathbb{C}^{4\times 4}$ be the matrix with entries $\sbarmn{j}{i}$ and $\mathbf{m}(\V{\omega})\in\mathbb{C}^4$ be the vector with entries $m_j(\V{\omega})$ in \eqref{eq: LS}, then \eqref{eq: LS} can be written as \[\widetilde{\mathbf{M}}(\V{\omega})\, \mathbf{m} (\V{\omega})= [1,0,0,0]^\top.\]
Begin with a pre-designed $\m{1}$ with desired propery, $\m{2}$ and $\m{3}$ are determined by symmetry. Lemma 2.2.2 in \cite{cohen1993compactly} then leads to
\begin{align}\label{eq: m0-sol}
m_0(\V{\omega}) &= D^{-1}%\propto 
\left|
\begin{matrix}
    \; \sbarmn{1}{1}  & \sbarmn{2}{1}  & \sbarmn{3}{1}\; \\
    \; \sbarmn{1}{2}  & \sbarmn{2}{2}  & \sbarmn{3}{2}\; \\
    \; \sbarmn{1}{3} & \sbarmn{2}{3} & \sbarmn{3}{3}\;
\end{matrix}
\right| \notag\\
&= D^{-1}\widetilde{\mathbf{M}}_{0,0}(\V{\omega}),
\end{align}
where $\widetilde{\mathbf{M}}_{0,0}(\V{\omega})$ is the minor of $\widetilde{\mathbf{M}}(\V{\omega})$ with respect to $\sbarm{0}$ and $ D \equiv \det(\widetilde{\mathbf{M}}(\V{\omega}))\in \mathbb{C}^* = \mathbb{C}\setminus\{0\}$ does not depend on $\V{\omega}$ in \cite{cohen1993compactly}, due to the symmetry of $\widetilde{m_j}$.
%{\it Remark.} 
%For \eqref{eq: m0-sol} to hold, $m_0(\mathbf{\omega})$ and $\det(\widetilde{\mathbf{M}}_{1,1}(\V{\omega}))$ having the same phase suffices, which is implied by the symmetry of $m_0$ and $\widetilde{m_j} $'s.\\ % Both $\mo{0}$ and $\det(\widetilde{\mathbf{M}}_{1,1}(\V{\omega}))$ are $\frac{2\pi}{3}-$rotation invariant. \\

Expanding $det(\widetilde{\mathbf{M}}(\V{\omega}))$ with respect to the first column leads to the following constraint on $\m{0}$, 
\begin{align}\label{eq: biorth-eq}
m_0\sbarm{0} + m_0\sbarmn{0}{1} + m_0\sbarmn{0}{2} + m_0\sbarmn{0}{3} = 1,
\end{align}
which is the same as the identity constraint \eqref{eq: identity-cond} in our setup. 
Once \eqref{eq: biorth-eq} is solved for $\widetilde{m_0}$, $m_1,m_2$ and $m_3$ are obtained by solving the linear system \eqref{eq: LS} as $\M(\V{\omega})$ has been determined.
%\subsubsection{Solving $\m{0}$}

\subsection{Adaptation to dyadic quincunx downsampling}\label{subsec: adapt-cohen}

Cohen et al's approach can be adapted to construct biorthogonal bases in different settings; We shall apply it to our framework, even though we work with different lattices, downsampling schemes and symmetries. In particular, we adapt their approach to solve \eqref{eq: LS-new} with constraint \eqref{eq: identity-cond} where $\widetilde{m_j},\, j= 1,\cdots,6$ are pre-designed. Furthermore, by exploiting the symmetric structure of \eqref{eq: LS-new} with respect to the shifts $\V{\pi}_i,\,i=0,\cdots,7$, we derive necessary conditions for \eqref{eq: LS-new} to have a unique solution. It turns out that these will, once again, force to exhibit lack of regularity in our biorthogonal scheme.

Since \eqref{eq: LS-new} takes the same form as \eqref{eq: LS}, we adopt, for the sake of simplicity and for the rest of this paper,  the matrix and vector notations $\widetilde{\mathbf{M}}(\V{\omega}),\,\mathbf{m}(\V{\omega}) $ that helped to simplify \eqref{eq: LS}. Accordingly, we rewrite \eqref{eq: LS-new} as  \[\widetilde{\mathbf{M}}(\V{\omega})\, \mathbf{m} (\V{\omega})= [1,0,0,0,0,0,0]^\top,\]  where $\widetilde{\mathbf{M}}(\V{\omega})\in\mathbb{C}^{8\times 7}$ and $\mathbf{m}(\V{\omega})\in\mathbb{C}^7$. In addition, let $\V{b}_k \in\mathbb{R}^8,\, 0\leq k\leq 7$, whose only non-zero entry is $\V{b}_k[k] = 1$, where the indexing starts with zero. Note that $\M(\V{\omega})\,\mathbf{m}(\V{\omega}) = \V{b}_0\in\mathbb{R}^8$ is over-determined; it has a unique solution of $m_j$ if and only if 

\begin{enumerate}[leftmargin=.5in]
\item[\mylabel{cond: 1}{(\ref{sec: solve-quincunx}.i)}] $\M(\V{\omega})$ is full rank,
\item[\mylabel{cond: 2}{(\ref{sec: solve-quincunx}.ii)}] $[\M(\V{\omega}), \V{b}_0]$ is singular,
\end{enumerate}
where we use the notation $[\;]$ for the concatenation of $\M(\V{\omega})$ and $\V{b}_0$ into a $8\times 8$ matrix. %denote by $[\M(\V{\omega}), \V{b}_0]$ the $8\times 8$ matrix obtained by adding the 8-th column $\V{b}_0$ to the $8\times 7$ matrix $\M(\V{\omega})$.
The matrix $\M(\V{\omega})$ is structured such that each row is associated with a shift $\V{\pi}_i,\,i=0\cdots,7$ and each column is associated with a dual function $\m{j},\,j=0,\cdots,7$. In particular, $\M(\V{\omega})$ depends on the value of $\widetilde{m_j}$ at $\V{\omega}$ and its shifts $\V{\omega}+\V{\pi}_i$. We denote a submatrix of $\M(\V{\omega})$ containing all but the row associated with $\V{\pi}_k$ (respectively, the column associated with $\m{k}$) as $\M[\widehat{k},:](\V{\omega})$ (respectively, $\M[:,\widehat{k}](\V{\omega})$).
In particular, we denote $\M[\widehat{0},\widehat{0}](\V{\omega})$ as $\Msub(\V{\omega})$.

We have the following observations for $\M(\V{\omega})$.
\begin{lemma}\label{lem: subM-singular}
 $\forall \V{\omega}\in S_0$, if \eqref{eq: LS-new} is solvable, then $\M[\widehat{0},:](\V{\omega})$ is singular.
\end{lemma}
\noindent{\it Proof.}
If \eqref{eq: LS-new} is solvable, then condition \ref{cond: 2} holds, which implies that $\det([\M(\V{\omega}),\V{b}_0]) = 0$. Expanding the determinant with respect to the last column $\V{b}_0$ yields $\det(\,\M[\widehat{0},:](\V{\omega})\,) = 0$.\qed
%If \eqref{eq: LS-new} has a solution, then $\forall \V{\omega}$,  $[1,0,\cdots,0]^\top\in \mathbb{R}^8$ is a linear combination of the columns of $\M$ hence the solution $\mathbf{m} \in Null(\M[\widehat{0},:])$ and it is non-zero. This implies that $\M[\widehat{0},:]$ is singular.\qed\\[1em]

\begin{lemma}\label{lem: M-symmetry}
$\M(\V{\omega}),\,\M(\V{\omega}+\V{\pi}_2),\,\M(\V{\omega}+\V{\pi}_4)$ and $\M(\V{\omega}+\V{\pi}_6)$ are the same up to row permutations. \eqref{eq: LS-new} holds $\forall \, \V{\omega}$ if and only if 
%\begin{align*}
$\M(\V{\omega}) \big[\,\mathbf{m}(\V{\omega}),\mathbf{m}(\V{\omega}+\V{\pi}_2),\mathbf{m}(\V{\omega}+\V{\pi}_4),\mathbf{m}(\V{\omega}+\V{\pi}_6)\,\big] = \big[\,\V{b}_0,\V{b}_2,\V{b}_4,\V{b}_6\,\big].$
%\end{align*}
\end{lemma}

\noindent{\it Remark}.
If we consider $\M(\V{\omega})$ a matrix-valued function of $\V{\omega}$, then the conditions \ref{cond: 1} and \ref{cond: 2} are both pointwise, yet Lemma \ref{lem: M-symmetry} shows that the set of points $\{\V{\omega},\V{\omega}+\V{\pi}_2,\V{\omega} + \V{\pi}_4, \V{\omega}+\V{\pi}_6\}$ are linked together by the symmetry in $\M(\V{\omega})$. 

Due to condition \ref{cond: 1}, $\forall\,\V{\omega}$, $\exists\, k_{\V{\omega}}$ depending on $\V{\omega}$ such that $\M[\widehat{k_{\V{\omega}}},:](\V{\omega})$ is non-singular. Lemma \ref{lem: subM-singular} implies that $\widehat{k_{\V{\omega}}}\neq 0$;\footnote{By symmetry, we have the stronger result $k_{\V{\omega}}\not\in\{0,2,4,6\}$. Indeed, Lemma \ref{lem: subM-singular} and Lemma \ref{lem: M-symmetry} together imply that $\M[\widehat{k},:](\V{\omega}),\,k=0,2,4,6$ are singular. Therefore, $\widehat{k_{\V{\omega}}} \in\{1,3,5,7\}$ and thus $\M[\widehat{k_{\V{\omega}}},:](\V{\omega})$ contains all rows associated with shifts $\V{\pi}_{2i},\, i = 0,\cdots,3$. 
} 
therefore we may apply Cramer's rule to $\M[\widehat{k_{\V{\omega}}},:](\V{\omega})$, as in Section \ref{subsec: cohen-summary}, and obtain
the following expression of $m_0(\V{\omega})$
%there is a unique row $\M[\widehat{k_{\V{\omega}}},:],\,k_\omega\in\{2,\cdots,8\}$ such that removing it from $\M$ gives a non-singular square matrix $\M[\widehat{k_{\V{\omega}}},:]$. By Cramer's rule, 
\begin{align}\label{eq: m0-cramer}
m_0(\V{\omega}) = \det(\,\Msub[\widehat{k_{\V{\omega}}},:](\V{\omega})\,)/\det(\,\M[\widehat{k_{\V{\omega}}},:](\V{\omega})\,).
\end{align}
Moreover, based on \eqref{eq: m0-cramer}, the identity condition \eqref{eq: identity-cond} on $m_0(\V{\omega})$ and $\m{0}$ can be derived in the same way as \eqref{eq: biorth-eq} by expanding $\det(\,\M[\widehat{k_{\V{\omega}}},:](\V{\omega})\,)$.

%We then discuss how to design $\widetilde{m_j}$ with more symmetry and solve the corresponding system \eqref{eq: LS-new}.

\subsection{Discontinuity of $\m{j}$}\label{subsec: discontinuity}
In this subsection, we show our main result that for \eqref{eq: LS-new} to be uniquely solvable, the pre-designed $\widetilde{m_j}$ have to be discontinuous as soon as they satisfy mild symmetry conditions and concentration of support on $C_j$.
%when the support of $\widetilde{m_j}$ concentrates in the direction of $C_j$ and a minimum symmetry of $\widetilde{m_j}$ is required.

We assume that
$|\m{1}|$ and $|\m{6}|$ are symmetric with respect to the diagonal $\omega_1=\omega_2$, i.e.
\begin{align}\label{eq: sym-m16}
 |\widetilde{m_1}(\V{\omega})| = |\widetilde{m_6}(\V{\omega}')|\quad \forall\, \omega_1=\omega_2', \,\omega_2=\omega_1',
\end{align}
and likewise for $\m{3}$ and $\m{4}$,
\begin{align}\label{eq: sym-m34}
|\widetilde{m_3}(\V{\omega})| = |\widetilde{m_4}(\V{\omega}')|\quad \forall\, \omega_1=-\omega_2', \,\omega_2=-\omega_1'.
\end{align}
%Same as in Section \ref{subsec: cohen-summary} , we first compute $m_0$ and assume that $\M$ is full rank, otherwise \eqref{eq: LS-new} has infinitely many solutions. Moreover, $\M[2:8,:]$ is singular. 
%Proposition \ref{prop: feasibility} provides a necessary condition such that the numerical optimization solving $\widetilde{m_0}$ is feasible.
\begin{figure}
\centering
\includegraphics[width = .3\textwidth]{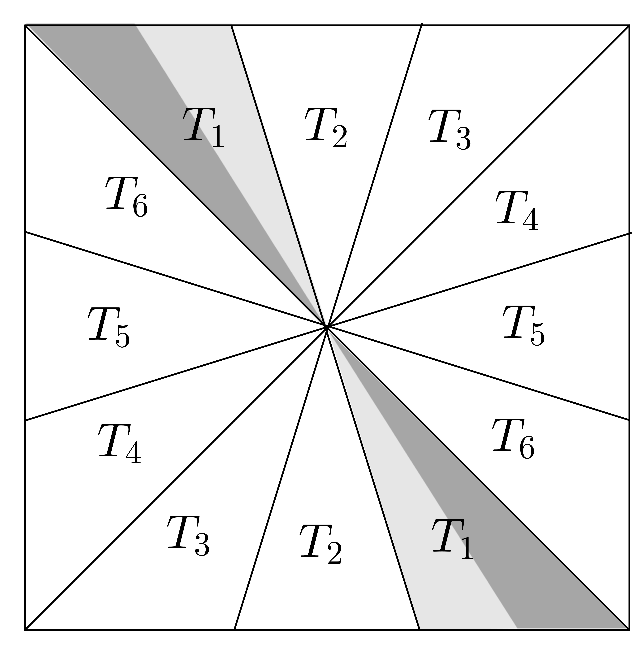}
\caption{Partition of frequency square in six directions, where the essential support of $\m{i}$ is contained in each pair of triangles $T_i$. The pair of dark grey triangles is $T_1^-$ and the light grey pair is $T_1^+$.}
\label{fig: partition 2}
\end{figure}
%Let pairs of triangles $T_i$ in Fig.\ref{fig: partition 2} contain the essential support of $\widetilde{m_i},\,i=1,\cdots,6$.
%\eqref{eq: LS-new} takes a similar form to \eqref{eq: LS}, but with $\M\in\mathbb{C}^{8\times 7}$, which is an over-determinant linear system.
In what follows, we introduce a triangular partition of $S_0=[-\pi,\pi)\times[-\pi,\pi)$ in the frequency plane and define formally the concentration of the support of the $\widetilde{m_j}$.

\noindent{\bf Definition.}
The {\it  domination-support} $\Omega_j$ of a function $\widetilde{m_j}$ (with respect to the other $m_i$, $i\neq j$) is the set $\{\V{\omega}:\,|\widetilde{m_j}(\V{\omega})|> |\widetilde{m_i}(\V{\omega})|,\,\forall i\neq j\}$. \vspace{.5em}

Let $T_j$ be pairs of triangles shown in Figure \ref{fig: partition 2}, defined such that $C_j\subset T_j,\, j = 1,\cdots,6.$ Consider the decompositions $T_j = T_j^-\bigcup T_j^+$, where $T_j^-, T_j^+$ are halves of $T_j$ adjacent to its neighboring triangles $T_i$ in the counter clockwise and clockwise directions respectively.\\[.5em]
\noindent{\bf Definition.}  $\widetilde{m_j}$ {\it concentrates} in $T_j$ for $j = 1,\cdots,6$ if 
\begin{itemize}
\item[(i)] $\Omega_j\subset T_j$;
\item[(ii)] $\text{supp}(\widetilde{m_j})\subset T_{j-1}^+\bigcup T_j\bigcup T_{j+1}^-$ and $\int_\Omega|\widetilde{m_j}| > \int_{\Omega'}|\widetilde{m_j}|, \forall\, \Omega\subset T_j\bigcap\text{supp}(\widetilde{m_j})$ s.t. $|\Omega|>0$, where $\Omega' \subset T_{j-1}^+\bigcup T_{j+1}^-$ is symmetric to $\Omega$ with respect to the boundary of $T_j$.
\end{itemize}
In other words, for $\widetilde{m_j}$ to concentrate in $T_j$, $\widetilde{m_j}$ should be ``mainly" supported in $T_j$ (condition (i)) and ``decay" properly outside of $T_j$ (condition (ii)).

\begin{figure}
\centering
\includegraphics[width = .5\textwidth]{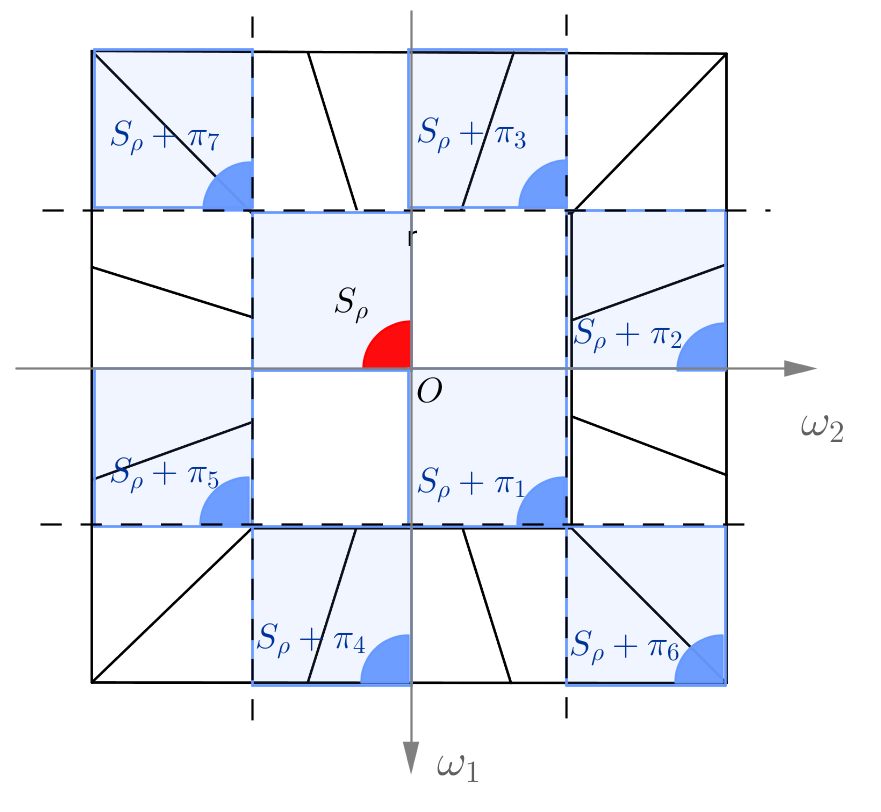}
\caption{$S_{\rho}$ and its shifts}
\label{fig: S-shifts}
\end{figure}
We say $\widetilde{m_0}$ concentrates in $C_0$ if $\Omega_0\subset C_0$.
For $m_0$, we impose the natural requirement that, for some (possibly small) $\rho>0$, we have $|m_0(\V{\omega})| > 0,\, \forall\, |\V{\omega}| < \rho$.
Given these constraints on the support of $\widetilde{m_j}$ and $m_0$, we examine the consequences of the singularity condition on $\M[\widehat{0},:](\V{\omega})$ from Lemma \ref{lem: subM-singular}, specifically in the domain $S_{\rho} = \{(\omega_1,\omega_2)|\;|\V{\omega}| < \rho, \omega_1 <0,\,\omega_2<0\}$, see the red zone in Figure \ref{fig: S-shifts}. 

Let $\mrow{i}(\V{\omega}) = [\widetilde{m_1}(\V{\omega}+\V{\pi}_i)\, \cdots,\,\widetilde{m_6}(\V{\omega}+\V{\pi}_i)]\in\mathbb{C}^6,\, i = 0,\cdots,7$ be the rows of $\M[:,\widehat{0}](\V{\omega})$. 
\begin{lemma}\label{lem: rank1}
If $\V{\omega}\in S_\rho$ s.t. \eqref{eq: identity-cond} holds and $\M[\widehat{0},:](\V{\omega})$ is singular, then  $rank(\mrow{1},\mrow{7}) = 1$ and $rank(\mrow{3},\mrow{5}) = 2$ or $rank(\mrow{3},\mrow{5}) = 1$ and $rank(\mrow{1},\mrow{7}) = 2$.
\end{lemma}
Lemma \ref{lem: rank1} can be proved by analyzing the linear dependency and independency between the $\mrow{i}$ on $S_\rho$, since the $\mrow{i}$ have known locations of zero entries when $\rho$ is small due to the concentration of the $\widetilde{m_j}$.
For the full proof of Lemma \ref{lem: rank1}, see Appendix \ref{app: lemmas}.

The concentration of $\m{3}$ and $\m{4}$ in $T_3$ and $T_4$ and their symmetry together imply that  $rank(\mrow{3},\mrow{5}) \neq 1\; a.e.$ on $S_\rho$ (see Lemma \ref{lem: full-rank-m35} in Appendix \ref{app: discontinuity}), hence $rank(\,\mrow{1}(\V{\omega}),\mrow{7}(\V{\omega})\,) = 1\; a.e.$ on $S_\rho$.
Therefore, $\mp{1}{1}, \mp{6}{1}$ in $\mrow{1}(\V{\omega})$ and the corresponding $\mp{1}{7}, \mp{6}{7}$ in $\mrow{7}(\V{\omega})$ on $S_\rho$ are linearly related. Furthermore, we can show that $\m{6}=0\; a.e.$ on $S_\rho + \V{\pi}_1\cap \{\omega_1 < \omega_2\}$ (see Proposition \ref{prop: zero-corner} in Appendix \ref{app: discontinuity}), if $\m{0}$, $\m{1}$ and $\m{6}$ concentrate in $C_0$, $T_1$ and $T_6$ respectively. Therefore, if $\m{6}$ is continuous, $\widetilde{m_6}(\dfrac{\pi}{2},\dfrac{\pi}{2}) = 0$; the same holds for $\m{1}$ and for $(-\dfrac{\pi}{2},-\dfrac{\pi}{2})$ as well by symmetry.% and by symmetry $\m{1}$ and $\m{6}$ will be zero at $(\frac{\pi}{2},\frac{\pi}{2})$ and $(-\frac{\pi}{2},-\frac{\pi}{2})$ if they are continuous.% at $(\frac{\pi}{2},\frac{\pi}{2})$ or $(-\frac{\pi}{2},-\frac{\pi}{2})$.
The following theorem summarizes our main result.
\begin{theorem}\label{thm: thm}
%\begin{proposition}\label{prop: continuity}
If the $\widetilde{m_j}$ concentrate in $T_j$ for $j= 1,3,4,6$ with symmetries \eqref{eq: sym-m16} and \eqref{eq: sym-m34} and $\widetilde{m_0}$ concentrates in $C_0$, then $\m{1},\m{6}$ cannot be continuous at both $(\dfrac{\pi}{2},\dfrac{\pi}{2})$ and $(-\dfrac{\pi}{2},-\dfrac{\pi}{2})$ for \eqref{eq: LS-new} to have a unique solution of $m_j$ s.t. $\exists\,\rho>0,\,m_0(\V{\omega})$ is non-zero on $|\V{\omega}|<\rho$.
\end{theorem}
%\end{proposition}
\noindent{\it Proof. }
If $\m{1}$ and $\m{6}$ are both continuous at $(\dfrac{\pi}{2},\dfrac{\pi}{2})$ and $(-\dfrac{\pi}{2},-\dfrac{\pi}{2})$, then $\widetilde{m_1}(\dfrac{\pi}{2},\dfrac{\pi}{2}) = \widetilde{m_1}(-\dfrac{\pi}{2}, -\dfrac{\pi}{2}) = \widetilde{m_6}(\dfrac{\pi}{2},\dfrac{\pi}{2}) = \widetilde{m_6}(-\dfrac{\pi}{2},-\dfrac{\pi}{2}) = 0$. %$\widetilde{m_1}(\dfrac{\pi}{2},\dfrac{\pi}{2}) = \lim_{\alpha\rightarrow 1^-}\widetilde{m_1}(\V{\omega}(\alpha)) = 0$, where $\{\V{\omega}(\alpha),\,0\leq \alpha<1\} \subset S_\rho + \V{\pi}_1$ and $\V{\omega}(1) = (\frac{\pi}{2},\frac{\pi}{2})$. By symmetry, we have $\widetilde{m_6}(\dfrac{\pi}{2},\dfrac{\pi}{2}) = 0$. Similarly, the continuity at $(-\dfrac{\pi}{2},-\dfrac{\pi}{2})$ implies $\widetilde{m_1}(-\dfrac{\pi}{2},-\dfrac{\pi}{2}) = \widetilde{m_6}(-\dfrac{\pi}{2},-\dfrac{\pi}{2}) = 0$. 
Therefore, $\mrow{1}(\V{0}) = \mrow{7}(\V{0}) = \mathbf{0}$ at the origin which results in contradiction with Lemma \ref{lem: rank1}.\qed%and from \eqref{eq: m0C} $m_0^C(0)=0$ so that $m_0(0)=0$, %  On the other hand, Proposition\ref{prop: origin-det} implies that $m_0(0) = 0$ as $a = |\widetilde{m_1}(\pi_1)| = 0$, which results in contradiction.
%The following theorem summarizes our main result.
%\begin{theorem}\label{thm: thm}
%If the $\m{j}$ concentrate in $T_j$, $\m{0}$ concentrates in $C_0$ and the symmetries \eqref{eq: sym-m16} and \eqref{eq: sym-m34} hold, then \eqref{eq: LS-new} does not have a unique solution of $m_j(\V{\omega})$ s.t. $\exists\,\rho>0,$ $|m_0(\V{\omega})|> 0, \, \forall |\V{\omega}|<\rho$, if $\m{1}$ and $\m{6}$ are continuous.
%\end{theorem}

\section{Numerical construction of biorthogonal bases}\label{sec: solve-quincunx}
In this section, we develop a numerical construction of biorthogonal bases on a dyadic quincunx lattice following an approach similar to Cohen et al. We first design $\m{j},\,j = 1,\cdots,6,$ on the canonical frequency square $S_0 = [-\pi,\pi)\times[-\pi,\pi)$associated with the lattice $\mathbb{Z}^2$, then solve for $m_0,\widetilde{m_0}$ and $m_j$ on $S_0$ in order with respect to \eqref{eq: LS-new} and \eqref{eq: identity-cond}.
%we focus on solving $m_i$'s and $\widetilde{m_0}$ in \eqref{eq: LS-new} given pre-designed $\m{i},\,i=1,\cdots,6$. %Assume $\m{i},\,i=1,\cdots,6$ satisfy weak constraints on the direction selectivity of their support.

\subsection{Design of input $\m{j}$}\label{sec: phase-design}
In this sub-section, we construct $\m{j},\, j = 1,\cdots,6,$ which concentrate in $T_i$.
Specifically, following the orthonormal construction in \cite{yin2014orthshear}, we consider $\m{j}$ in the form 
\begin{align}\label{eq: m-form}
\m{j} = e^{-i\V{\eta}_j^\top\V{\omega}}|\m{j}|,\quad j = 1,\cdots,6,
\end{align}
where $\V{\eta}_j\in\mathbb{Z}^2$ is the phase constant of $\widetilde{m_j}$. In addition to the symmetry of pairs $(|\widetilde{m_1}|, |\widetilde{m_6}|)$ and $(|\widetilde{m_3}|, |\widetilde{m_4}|)$ assumed in Section \ref{subsec: discontinuity}, we further require that $|\widetilde{m_2}|$ and $|\widetilde{m_5}|$ are symmetric with respect to the $\omega_1$-axis and $\omega_2$-axis accordingly.
% where $|\m{j}|$'s have symmetries with respect to the two diagonals and the two axis. 
Figure \ref{fig: mjdual} shows a design of $|\m{j}|$ that has these strong symmetries.
 
The symmetries of $(|\widetilde{m_1}|, |\widetilde{m_6}|)$ leads to constraints on the phase constants $\V{\eta}_j$ introduced in \eqref{eq: m-form}.
% We want to design the phase $\V{\eta}_k$ such that $m_0(\V{\omega}) > 0, \; \forall \omega\in S_1$. This is the same as requiring $\Msub$ to be full rank. We first show the necessary conditions on phases $\V{\eta}$ of the full rank requirement on $\Msub$.
 
\begin{lemma}\label{lem: phase-ineq}
If $\exists\,\V{\omega}\in D_1:=\{\omega_1=\omega_2,\,\omega_1\in(-\frac{\pi}{2},0)\},\,s.t. \,|m_0(\V{\omega})|\neq 0,$ then $(\V{\eta}_1-\V{\eta}_6)^\top (\V{\pi}_6-\V{\pi}_7)\neq 0(\text{mod}\,2\pi)$. 
\end{lemma} 

Because $m_0(\V{\omega})$ can be expressed as in \eqref{eq: m0-cramer}, $|m_0(\V{\omega})|\neq 0$ is equivalent to $\det(\,\Msub[\widehat{k_{\V{\omega}}},:](\V{\omega})\,)\neq 0$, i.e. $\Msub(\V{\omega})$ is full rank. The constraint on $\V{\eta}_1$ and $\V{\eta}_6$ then follows from substituting non-zero entries of $\Msub(\V{\omega})$ by \eqref{eq: m-form} and consider the linear dependency of the columns in $\Msub(\V{\omega})$. For the full proof of Lemma \ref{lem: phase-ineq}, see Appendix \ref{app: input design}.

Similarly, if $\exists\,\V{\omega}\in \{\omega_1 = \omega_2,\, \omega_1\in(0,\frac{\pi}{2})\},\, s.t.\, |m_0(\V{\omega})| \neq 0$, then $(\V{\eta}_1-\V{\eta}_6)^\top (\V{\pi}_6-\V{\pi}_1)\neq 0(\text{mod}\,2\pi)$. These two conditions are equivalent to 
\begin{align*}
(\V{\eta}_1-\V{\eta}_6)^\top(\pi/2,\pi/2)\neq 0 (\text{mod}\,2\pi)\tag{\bf c1.1}
\end{align*}
since $\V{\eta}_1,\,\V{\eta}_6\in\mathbb{Z}^2$.
Considering the other diagonal segment $\{\omega_2 = -\omega_1, |\omega_1| <\frac{\pi}{2}\}$ and the symmetry of $(|\widetilde{m_3}|, |\widetilde{m_4}|)$, we similarly obtain
\begin{align*}
(\V{\eta}_3-\V{\eta}_4)^\top(-\pi/2,\pi/2)\neq 0 (\text{mod}\, 2\pi)\tag{\bf c1.2}
\end{align*}
%from the full rank condition.

Next, we consider $\widetilde{m_0}(\V{0})$ and investigate $\Msub(\V{\omega})$ at the origin.
\begin{proposition}\label{prop: origin-det}
%If $|\widetilde{m_1}(\V{\pi}_1)| = |\widetilde{m_1}(\V{\pi}_7)| = |\widetilde{m_3}(\V{\pi}_3)| = |\widetilde{m_3}(\V{\pi}_5)|$ and $|\widetilde{m_1}(\V{\pi}_6)|= | \widetilde{m_3}(\V{\pi}_6)|$, then 
If $|\widetilde{m_0}(\V{0})|\neq 0,$ then $\V{\pi}_1^\top(\V{\eta}_1-\V{\eta}_6)\neq \pi(\text{mod}\,2\pi)$ or $\V{\pi}_3^\top(\V{\eta}_3-\V{\eta}_4)\neq \pi(\text{mod}\,2\pi)$. 
\end{proposition}
\noindent{\it Remark.}
The proof of Proposition \ref{prop: origin-det} is similar to that of Lemma \ref{lem: phase-ineq} but more involved. See Appendix \ref{app: input design} for the full proof.

We propose the following set of phases such that ({\bf c1.1}) and ({\bf c1.2}) as well as the necessary condition from Proposition \ref{prop: origin-det} are all satisfied,
\begin{align}\label{eq: phase-sol}
\V{\eta}_1 = (0,0),\; \V{\eta}_2 = (-1,1),\; \V{\eta}_3 = (0,2),\notag\\
\V{\eta}_4 = (1,0),\; \V{\eta}_5 = (0,-1),\; \V{\eta}_6 = (0,1).
\end{align}

The design of $\m{j}$ in the form of \eqref{eq: m-form} with phases \eqref{eq: phase-sol} introduced here do not guarantee that \eqref{eq: LS-new} has a unique solution. We will see the necessary and sufficient conditions that $\m{j}$ have to satisfy in the next subsection given by Proposition \ref{prop: m0_formula}.
\begin{comment}
where
\begin{align*}
%\label{eq: phase-constraint}
%&(\V{\eta}_1-\V{\eta}_2)^\top(-\pi/2, \pi/2) = (\V{\eta}_5-\V{\eta}_6)^\top(-\pi/2,\pi/2) = \pi\, (\text{mod}\, 2\pi)\notag\\
%&(\V{\eta}_2-\V{\eta}_3)^\top(\pi/2,\pi/2) = (\V{\eta}_4-\V{\eta}_5)^\top(\pi/2,\pi/2) = \pi\, (\text{mod}\, 2\pi)\\
%&
(\V{\eta}_3-\V{\eta}_4)^\top(-\pi/2, \pi/2) &=-\pi/2\,(\text{mod}\,2\pi)\notag\\
 (\V{\eta}_6 - \V{\eta}_1)^\top(\pi/2,\pi/2) &= \pi/2\, (\text{mod}\,2\pi)\notag
\end{align*}
%where we require strong orthogonal constraints on pair of shifts corresponding to $\widetilde{m}$ function with non-diagonal common boundary and weaker constraints on $(\V{\eta}_1,\V{\eta}_6)$ and $(\V{\eta}_3,\V{\eta}_4)$. A solution to \eqref{eq: phase-constraint} is 
\end{comment}

\subsection{Solving \eqref{eq: LS-new} and \eqref{eq: identity-cond} for $m_0,\widetilde{m_0}$ and $m_j$}\label{subsec: compute-m0}
Once $\m{j}, j= 1,\cdots, 6$ are fixed on $S_0$, \eqref{eq: LS-new} can be reformulated as follows,
\begin{align}\label{eq: mj-eq}
\M[:,\widehat{0}](\V{\omega}) 
\begin{bmatrix}
m_1(\V{\omega})\\
m_2(\V{\omega})\\
m_3(\V{\omega})\\
m_4(\V{\omega})\\
m_5(\V{\omega})\\
m_6(\V{\omega})
\end{bmatrix}
= \V{b}_0 - m_0(\V{\omega})
\begin{bmatrix}
 \sbarm{0}\\
 0\\
\sbarmp{0}{2}\\
0\\
\sbarmp{0}{4}\\
0\\
\sbarmp{0}{6}\\
0
\end{bmatrix} \doteq \V{b}_0'(\V{\omega}),
\end{align}
where $\M[:,\widehat{0}](\V{\omega})$ is completely determined by $\m{j},\,j = 1,\cdots,6$ 
and $m_j,\,j=1,\cdots,6$ can be uniquely solved on $S_0$ if and only if $\forall\V{\omega}\in S_0$
\hspace{-1em}
\begin{enumerate}[leftmargin=.5in]
\item[\mylabel{cond: i}{(\ref{subsec: compute-m0}.i)}] $\M[:,\widehat{0}](\V{\omega})$ is full rank,%
\item[\mylabel{cond: ii}{(\ref{subsec: compute-m0}.ii)}] $\V{b}_0'(\V{\omega})$ is in $col\big(\M[:,\widehat{0}](\V{\omega})\big)$, the column space of $\M[:,\widehat{0}](\V{\omega})$.%
\end{enumerate}
%Because , \ref{cond: i} can be easily checked by computing $\det(\M[:,\widehat{0}]^\top\M[:,\widehat{0}])$ and see if it is non-zero. Next, we show that (ii) provides an explicit way of computing $m_0(\V{\omega})$ without knowing $\m{0}$.
Next, we show that \ref{cond: ii} breaks down to constraints on two submatrices of $\M[:,\widehat{0}](\V{\omega})$ and quadruples $\big( m_0(\V{\omega}),m_0(\V{\omega}+\V{\pi}_2), m_0(\V{\omega} +\V{\pi}_4), m_0(\V{\omega}+\V{\pi}_6) \big)$, $\big( m_0(\V{\omega}+\V{\pi}_1),m_0(\V{\omega}+\V{\pi}_3), m_0(\V{\omega} +\V{\pi}_5), m_0(\V{\omega}+\V{\pi}_7) \big)$.
\begin{proposition}\label{prop: m0_formula}
Let $\M[odd,\widehat{0}](\V{\omega}),\M[even,\widehat{0}](\V{\omega})\in\mathbb{C}^{4\times6}$ be the submatrices of $\M[:,\widehat{0}](\V{\omega})$ consisting of odd and even indexed rows respectively. $\forall\V{\omega}\in S_0$, suppose {\rm\ref{cond: i}} holds, then {\rm\ref{cond: ii}} holds if and only if $rank(\,\M[odd,\widehat{0}](\V{\omega})\,) = rank(\,\M[even,\widehat{0}](\V{\omega})\,) = 3$ and 
\begin{align}\label{eq: m0-even-null}
[m_0(\V{\omega}),m_0(\V{\omega}+\V{\pi}_2), m_0(\V{\omega} +\V{\pi}_4), m_0(\V{\omega}+\V{\pi}_6)]\, \M[even,\widehat{0}](\V{\omega}) = \V{0},
\end{align}
\begin{align}\label{eq: m0-odd-null}
[m_0(\V{\omega}+\V{\pi}_1),m_0(\V{\omega}+\V{\pi}_3), m_0(\V{\omega} +\V{\pi}_5), m_0(\V{\omega}+\V{\pi}_7)] \,\M[odd,\widehat{0}](\V{\omega}) = \V{0}.
\end{align}
\end{proposition}

For the proof of Proposition \ref{prop: m0_formula}, see Appendix \ref{app: solving}.

\noindent{\it Remark.} 
Note that the submatrices $\M[odd,\widehat{0}](\V{\omega})$ and $\M[even,\widehat{0}](\V{\omega})$ are dual to each other under the shift of variable $\V{\omega}\mapsto \V{\omega}+\V{\pi}_i$, when $i$ is odd. Therefore, the constraints $rank(\,\M[even,\widehat{0}](\V{\omega})\,) = 3$ and \eqref{eq: m0-even-null} from Proposition \ref{prop: m0_formula} are sufficient for \ref{cond: ii} to hold on $S_0$. Furthermore, because $\M[even, \widehat{0}](\V{\omega})$ and $(\V{\omega},\V{\omega}+\V{\pi}_2,\V{\omega}+\V{\pi}_4, \V{\omega} + \V{\pi}_6)$ are invariant to the shift of variable $\V{\omega}\mapsto\V{\omega}+\V{\pi}_i$ when $i$ is even, we only need to consider the constraints above on the subset $[-\pi,0)\times[-\pi,0)$ of $S_0$.
%Both conditions \eqref{eq: m0-even-null} and \eqref{eq: m0-odd-null} hold $\forall \V{\omega}\in[-\pi,\pi)\times[-\pi,\pi)$ if and only if \eqref{eq: m0-even-null} holds $\forall \V{\omega}\in [-\pi,0)\times[-\pi,0)$.

In summary, $\M[:,\widehat{0}](\V{\omega})$ (or equivalently $\widetilde{m_j}$) has to satisfy the following rank constraints on $[-\pi,0)\times[-\pi,0)$ for \eqref{eq: mj-eq} to be uniquely solvable on $S_0$, 
\begin{align}\label{eq: Mranks}
rank(\,\M[:,\widehat{0}](\V{\omega})\,) = 6,\, rank(\,\M[even,\widehat{0}](\V{\omega})\,) = 3.
\end{align}
In practice, the rank constraints are hard to impose while designing $\widetilde{m_j}$, in our numerical experiments, we therefore first construct $\widetilde{m_j}$ following the design in Section \ref{sec: phase-design} and then check if these rank constraints are satisfied, see step 1. in \ref{alg}.

If \eqref{eq: Mranks} holds, the vector $[\, m_0(\V{\omega}),m_0(\V{\omega}+\V{\pi}_2), m_0(\V{\omega} +\V{\pi}_4), m_0(\V{\omega}+\V{\pi}_6) \,]$ can be uniquely determined by \eqref{eq: m0-even-null} {\it up to a constant factor $a_{\V{\omega}} $}, since it is orthogonal to the column space of $\M[even,\widehat{0}](\V{\omega})$ of co-dimension 1. In particular, we obtain $m_0(\V{\omega})$ on $S_0$ by solving \eqref{eq: m0-even-null} independently at each $\V{\omega}$ on $[-\pi,0)\times[-\pi,0)$, see step 2. in \ref{alg}. As the constant $a_{\V{\omega}}$ can change drastically as $\V{\omega}$ changes, there is potential lack of regularity of $m_0(\V{\omega})$ as an artifact of the algorithm. Figure \ref{fig: m_0} shows an $m_0(\V{\omega})$ computed in this way, which has discontinuous phase due to $a_{\V{\omega}}$. Fortunately, this irregularity is an artifact that can be removed as suggested by the following proposition.

\begin{proposition}\label{prop: mc}
If $\m{j}, m_j(\V{\omega}),\,  j = 0,1,...,6$ satisfy \eqref{eq: LS-new} and \eqref{eq: identity-cond}, 
%$m_0^C(\V{\omega}),$ $\,\mc{0}, m_i(\V{\omega}),\,i = 1,...,6$ do. More generally, 
%$m_0(\V{\omega})c(\V{\omega}),\,\m{0}c(\V{\omega})^{-1}, m_i(\V{\omega}),\,i=1,...,6$ 
then $m_0'(\V{\omega})\doteq m_0(\V{\omega})c(\V{\omega}), \widetilde{m_0}'(\V{\omega})\doteq \m{0}\overlinespace{c}(\V{\omega})^{-1}$ together with the same $m_j(\V{\omega}),\m{j}, \, j = 1,\cdots,6$ 
satisfy \eqref{eq: LS-new} and \eqref{eq: identity-cond} if $ c(\V{\omega}) = c(\V{\omega}+\V{\pi}_2)=c(\V{\omega}+\V{\pi}_4) = c(\V{\omega}+\V{\pi}_6) \neq 0$, i.e. $c(\V{\omega})$ is $\pi$-periodic in both $\omega_1$ and $\omega_2$.
\end{proposition}
\noindent{\it Proof. }
It follows from the observation that $m_0'(\V{\omega})\overline{\widetilde{m_0}'}(\V{\omega} + \V{\pi}_i) = m_0(\V{\omega})\sbarmp{0}{i},\,$ when $i$ is even.\qed

%\begin{align*}
%m_0^c(\V{\omega})\overline{\widetilde{m_0}^c}(\V{\omega} + \V{\pi}_i) = m_0(\V{\omega})\sbarmp{0}{i},\quad & i \text{ is even,}\hspace{4em}\qed
%\hspace{1em}m_j^c(\V{\omega})\overline{\widetilde{m_j}^c}(\V{\omega} + \V{\pi}_i) = \big(m_j(\V{\omega})\sbarmp{j}{i}\big)c(\V{\omega})c(\V{\omega}+\V{\pi}_1)^{-1},\quad & i \text{ is odd.}\qedhere\hspace{2em}\qed
%\end{align*}

\noindent{\it Remark.}
In practice, we use Proposition \ref{prop: mc} compensate for irregularities introduced by the arbitrary $a_{\V{\omega}}$; After $m_0(\V{\omega})$ is solved, we can choose $c(\V{\omega})$ $\pi$-periodic in both $\omega_1,\,\omega_2$ such that $m_0'(\V{\omega})$ has improved regularity and use $m_0'(\V{\omega})$ as the ``regularized" $m_0(\V{\omega})$ for the rest of the construction.
%we first solve $\m{0}$, step 3. in \ref{alg}, then we can choose $c(\V{\omega})$ $\pi$-periodic in both $\omega_1,\,\omega_2$ such that $(m_0',\widetilde{m_0}')$ have improved regularity and replace $m_0(\V{\omega})$ and $\widetilde{m_0}(\V{\omega})$ by $m_0'(\V{\omega})$ and $\widetilde{m_0}'(\V{\omega})$ respectively.

To obtain $\m{0}$ on $S_0$, we solve the identity condition \eqref{eq: identity-cond} on $[-\pi, 0)\times[-\pi, 0 )$ for the quadruple $(\m{0}, \mp{0}{2}, \mp{0}{4}, \mp{0}{6})$. Note that \eqref{eq: identity-cond} is the same as \eqref{eq: biorth-eq} in Section \ref{subsec: cohen-summary}. 
According to Lemma 3.2.1 in \cite{cohen1993compactly}, by {\it Hilbert's Nullstellensatz} \eqref{eq: biorth-eq} has a solution if and only if there does not exist $(z_1,z_2)\in (\mathbb{C}^*)^2,\, \mathbb{C}^* = \mathbb{C}\setminus\{0\}$\, s.t. $(\pm z_1,\pm z_2)$ are all vanishing points of the $z$-transform of $m_0$. Unfortunately this is not very constructive: in general, there is no efficient algorithm to solve {\it Hilbert's Nullstellensatz}.%, and how \eqref{eq: biorth-eq} is solved exactly is not mentioned in \cite{cohen1993compactly}.
%In Section \ref{sec: numerics}, we propose an optimization problem that solves $\widetilde{m_0}$, where \eqref{eq: m0-sol} (which is the same as \eqref{eq: identity-cond}) serves as a linear constraint.
%As mentioned at the end of Section \ref{subsec: cohen-summary}, we are not aware of existing efficient algorithms that solve \eqref{eq: identity-cond}. 

Our approach here is to reformulate solving $\m{0}$ under the condition \eqref{eq: identity-cond} as an optimization problem where \eqref{eq: identity-cond} serves as a linear constraint. In particular, on a $2N\times 2N$ regular grid $\G=\{\V{\omega}_i\}_{i=1}^{4N^2}$ of $[-\pi, \pi)\times[-\pi, \pi), $ \eqref{eq: identity-cond} can be rewritten as
\begin{align}
\V{A}\, \overlinespace{\widetilde{\mathbf{m}_0}}= \mathbf{1}_{N^2}, \label{eq: m0-A}%\\ 
%\widetilde{\mathbf{m}_0} = [\widetilde{m_0}(\V{\omega}_i)]_{\,i=1,\hdots,4N^2}&\in\mathbb{C}^{4N^2} \notag
\end{align}
where $\widetilde{\mathbf{m}_0} = [\widetilde{m_0}(\V{\omega}_i)]_{\,i=1}^{4N^2}$ and $\V{A}\in \mathbb{C}^{N^2\times 4N^2}$ is a sparse matrix with entries 
$$\V{A}_{i,j} = m_0(\V{\omega}_j)\sum_{k=0}^3\delta(\V{\omega}_j-\V{\omega}_i-\V{\pi}_{2k}), \quad \V{\omega}_j \in [-\pi,0)\times[-\pi,0).$$

%We optimize the regularity of $\sbarm{0}m_0(\V{\omega})$ rather than $\sbarm{0}$; its regularity controls the regularity of $(\,m_0'(\V{\omega}), \widetilde{m_0}'(\V{\omega})\,)$ that can be achieved by a judicious application of Proposition \ref{prop: mc}.
Note that $m_0(\V{\omega})$ in $\V{A}$ here has been regularized by $c(\V{\omega})$, hence we expect the corresponding $\m{0}$ that satisfies \eqref{eq: identity-cond} (or equivalently \eqref{eq: m0-A} on the grid $\mathcal{G}$) to be regular as well. To optimize the regularity of $\m{0}$,  we choose the squared $\ell_2$ norm of the gradient of $\m{0}$ as the objective function, although other forms of regularity may be imposed by different objective functions.

We thus solve the following quadratic minimization problem with linear constraint,
%On a $2N\times 2N$ grid $\G$ of $S_0 = [-\pi, \pi)\times[-\pi, \pi)$, s.t. $\forall \V{\omega}_j \in \G, \; \V{\omega}_j+\V{\nu}_1,\,\V{\omega}_j+\V{\nu}_2,\,\V{\omega}_j+\V{\nu}_3 \in \G$, \eqref{eq: m0-sol} is reformulated as
\begin{align}\label{eq: opt}
\min_{\xvec}\; \Vert \V{D}\xvec\Vert^2,\quad 
%\min_{\xvec\in\mathbb{C}^{4N^2}}\; \Vert \V{D}(\mathbf{m}_0\circ\xvec)\Vert^2,\quad 
s.t. \; \V{A}\xvec = \mathbf{1},
\end{align}
where $\V{D}$ is the gradient operator, $\circ$ is the Hadamard product and $\V{A}$ is the linear operator from \eqref{eq: identity-cond}. 
%Because the set $\{\V{\omega},\, \V{\omega}+\V{\nu}_k,k=1,2,3\}$ is invariant under the shift $\V{\nu}_i,\, i = 1,2,3,$ the rows of $\V{A}$ corresponding to $\V{\omega}$ and $\V{\omega}+\V{\nu}_i$ are identical and 
%Because we only need to consider \eqref{eq: identity-cond} on $[-\pi,0)\times[-\pi,0)$, 
%rows corresponds to $\V{\omega}\in [-\pi,\pi)\times[-\pi,\pi)/\{\V{0},\,\V{\nu}_i,i=1,2,3\}$. Therefore, after removing the duplicate rows, 
%$\V{A}\in \mathbb{C}^{N^2\times 4N^2}$ and \eqref{eq: m0-A} is under-determinant. \\
%We thus use \eqref{eq: m0-A} as a linear constraint in quadratic optimization to solve $\mathbf{\widetilde{m}_0}$. Suppose that $\m{0}$ is smooth, then we build a differential operator $\V{D}$ and solve the following minimization problem:
%or
%\begin{align}
%&\min_{\mvec{0}}\; \Vert \V{D}\mvec{0}\Vert^2 + \lambda \Vert \V{A}\mvec{0} - \mathbf{1}\Vert^2 \label{eq: m0-smooth-relaxed}
%\end{align}
%The solution of \eqref{eq: m0-smooth-relaxed} is $\mvec{0} = \lambda(\lambda \V{A}^\top \V{A} + \V{D}^\top \V{D})^{-1}\V{A}^\top\mathbf{1}$.

%Or suppose $\m{0}$ decays away from the origin, then we build a diagonal weighting operator $\V{W}$, and solve the following minimization problem:
%\begin{align}\label{eq: opt-weight}
%&\min_{\mvec{0}}\; \Vert \V{W}\mvec{0}\Vert^2,\quad s.t. \; \V{A}\mvec{0} = \mathbf{1}
%\end{align}
Supplementary numerical results on solving $\m{0}$ by optimization are provided in Appendix \ref{app: supp-numerical}, where we test this optimization method on known biorthogonal filters $m_0$ and $\widetilde{m_0}$ and compare the solution from the optimization with the ground truth.

Finally, we plug $m_0(\V{\omega})$ and $\m{0}$ into $\V{b}_0'(\V{\omega})$ on the right of \eqref{eq: mj-eq} and solve the linear system for the $m_j$, which has a guaranteed unique solution.

To sum up, we propose \ref{alg} for biorthogonal directional filter construction with dyadic quincunx downsampling scheme.

%According to Proposition \ref{prop: mc}, we can first solve $\mc{0}$ and $m_0^C(\V{\omega})$ and then construct $c(\V{\omega})$ for optimal $\m{0}$ and $m_0(\V{\omega})$. 
%In particular, $m_0^C$ can be computed without knowing $\widehat{k_{\V{\omega}}}$,
%\begin{align}\label{eq: m0C}
%m_0^C(\V{\omega}) = m_0(\V{\omega})|C_{\V{\omega}}| = |det(\M_{1,1}[\widehat{k_{\V{\omega}}},:])| = \prod_{i=1}^6\sigma_i(\M_{1,1}[\widehat{k_{\V{\omega}}},:]) = \prod_{i=1}^6\sigma_i(\M_{1,1}).
%\end{align}
%In practice, we first perform QR decomposition on $\Msub:=\M_{1,1}$ and then take the absolute value of the product of the diagonal entries of the upper triangular matrix, $diag(R)$. 
\vspace{.5em}
\noindent\begin{minipage}{\linewidth}
\rule{\textwidth}{.5pt}
\vspace*{-2em}
\begin{description}% prevent items from splitting
\item[\mylabel{alg}{Algorithm 1}. Construction of $m_0,\widetilde{m_0}$ and $\widetilde{m_j}$ in biorthogonal basis]\
\begin{itemize}
\item[Input:] $\m{j},\,j=1,...,6$, a $2N\times 2N$ regular grid $\G=\{\V{\omega}_i\}_{i=1}^{4N^2}$ over $[-\pi, \pi)\times[-\pi, \pi), $
\item[step 1.] construct $\M[:,\widehat{0}](\V{\omega})$ on the subgrid $[-\pi,0)\times[-\pi,0)$ and check rank constraints \eqref{eq: Mranks},%compute $m_0^C(\V{\omega}) = \left|det(\M_{1,1}[\widehat{k_{\V{\omega}}},:])\right|$
\item[step 2.] solve quadruple $\big( m_0(\V{\omega}),m_0(\V{\omega}+\V{\pi}_2), m_0(\V{\omega} +\V{\pi}_4), m_0(\V{\omega}+\V{\pi}_6) \big)$ using \eqref{eq: m0-even-null} on the subgrid in $[-\pi,0)\times[-\pi,0)$,%compute $\mc{0}$, such that \eqref{eq: LS-new} is solvable and \eqref{eq: identity-cond} holds
\item[step 3.] choose appropriate $\pi$-periodic $c(\V{\omega})$ and replace $m_0(\V{\omega})$ by $m_0'(\V{\omega}) = c(\V{\omega})m_0(\V{\omega})$,% and $\widetilde{m_0}'(\V{\omega})=\m{0}\overline{c}(\V{\omega})^{-1}$,
%solve $m_i(\V{\omega}),\, i=1,...,6$ according to \eqref{eq: LS-new}
\item[step 4.] solve the optimization \eqref{eq: opt} for $\m{0}$ on $[-\pi,\pi)\times[-\pi,\pi)$,
\item[step 5.] solve the reduced linear system \eqref{eq: mj-eq} for $m_j(\V{\omega}),\, j = 1,\cdots,6$.%design $c(\V{\omega})$ and let $m_0(\V{\omega}) = m_0^C(\V{\omega})c(\V{\omega}),\,\m{0} = \mc{0}\overline{c}(\V{\omega})^{-1}$
\end{itemize}
\end{description}
\vspace*{-1em}
\rule{\textwidth}{.5pt}
\end{minipage}\\[.5em]

\noindent{\it Remarks.}
\begin{itemize}
\item[1.] Since $\widetilde{m_j},\, j=1,\cdots,6$ are pre-designed, it is relatively easy to control their regularity. In addition, the regularity of $\widetilde{m_0}$ is optimized by \eqref{eq: opt}. Therefore, according to \eqref{eq: mj_dual}, we may hope to obtain dual wavelets with good regularity.
%\item[1.] The order of step 4. and step 5. in \ref{alg} can be reversed because \eqref{eq: mj-eq} only depends on $m_0(\V{\omega})\overline{\widetilde{m_0}}(\V{\omega} + \V{\pi}_i)$ when $i$ is even, which remains the same for $m_0'(\V{\omega}),\widetilde{m_0}'(\V{\omega})$ regardless of $c(\V{\omega})$ by construction.
\item[2.] In principle, one could formulate an optimization for $c(\V{\omega})$ in step 3. and $\m{0}$ in step 4. jointly in order to obtain optimal smoothness for $\m{0}$ given $m_0(\V{\omega})$ solved in step 2. Instead of solving a linearly constrained quadratic program like \eqref{eq: opt}, one solves a {\it quadratically} constrained quadratic program (QCQP), which is non-convex and in general NP-hard. Such a QCQP can be relaxed to a convex semidefinite program (SDP) that can be efficiently solved although the solution is not exact. See Appendix \ref{app: QCQP} for more details. In Section \ref{sec: numerics}, we discuss how to choose $c(\V{\omega})$ for an $m_0(\V{\omega})$ solved from a specific set of input $\widetilde{m_j}$.
\item[3.] Once can also manipulate pairs of $(m_j,\widetilde{m_j})$ according to the generalization of Proposition \ref{prop: mc} below.
\end{itemize}
\begin{proposition}\label{prop: mjc}
If $\m{j}, m_j(\V{\omega}),\,  j = 0,1,...,6$ satisfy \eqref{eq: LS-new} and \eqref{eq: identity-cond}, 
$m_j^c(\V{\omega})\doteq m_j(\V{\omega})c_j(\V{\omega}), \widetilde{m_j}^c(\V{\omega})\doteq \m{j}\overlinespace{c_j}(\V{\omega})^{-1} \, j = 0,\cdots,6$ 
satisfy \eqref{eq: LS-new} and \eqref{eq: identity-cond} if $ c_0(\V{\omega}) = c_0(\V{\omega}+\V{\pi}_{2k})\,,\forall k = 0,\cdots,3$ and $c_j(\V{\omega}) = c_j(\V{\omega} + \V{\pi}_k),\, \forall \;k = 0,\cdots,7,\, j = 1,\cdots,6$.
\end{proposition}

\begin{comment}
\subsection{solving $m_i$}
In the final step, we substitute $\mc{0}$ and $m_0^C(\V{\omega})$ into \eqref{eq: LS-new} and rewrite it into the following linear system,
\begin{align}\label{eq: mi}
\overline{\M}[:,2:7]\,\mathbf{m}[2:7](\V{\omega}) = 
\begin{bmatrix}
1-m_0^C\overline{\widetilde{m_0}^C}(\V{\omega})\\
0\\
-m_0^C\overline{\widetilde{m_0}^C}(\V{\omega}+\V{\pi}_2)\\
\vdots \\
0
\end{bmatrix}
=:\mathbf{b}(\V{\omega}).
\end{align}
The solution of \eqref{eq: mi} depends only on $m_0^C\overline{\widetilde{m_0}^C}$, or equivalently $m_0\overline{\widetilde{m_0}}$. 
\end{comment}

\section{Numerical Experiments}\label{sec: numerics}

In this section, we demonstrate the numerical construction of biorthogonal directional wavelets on a quincunx lattice using our proposed \ref{alg} implemented in Matlab.

\begin{figure}
\centering
\includegraphics[width=.8\textwidth]{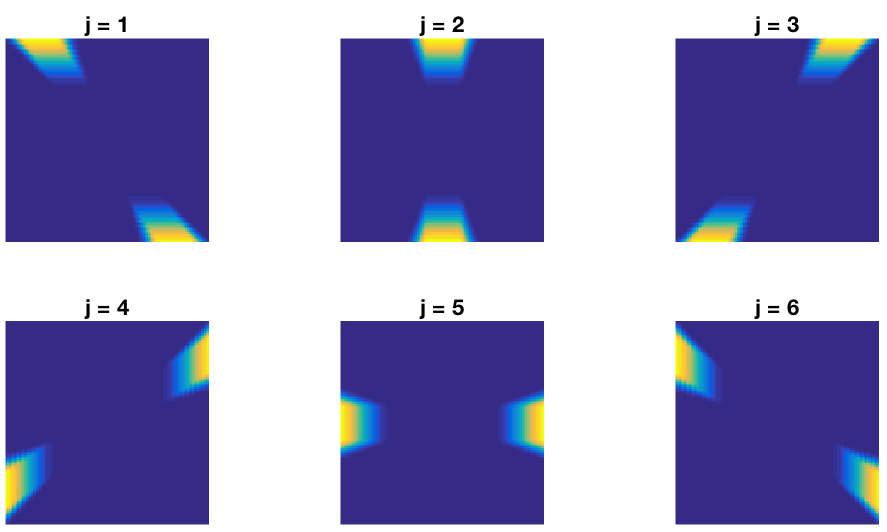}
\caption{ Input $|\m{j}|$ constructed in the same way as shearlets.}
\label{fig: mjdual}
\end{figure}

For the input of \ref{alg}, we use $\widetilde{m_j}$ in the form of \eqref{eq: m-form}, with phases in \eqref{eq: phase-sol} and amplitudes $|\widetilde{m_j}|$ shown in Figure \ref{fig: mjdual} constructed as follows. We start with a symmetric $|\widetilde{m_2}|$, then compute $|\widetilde{m_1}|$ and $|\widetilde{m_3}|$ by shearing $|\widetilde{m_2}|$ counter-clockwise and clockwise respectively. $|\widetilde{m_4}|, |\widetilde{m_5}|$ and $|\widetilde{m_6}|$ are obtained by symmetry with respect to the diagonal. This is the same approach used in the shearlet construction in \cite{kutyniok2012digital}. Furthermore, we set $\m{j} = 0,\,\forall \V{\omega}\in C_0 = [-\pi/2,\pi/2)\times [-\pi/2, \pi/2)$ and according to Theorem \ref{thm: thm}, we enforce $|\widetilde{m_1}(\dfrac{\pi}{2},\dfrac{\pi}{2})|\neq 0$ and $|\widetilde{m_6}(\dfrac{\pi}{2},\dfrac{\pi}{2})|\neq 0$. As the first step, we numerically verify that this particular design of $\widetilde{m_j}$ satisfies the rank constraints \eqref{eq: Mranks}.\footnote{In practice, we find it hard for $\widetilde{m_j}$ to satisfy the rank constraint \eqref{eq: Mranks} without enforcing $\widetilde{m_j}$ to be zero on $C_0$. This may indicate topological obstruction in our biorthogonal scheme}

\begin{figure}
\centering
\includegraphics[width=.9\textwidth]{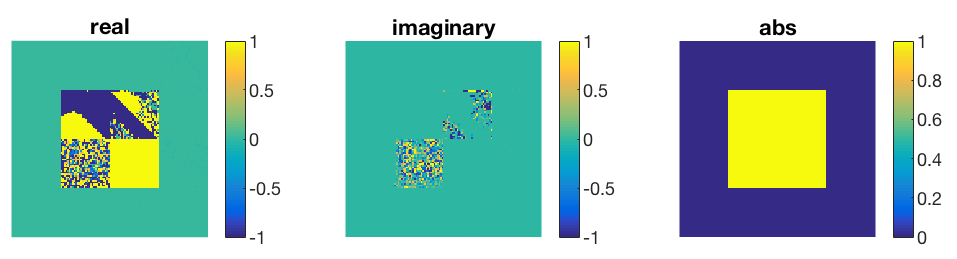}
\caption{ $m_0(\V{\omega})$ constructed from $\widetilde{m_j}$. Left to right: $Re(m_0(\V{\omega})),\, Im(m_0(\V{\omega}))$ and $|m_0(\V{\omega})|$.}
\label{fig: m_0}
\end{figure}

We proceed to solve $m_0(\V{\omega})$ in quadruple separately for each $\V{\omega}$ in $[-\pi,0)\times [-\pi,0)$. As pointed out earlier, these solutions still have an unconstrained degree of freedom in the form of a constant $a_{\V{\omega}}$; the result is shown in Figure \ref{fig: m_0} for one implementation using Matlab solvers. This solution $m_0(\V{\omega})$ has both inherent irregularity of the biorthogonal construction from the input and artificial irregularity from the algorithm: the amplitude $|m_0(\V{\omega})|$ is supported on $C_0$, where $|m_0(\V{\omega})| = 1$ and its discontinuity at $\partial C_0$ corresponds to that of the input $\m{j}$; however, the phase of $m_0(\V{\omega})$ is discontinuous even on the interior of $C_0$ due to $a_{\V{\omega}}$, an artificial irregularity we remove in the next step by introducing $c(\V{\omega})$.

%We then solve \eqref{eq: opt} for $\m{0}$, and the solution is coupled with $m_0(\V{\omega})$ in the previous step. The left figure in Figure \ref{fig: m_0_m0dual} shows that $\m{0}$ is also supported on $C_0$ as $m_0(\V{\omega})$. Minimizing the objective function in \eqref{eq: opt} optimizes the smoothness of $m_0(\V{\omega})\sbarm{0}$; this quantity is already $\equiv 0$ outside $C_0$ since $m_0(\V{\omega})=0,\,\forall\V{\omega}\not\in C_0$. It follows that the behavior of the optimizing $\widetilde{m_0}$ outside $C_0$ is mostly determined by the linear constraint in \eqref{eq: opt}-- it is thus essentially the identity condition \eqref{eq: identity-cond} that causes $\m{0}$ to vanish outside $C_0$.%This is a result of the linear constraint in \eqref{eq: opt} from the identity condition \eqref{eq: identity-cond}, not because of the regularization by the objective function where $m_0(\V{\omega})\sbarm{0}\equiv 0,\, \forall\V{\omega}\not\in C_0$. 
%\textcolor{red}{Furthermore, the middle and right figures in Figure \ref{fig: m_0_m0dual} suggest that $m_0(\V{\omega})\sbarm{0}=\V{1}_{C_0}$ numerically, as they should be.}

\begin{figure}
\centering
\includegraphics[width=.5\textwidth]{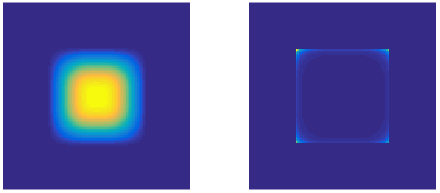}
\caption{Left: $\widetilde{m_0}$, designed smooth function supported on the central square $C_0$, right: $m_0'$, where %$ m_0'\overlinespace{\widetilde{m_0}}(\V{\omega}) = 
$\overline{m_0'}\m{0} = \mathbbm{1}_{C_0}(\V{\omega})$.} 
\label{fig: smooth_m0dual}
\end{figure}

To regularize $m_0(\V{\omega})$, we multiply it by an appropriate $\pi$-periodic $c(\V{\omega})$. In particular, we can first construct $c(\V{\omega})$ on $C_0$ freely and then extend it to $S_0$ by its $\pi$-periodicity in both $\omega_1$ and $\omega_2$. It turns out that in this specific numerical example we consider here, we can explicitly design the regularized $m_0(\V{\omega})$ ($m_0'(\V{\omega})$) and the corresponding $\m{0}$. Since $m_0$ is only supported on $C_0$, $m_0'(\V{\omega}) = m_0(\V{\omega})c(\V{\omega})$ is determined by the value of $c(\V{\omega})$ on $C_0$. Therefore, $m_0'(\V{\omega})$ can be any continuous function on $C_0$. On the other hand, $m_0'\sbarm{0}\equiv 0,\,\forall\V{\omega}\not\in C_0$, and \eqref{eq: identity-cond} (correspondingly the linear constraint \eqref{eq: m0-A}) reduces to $m_0'\sbarm{0} = 1,\, \forall\V{\omega}\in C_0$. In other words, $\m{0}$ is uniquely determined on $C_0$ by $m_0'(\V{\omega})$ or vice versa. Because we want $\m{0}$ to be smooth and has fast decay from the origin such that the corresponding dual wavelets $\widetilde{\psi^j}$ have good spatial locality, we can actually first design $\m{0}$ on $S_0$ and then construct $m_0'(\V{\omega}) = \m{0}^{-1}$ on $C_0$.
In particular, we let $\widetilde{m_0}'$ be the low pass filter of a 2D tensor wavelets, see Figure \ref{fig: smooth_m0dual}.

\begin{figure}[h]
\centering
\includegraphics[width = .5\textwidth]{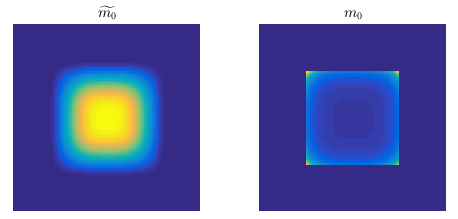}
\caption{ Left: $\widetilde{m_0}$, with support outside $C_0$, right: $m_0'$, where $ \overline{m_0'}\m{0}  = \mathbbm{1}_{C_0}(\V{\omega})$.}
\label{fig: smooth_m0dual-slow}
\end{figure}

\begin{figure}
\includegraphics[width = \textwidth]{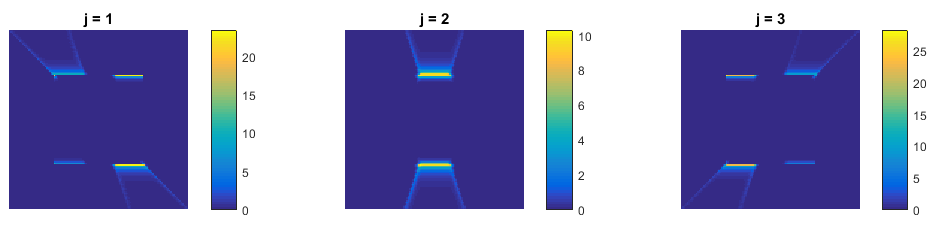}
\caption{$|m_j(\V{\omega})|$ solved from \eqref{eq: mj-eq} given $\widetilde{m_j}$ in Figure \ref{fig: mjdual}, $m_0'$ and $\widetilde{m_0}$ in Figure \ref{fig: smooth_m0dual-slow}}
\label{fig: m_j-slow}
\end{figure}

\noindent{\it Remarks.} 
\begin{itemize}
\item[1.]If we use the above $m_0'$ derived from a known $\m{0}$ and solve \eqref{eq: opt} for $\m{0}$ as in step 4. of \ref{alg}, we obtain a solution $\widetilde{m_0}'(\V{\omega})$ not exact the same but close to the known $\m{0}$. Moreover, we numerically verify that $m_0'(\V{\omega})\overlinespace{\widetilde{m_0}'}(\V{\omega})=\mathbbm{1}_{C_0}$ as they should be.
\item[2.]There is no restriction on the support of $\m{0}$ as long as \eqref{eq: identity-cond} is satisfied. Although a slower decay of $\m{0}$ on $S_0$ increases the regularity $m_0'(\V{\omega})$ on $C_0$, see Figure \ref{fig: smooth_m0dual-slow}, the resulting $m_j$ solved in the final step do not have ideal direction selectivity, see Figure \ref{fig: m_j-slow}.
\end{itemize}

\begin{figure}
\centering
\includegraphics[width = .9\textwidth]{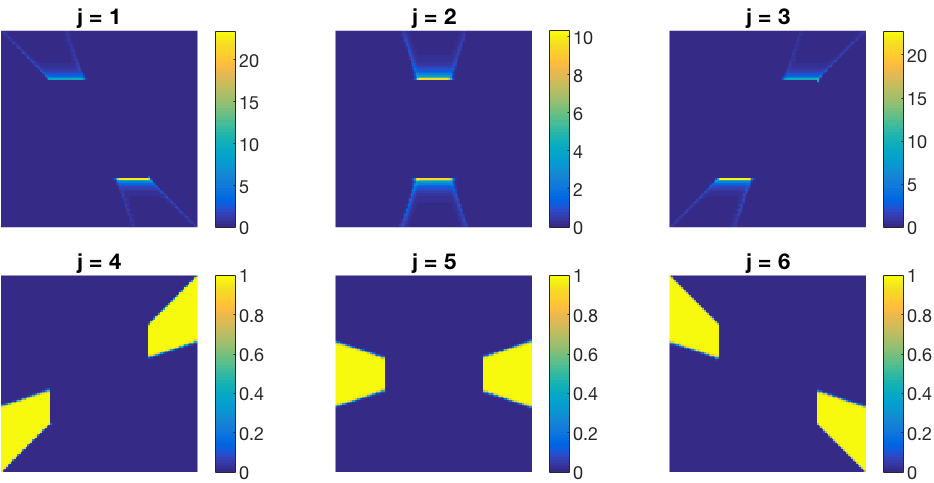}
\caption{Top: $|m_j(\V{\omega})|,\, j= 1,2,3$,
Bottom: $|m_j(\V{\omega})\m{j}|,\, j=4,5,6$, where $m_j(\V{\omega})$ is solved from \eqref{eq: mj-eq} given $\widetilde{m_j}$ in Figure \ref{fig: mjdual}, $m_0'$ and $\widetilde{m_0}$ in Figure \ref{fig: smooth_m0dual}.}
\label{fig: m_j}
\end{figure}

Finally, we solve \eqref{eq: mj-eq} for $m_j$. As shown in the top row Figure \ref{fig: m_j}, the energy of $m_j$ concentrates at $\partial C_0$, where $m_j$ decay to near zero. Moreover, the bottom row of Figure \ref{fig: m_j} shows that $|m_j\m{j}|$ are close to constant on $C_j$. Such irregularity roots in the irregularity of biorthogonal bases construction we show in Section \ref{subsec: discontinuity}, which prevents input $\widetilde{m_j}$ to be continuous in the first place.
 We also numerically verify that $m_j(\V{\omega})$ and $\m{j}$ have the same phase, i.e. $m_j\sbarm{j}\in \mathbb{R}$.

\begin{figure}
\centering
\includegraphics[width=.7\textwidth]{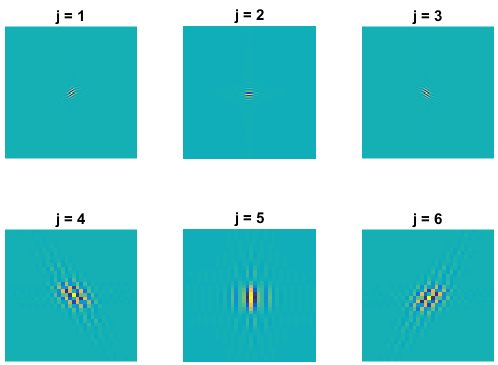}
\caption{ Real part of $\widetilde{\psi^j}$ constructed from $\widetilde{m_j},\, j=1,\cdots, 6$ in Figure \ref{fig: mjdual} and $\widetilde{m_0}'$ in Figure \ref{fig: smooth_m0dual} using \eqref{eq: mj_dual}. Top: $\widetilde{\psi^j}$ without scaling, bottom: $\widetilde{\psi^j}$ with eight time zoom-in } 
\label{fig: wavelets}
\end{figure}

So far, we construct a set of $(m_j,\widetilde{m_j})_{j = 0,\cdots,6}$ that satisfies \eqref{eq: LS-new} and \eqref{eq: identity-cond}, thus it can be used to construct biorthogonal wavelets based on \eqref{eq: mj} and \eqref{eq: mj_dual}. Figure \ref{fig: wavelets} shows the dual wavelets  $\widetilde{\psi^j}$ in \eqref{eq: bi-orth MRA} constructed using \eqref{eq: mj_dual}. Because of the regularity we impose on $\widetilde{m_j}$ and $\widetilde{m_0}'$, the dual wavelets are spatially localized and have good direction selection. The wavelets and scaling functions in \eqref{eq: bi-orth MRA} can be constructed using \eqref{eq: mj} similarly, but with much poorer regularity originated in $m_j$ and $m_0'$.

Although using a different set of $\widetilde{m_j}$ as input paired with a carefully tweaked $\widetilde{m_0}'$ might improve the regularity of the dual wavelets $\widetilde{\psi^j}$, the intrinsic irregularity of the corresponding wavelets $\psi^j$ shall remain.

\begin{comment}
\begin{figure}
\centering
\includegraphics[width=\textwidth]{m0_m0dual_new.png}
\caption{$|\m{0}|$ and $m_0\sbarm{0}$, where $\m{0}$ is solved by optimization \eqref{eq: opt}, given $\widetilde{m_j}$ in Figure \ref{fig: mjdual} and $m_0$ in Figure \ref{fig: m_0}. Left to right: $|\m{0}|$, $Re(m_0\sbarm{0})$ and $Im(m_0\sbarm{0})$. }
\label{fig: m_0_m0dual}
\end{figure}
\end{comment}

\begin{comment}
\begin{figure}
\centering
\includegraphics[width=\textwidth]{m.png}
\caption{ $|m_j(\V{\omega})|$, where $m_j(\V{\omega})$ is solved from \eqref{eq: mj-eq} given $\widetilde{m_j}$ in Figure \ref{fig: mjdual}, $m_0'$ and $\widetilde{m_0}'$ in Figure \ref{fig: smooth_m0dual}. }
\label{fig: m_j}
\end{figure}

\begin{figure}
\centering
\includegraphics[width=\textwidth]{m_mdual.png}
\caption{ $|m_j(\V{\omega})\m{j}|$ for $m_j$ in Figure \ref{fig: m_j} and $\widetilde{m_j}$ in Figure \ref{fig: mjdual}. }
\label{fig: m_j_mjdual}
\end{figure}
\end{comment}

%% file: conclusion.tex
\section{Conclusion and future work}\label{sec: end}
In this paper, we consider directional wavelet schemes on a dyadic quincunx sub-lattice and analyze their regularity. We show that filters in bi-orthogonal bases have the same discontinuity in the frequency domain as the orthonormal bases at the corners of $C_0 = [-\pi/2,\pi/2)\times[-\pi/2,\pi/2)$. 

%Our analysis is closely related to our proposed bases construction algorithms, 
%and we show that the construction method of orthonormal bases can be easily extended to build frames construction of redundancy 2, which achieve much better time frequency localization and thus practically useful.
 We propose a different approach to construct biorthogonal wavelets from our previous approach for the orthonormal bases construction \cite{yin2014orthshear}. The directional dual filters $\widetilde{m_j}$ are first designed such that they can be extended to a bi-orthogonal frame and the remaining filters are obtained by solving linear systems and a constrained quadratic optimization derived from the identity summation and shift cancellation conditions for a biorthogonal MRA. We show numerically that regularized dual wavelets $\widetilde{\psi^j}$ can be constructed, yet their corresponding wavelets $\psi^j$ are still discontinuous in frequency domain, which is unavoidable according to our analysis.

We have looked at extensions of orthonormal bases in two different directions: tight frames (which are self-dual but redundant) with low redundancy and bi-orthogonal bases (which remain non-redundant but are no longer self-dual). In both cases we can gain some regularity. The extension of the biorthogonal bases to low-redundancy dual frame construction is not studied here, achieve at least the same regularity as low-redundancy tight frames, but with more flexibility in the construction. We leave this further generalization to future work.

%% file: parseval-cond.tex
\section{Proof of Theorem \ref{thm: conds}}\label{app: cond-thm}
Take the Fourier transform of both sides of \eqref{eq: PR}, we have 
\begin{align*}
\sum_{\V{k}}\langle f,\phi_{\V{k}}\rangle\hat{\phi}(\V{\omega})e^{-i\V{\omega}^\top\V{k}} = \sum_{\V{k}}&\langle f,\phi_{1,\V{k}}\rangle e^{-i\V{\omega}^\top\V{Dk}}|\V{D}|^{1/2}\hat{\phi}(\V{D}^T\V{\omega}) \\
&+ \sum_{j=1}^J\sum_{\V{k}}\langle f,\psi^j_{1,\V{k}}\rangle e^{-i\V{\omega}^\top\V{QDk}}|\V{QD}|^{1/2}\hat{\phi}(\V{D}^\top\V{\omega}).
\end{align*}
We use $\sum_{\V{k}}$ for summation over $\mathbb{Z}^2$ without specifying the set $\mathbb{Z}^2$.
Suppose $m_j$ are trigonometric series
\begin{align}\label{eq: mra1}
m_0(\V{\omega}) = \sum_{\V{k}} c_{\V{k}}e^{-i \V{\omega}^\top\V{k}} \quad
m_j(\V{\omega}) = \sum_{\V{k}} g_{\V{k}}e^{-i \V{\omega}^\top\V{k}},\quad j=1,\cdots,J.
\end{align}
The first term on the right hand side can be represented by $\hat{\phi}(\V{\omega})$ and $\langle f,\phi_k\rangle$ using \eqref{eq: m0} and \eqref{eq: mra1}.

\begin{align*}
\text{the first term on R.H.S. } = \sum_{\V{k}}\langle f,\phi_{1,\V{k}}\rangle e^{-i\V{\omega}^\top\V{Dk}}|\V{D}|^{1/2}m_0(\V{\omega})\hat{\phi}(\V{\omega}) \\= \sum_{\V{k}}\Big(\sum_{\V{k}'}\langle f,\phi_{\V{k}'}\rangle\overlinespace{c_{\V{k'-Dk}}}|\V{D}|^{1/2}\Big)e^{-i\V{\omega}^\top\V{Dk}}|\V{D}|^{1/2}m_0(\V{\omega})\hat{\phi}(\V{\omega})\\
=\sum_{\V{k}'}\langle f,\phi_{\V{k}'}\rangle\Big(|\V{D}|\sum_{\V{k}}\overlinespace{c_{\V{k'-Dk}}}e^{i\V{\omega}^\top(\V{k'-Dk})}\Big)e^{-i\V{\omega} ^\top\V{k}'} m_0(\V{\omega})\hat{\phi}(\V{\omega}).
\end{align*}

Let $\{\V{\beta}\}\doteq \V{D}\mathbb{Z}^2+\V{\beta}$ for $\V{\beta}\in B, \,s.t.\,\bigcup_{\V{\beta}\in B} \{\V{\beta}\} = \mathbb{Z}^2$.\footnote{The choice of $B$ is not unique and one choice is $ \{ (0,0),\,(1,0),\,(0,1),\,(1,1)\}$.}  The sum over $\mathbb{Z}^2$ can then be written as a double sum $\sum_{\V{\beta}\in B}\sum_{\V{k}'\in \{\V{\beta}\}}$,
\begin{align*}
\sum_{\V{\beta}\in B}\sum_{\V{k}'\in\{\V{\beta}\}} \langle f,\phi_{\V{k}'}\rangle\sum_{\V{k}}\overlinespace{c_{\V{k'-Dk}}}e^{i\V{\omega}^\top(\V{k}'-\V{Dk})}e^{-i\V{\omega}^\top\V{k}'}|\V{D}|m_0(\V{\omega})\hat{\phi}(\V{\omega})\\
=\sum_{\V{\beta}\in B}\sum_{\V{k}'\in\{\V{\beta}\}} \langle f,\phi_{\V{k}'}\rangle\Big(\,\sum_{\V{k}\in\{\V{\beta}\}}\overlinespace{c_{\V{k}}}e^{i\V{\omega}^\top\V{k}}\,\Big)e^{-i\V{\omega}^\top\V{k}'}|\V{D}|m_0(\V{\omega})\hat{\phi}(\V{\omega}).
\end{align*}
Due to the identity $\sum_{\V{\pi}\in\Gamma_0}e^{i\V{\beta}^\top\V{\pi}} = |\Gamma_0| \, \scalebox{1.3}{$\chi$}_{\V{D}\mathbb{Z}^2}(\V{\beta})$, the sum $\sum_{\V{k}\in \{\V{\beta}\}}c_{\V{k}}e^{-i\V{\omega}^\top\V{k}}$ equals to a linear combination of  $m_0$ with shifts in $\Gamma_0$,
%, similar to the trigonometric expansion of $m_0$ in \eqref{eq: mra1}, but $\V{k}$ takes value on the shifted sub-lattice $\{\V{\beta}\}$ instead of $\mathbb{Z}^2$, equals to instead a linear combination of $m_0$ with shifts $\Gamma_0$,
\begin{align}\label{eq:eq1}
\sum_{\V{k}\in \{\V{\beta}\}}c_{\V{k}}e^{-i\V{\omega}^\top\V{k}}
= \frac{1}{|\Gamma_0|}\;\sum_{\V{\pi}\in\Gamma_0}m_0(\V{\omega}+\V{\pi})\;e^{i\V{\beta}^\top\V{\pi}} .
\end{align}
Substitute \eqref{eq:eq1} into the previous expression and notice $|\Gamma_0| = |D|=4$, we have
\begin{align*}
\sum_{\V{\beta}\in B}\sum_{\V{k}'\in \{\V{\beta}\}}\langle f,\phi_{\V{k}'}\rangle\sum_{\V{\pi}\in\Gamma_0}\overlinespace{m_0}(\V{\omega}+\V{\pi})\;e^{-i\V{\beta}^\top\V{\pi}}\,e^{-i\V{\omega}^\top\V{k}'}m_0(\V{\omega})\hat{\phi}(\V{\omega}).
\end{align*}
Since $e^{i \V{\pi}^\top\V{\beta}}=e^{i\V{\pi}^\top\V{k}'},\; \forall \V{k}'\in \{\V{\beta}\} $, we can rewrite the double sum $\sum_{\V{\beta}\in B}\sum_{\V{k}'\in \{\V{\beta}\}} $  back to a unit sum over $\mathbb{Z}^2$ as follows.
\begin{align*}
\sum_{\V{k}'}\langle f,\phi_{\V{k}'}\rangle e^{-i\V{\omega} ^\top\V{k}'}\hat{\phi}(\V{\omega})\Big(\sum_{\V{\pi}\in\Gamma_0}\overlinespace{m_0}(\V{\omega}+\V{\pi})m_0(\V{\omega})e^{-i\V{\pi}^\top\V{k}'} \Big).
\end{align*}

Similarly, the second term on the R.H.S. of \eqref{eq: PR} equals to 
\begin{align*}
\sum_{j=1}^J\sum_{\V{k}'}\langle f,\phi_{\V{k}'}\rangle e^{-i\V{\omega}^\top \V{k}'}\hat{\phi}(\V{\omega})\Big(\sum_{\V{\pi}\in\Gamma_1} \overlinespace{m_j}(\V{\omega}+\V{\pi})m_j(\V{\omega})e^{-i\V{\pi}^\top\V{k}'} \Big)
\end{align*}
based on the following equality analogous to \eqref{eq:eq1}
\begin{align}
\sum_{\V{k}\in\{\V{\alpha}\}} g_{\V{k}'} e^{-i\V{\omega}^\top\V{k}} = \frac{1}{|\Gamma_1|}\, \sum_{\V{\pi}\in\Gamma_1} m_j(\V{\omega} + \V{\pi})e^{i\V{\alpha}^\top\V{\pi}},
\end{align}
where $\{\V{\alpha}\} \doteq \V{QD}\mathbb{Z}^2 + \V{\alpha}$ for $\V{\alpha}\in A,\, s.t.\, \bigcup_{\V{\alpha}\in A} \{\V{\alpha}\} = \mathbb{Z}^2$.
(For Theorem \ref{thm: frame-conds} on frame construction, the summation of shifts $\V{\pi}$ is over $\Gamma_0$ instead of $\Gamma_1$.) 
Combining the two terms on the R.H.S. of \eqref{eq: PR}, and compare the coefficients of $\langle f,\phi_{\V{k}'}\rangle e^{-i\V{\omega}^\top \V{k}'}\hat{\phi}(\V{\omega})$ on both sides, the perfect reconstruction condition is then equivalent to $\forall \V{k}'$,
\begin{align*}
\sum_{\V{\pi}\in\Gamma_0}e^{-i\V{\pi}^\top\V{k}'}\overlinespace{m_0}(\V{\omega}+\V{\pi})m_0(\V{\omega}) + \sum_j\sum_{\V{\pi}\in\Gamma_1} e^{-i\V{\pi}^\top\V{k}'}\overlinespace{m_j}(\V{\omega}+\V{\pi})m_j(\V{\omega}) = 1. 
%\sum_{l=0}^{3}e^{-i\gamma_l^\top(k'-k_0)}\overline{M_0(\xi+\gamma_l)}M_0(\xi) + \sum_j\sum_{s=0}^7 e^{-i\nu_s^\top(k'-k_j)}\overline{M_j(\xi+\nu_s)}M_j(\xi) = 1. 
\end{align*} 
This is equivalent to 
\begin{align*}
&|m_0(\V{\omega})|^2 + \sum_j|m_j(\V{\omega})|^2 = 1
\end{align*}
and
\begin{align*}
\sum_{j=0}^J\overlinespace{m_j}(\V{\omega}+\V{\pi})m_j(\V{\omega}) = 0, 
%+ \overline{m_0(\V{\omega}+\V{\pi})}m_0(\V{\omega}) = 0, 
\,\V{\pi}\in \Gamma_0\setminus\{\V{0}\}\\
\sum_{j=1}^J\overlinespace{m_j}(\V{\omega}+\V{\pi})m_j(\V{\omega}) = 0,\,\V{\pi}\in \Gamma_1\setminus \Gamma_0
\end{align*}
\qed

\noindent{\it Remark}.
If we have a shift $\V{k}_0$ in the down-sample scheme for $\phi_1$, i.e. $\V{D}\mathbb{Z}^2 - \V{k}_0$ instead of $\V{D}\mathbb{Z}^2$, so that we obtain coefficient of $\tilde{\phi}_{1,\V{k}} = \phi_{1,\V{k}+\V{k}_0}$ instead of $\phi_{1,\V{k}}$, and $\tilde{\phi}_1(\V{x}) =\phi_1(\V{x}-\V{k}_0)= |\V{D}|^{1/2}\sum_{\V{k}}c_{\V{k}}\phi(\V{x-k-k}_0) = |\V{D}|^{1/2}\sum_{\V{k}}c_{\V{k}-\V{k}_0}\phi(\V{x-k})$. This change of down-sample scheme results in an extra phase term $e^{-i\V{\omega}^\top \V{k}_0}$ in $m_0$. 
Similarly, if we downsample $\psi_1^j$ on a shifted sub-lattice $\V{QD}\mathbb{Z}^2-\V{k}_j$, we then have an extra phase $e^{i\V{\pi}^\top\V{k}_j}$ before $\overlinespace{m_j}(\V{\omega}+\V{\pi})m_j(\V{\omega})$ in shift cancellation condition. This provides additional freedom in the construction yet it is not substantial. Here, we use the down-sample scheme without translation.

%% file: lemmas.tex
\pdfoutput=1
\section{Proof of lemmas and propositions for biorthogonal schemes}\label{app: lemmas}

\subsection{Discontinuity of $\m{j}$}\label{app: discontinuity}
\begin{lemma}\label{lem: subM-singular-sys}
Define $d_{i,j}(\V{\omega}) = \det([\mrow{k_1}(\V{\omega})^\top,\cdots,\mrow{k_6}(\V{\omega})^\top]),\;$ where $0\leq k_1<\cdots<k_6\leq 7,\, s.t.\, k_l\neq i,j.$
\eqref{eq: LS-new} is solvable $\forall \V{\omega}$ if and only if \vspace{.5em}
\begin{align}
\label{eq: singular-cond}
\mathfrak{D}(\V{\omega})\begin{bmatrix}
\sbarm{0}\\
\sbarmp{0}{2}\\
\sbarmp{0}{4}\\
\sbarmp{0}{6}
\end{bmatrix}
\doteq
\begin{bmatrix}
0 & d_{0,2} & d_{0,4} & d_{0,6}\\
-d_{0,2} & 0 & d_{2,4} & d_{2,6}\\
-d_{0,4} & -d_{2,4} & 0 & d_{4,6}\\
-d_{0,6} & -d_{2,6} & -d_{4,6} & 0
\end{bmatrix}
\begin{bmatrix}
\sbarm{0}\\
\sbarmp{0}{2}\\
\sbarmp{0}{4}\\
\sbarmp{0}{6}
\end{bmatrix}
= \begin{bmatrix}
0\\0\\0\\0
\end{bmatrix}.
\end{align}
\end{lemma}
\noindent{\it Proof.}
By Lemma \ref{lem: subM-singular} and Lemma \ref{lem: M-symmetry},  $\M[\widehat{k},:],\,k=0,2,4,6$ are singular,
The singularity condition on  $\M[\widehat{0},:](\V{\omega})$ can be rewritten as follows,
\begin{align}\label{eq: singular-omega}
0 &=\det(\M[\widehat{0},:]) \notag\\
&=  \sbarmp{0}{2}\cdot\det(\Msub[\widehat{2},:])\notag\\
&\quad+ \,\sbarmp{0}{4}\cdot\det(\Msub[\widehat{4},:])
+ \sbarmp{0}{6}\cdot\det(\Msub[\widehat{6},:])\notag\\
&= 0\cdot\sbarm{0}\,+\,d_{0,2}\cdot\sbarmp{0}{2} \notag\\
&\quad+\,d_{0,4}\cdot \sbarmp{0}{4}\,+\, d_{0,6}\cdot\sbarmp{0}{6}
\end{align}
%This is the first equation in the linear system \eqref{eq: singular-cond}. Substitute $\V{\omega}$ by $\V{\omega + \pi_2}$ in \eqref{eq: singular-omega} and use the $2\pi-$periodicity of $\V{\omega}$, we have the singularity condition on $\M[-1,:](\V{\omega+\pi_2})$ as follows
%then the above singularity condition on $\M[-1,:]$ at $\V{\omega}$ can be rewritten as follows,
%\begin{align*}
%[0,\, d_{0,2}(\V{\omega}),\, d_{0,4}(\V{\omega}),\, d_{0,6}(\V{\omega})]\,[\widetilde{m_0}(\V{\omega}),\,\widetilde{m_0}(\V{\omega}+\V{\pi}_2),\, \widetilde{m_0}(\V{\omega}+\V{\pi}_4),\,\widetilde{m_0}(\V{\omega}+\V{\pi}_6)]^\top = 0
%\end{align*}
%It is easy to verify that the above singular condition at $\V{\omega}+\V{\pi}_2$ is equivalent to 
%\begin{align*}
%-d_{0,2}(\V{\omega})\cdot \widetilde{m_0}(\V{\omega}) + d_{2,4}(\V{\omega})\cdot \widetilde{m_0}(\V{\omega}+\V{\pi}_4) + d_{2,6}(\V{\omega})\cdot\widetilde{m_0}(\V{\omega}+\V{\pi}_6) = 0,
%[-d_{0,2}(\V{\omega}),\, 0,\,d_{2,4}(\V{\omega}),\,d_{2,6}][\widetilde{m_0}(\V{\omega}),\,\widetilde{m_0}(\V{\omega}+\V{\pi}_2),\, \widetilde{m_0}(\V{\omega}+\V{\pi}_4),\,\widetilde{m_0}(\V{\omega}+\V{\pi}_6)]^\top = 0,
%\end{align*}
%which is the second linear equation in  \eqref{eq: singular-cond}.
Similarly, the second to fourth equations can be obtained by rewriting the singularity condition on $\M[\widehat{2},:]$, $\M[\widehat{4},:]$ and $\M[\widehat{6},:]$ respectively.\qed
% at $\V{\omega}+\V{\pi}_4$ and $\V{\omega}+\V{\pi}_6$ in the coordinate of $\V{\omega}$.\qed
%where $\mathfrak{D}(\V{\omega})$ is anti-symmetric. Because $\mathfrak{D}(\V{\omega})$ is independent of $m_0(\V{\omega})$, \eqref{eq: singular-cond} holds for $\mc{0}$ as well.

%On the other hand, given $m_0$, $\widetilde{m_0}$ has to satisfy the identity constraint \eqref{eq: identity-cond}.
The identity constraint \eqref{eq: identity-cond} on $m_0$ and the singularity condition \eqref{eq: singular-cond} together imply the following proposition,
%Due to the periodic wrapping of the frequency square $S_0$, we only need to consider \eqref{eq: singular-cond} and \eqref{eq: identity-cond} on $S_1$ and they imply the following proposition,
\begin{proposition}\label{prop: feasibility}
Given $\widetilde{m_i}, i = 1,\cdots,6$, \eqref{eq: identity-cond} has no solution for $\widetilde{m_0}$, if $\exists\,\V{\omega}, \,s.t. \; [m_0(\V{\omega}), m_0(\V{\omega}+\V{\pi}_2),m_0(\V{\omega}+\V{\pi}_4),m_0(\V{\omega}+\V{\pi}_6)]$ is a linear combination of the rows of $\mathfrak{D}(\V{\omega})$.% in \eqref{eq: singular-cond}.
\end{proposition}

\noindent{\bf Proof of Lemma \ref{lem: rank1}}:\\[.2em]
{\bf Lemma \ref{lem: rank1}}. {\it 
If $\V{\omega}\in S_\rho$ s.t. \eqref{eq: identity-cond} holds and $\M[\widehat{0},:](\V{\omega})$ is singular, then  $rank(\mrow{1},\mrow{7}) = 1$ and $rank(\mrow{3},\mrow{5}) = 2$ or $rank(\mrow{3},\mrow{5}) = 1$ and $rank(\mrow{1},\mrow{7}) = 2$.
}\\[1em]
\noindent{\it Proof.}
When $\rho$ is small enough, due to the concentration property, $\m{i}$ is zero on all but a few sets $S_\rho + \V{\pi}_j$ (see Fig.\ref{fig: S-shifts} for reference of $S_\rho$ and its shifts), thus $\mrow{i}(\V{\omega})$ is sparse on $S_\rho$ and $\M[:,\widehat{0}]$ takes the following form
\begin{align}
\label{eq: sparse-mat}
\M[:,\widehat{0}](\V{\omega})=
\begin{bmatrix}
\mrow{0}\\
\mrow{1}\\
\mrow{2}\\
\mrow{3}\\
\mrow{4}\\
\mrow{5}\\
\mrow{6}\\
\mrow{7}
\end{bmatrix}
=
\begin{bmatrix}
0 & 0 & 0 & 0 & 0 & 0\\
* & 0 & 0 & 0 & 0 & *\\
0 & 0 & 0 & * & * & 0\\
0 & 0 & * & * & 0 & 0\\
0 & * & * & 0 & 0 & 0\\
0 & 0 & * & * & 0 & 0\\
* & 0 & * & * & 0 & *\\
%0 & * & * & 0 & 0 & 0\\
%0 & 0 & 0 & * & * & 0\\
%* & 0 & 0 & 0 & 0 & *\\
%* & 0 & 0 & 0 & 0 & *\\
%0 & 0 & * & * & 0 & 0\\
%0 & 0 & * & * & 0 & 0\\
* & 0 & 0 & 0 & 0 & *
\end{bmatrix}
%=\V{P}\,\widetilde{\mathbf{M}}[:,2:7],
\end{align}
where $*$ denote possible non-zero entries.
%where $\V{P}$ is a row permutation matrix. 
We make the following observation of $\mrow{i}$:
\begin{itemize}
\item[(i)] $\mrow{0}$ is a zero vector
\item[(ii)] $\mrow{2}$ and $\mrow{4}$ are linearly independent of each other and the rest of $\mrow{i}$
\item[(iii)] $span\{\mrow{1},\mrow{7}\} \perp span\{\mrow{3},\mrow{5}\}$ and $rank(\mrow{1},\mrow{7}) \leq 2$, \\$rank(\mrow{3},\mrow{5})\leq 2$
\item[(iv)] $span\{\mrow{1}, \mrow{7}, \mrow{3},\mrow{5},\mrow{6}\} \leq 4$
\end{itemize}
Since $m_0(\V{\omega})\neq 0$ on $S_\rho$, \eqref{eq: m0-cramer} then implies that $\det(\Msub[\widehat{k_{\V{\omega}}},:])\neq 0$. Therefore, $\Msub$ is full rank, or equivalently, $rank(\M[:,\widehat{0}]) = 6$. It follows from  (ii) and (iv) that $rank(\mrow{1},\mrow{6},\mrow{7},\mrow{3},\mrow{5})= 4$.\\
On the other hand, (ii) and (iv) imply that $$rank(\Msub(\V{\omega}+\V{\pi}_2))=rank(\mrow{0},\mrow{4},\mrow{6},\mrow{1},\mrow{3},\mrow{5},\mrow{7})= 5$$ and likewise $$rank(\Msub(\V{\omega}+\V{\pi}_4))=rank(\mrow{0},\mrow{2},\mrow{6},\mrow{1},\mrow{3},\mrow{5},\mrow{7})= 5.$$ Therefore, $\det(\Msub(\V{\omega} + \V{\pi}_2)) = \det(\Msub(\V{\omega} + \V{\pi}_4)) = 0$ and \eqref{eq: m0-cramer} implies $m_0(\V{\omega}+\V{\pi}_2) = m_0(\V{\omega}+\V{\pi}_4) = 0$.\\
If $\mrow{1}$ and $\mrow{7}$ are linearly independent and so are $\mrow{3}$ and $\mrow{5}$, then $$rank(\Msub(\V{\omega}+\V{\pi}_6))=rank(\mrow{2},\mrow{4},\mrow{1},\mrow{3},\mrow{5},\mrow{7}) = 6,$$ hence $m_0(\V{\omega}+\V{\pi}_6)\neq 0$. Therefore, $$[m_0(\V{\omega}),m_0(\V{\omega}+\V{\pi}_2),m_0(\V{\omega}+\V{\pi}_4),m_0(\V{\omega}+\V{\pi}_6)] = [*,0,0,*].$$ In addition, $d_{i,j} = 0,\, \forall(i,j)$ except $(0,6)$, so in \eqref{eq: singular-cond} $$\mathfrak{D}(\V{\omega}) = [d_{0,6}, 0, 0,0]^\top [0,0,0,1] + [0,0,0,d_{0,6}]^\top [-1,0,0,0].$$  By Proposition \ref{prop: feasibility}, \eqref{eq: identity-cond} cannot be satisfied, hence $rank(\mrow{1},\mrow{7})\leq 1$ or $rank(\mrow{3},\mrow{5})\leq 1$.\\
As $rank(\mrow{1},\mrow{6},\mrow{7},\mrow{3},\mrow{5}) = 4$, we must have $rank(\mrow{1},\mrow{7}) = 1$ and $rank(\mrow{3},\mrow{5}) = 2$ or $rank(\mrow{3},\mrow{5}) = 1$ and $rank(\mrow{1},\mrow{7}) = 2$.\qed

\begin{lemma}\label{lem: full-rank-m35}
Let $\widetilde{S_\rho} = S_\rho\cap\{\V{\omega}:\, rank(\,\mrow{3}(\V{\omega}),\,\mrow{5}(\V{\omega})\,)= 1\}$, if $\m{3}$ and $\m{4}$ concentrate in $T_3$ and $T_4$ respectively, then $|\widetilde{S_\rho}|=0$.
\end{lemma}
\noindent{\it Proof.}
Let $\widetilde{S_\rho} + \V{\pi}_3 = \{\V{\omega} + \V{\pi}_3,\, \V{\omega}\in \widetilde{S_\rho}\}$ and $\Omega'$ be the set symmetric to a set $\Omega\subset S_0$ with respect to the diagonal $\omega_1 = -\omega_2$. If $|\widetilde{S_\rho}|>0$, by the concentration of $\m{3}$ in $T_3$, $\forall\, \Omega\subset \widetilde{S_\rho} + \V{\pi}_3\subset T_3$ s.t. $|\Omega|> 0$, $\int_{\Omega}|\widetilde{m_3}| > \int_{\Omega'}|\widetilde{m_3}|$. Due to the symmetry between $|\widetilde{m_3}|$ and $|\widetilde{m_4}|$ defined in \eqref{eq: sym-m34}, $\int_{\Omega'}|\widetilde{m_3}| = \int_{\Omega}|\widetilde{m_4}|$. Therefore, $\int_{\Omega}|\widetilde{m_3}| >  \int_\Omega|\widetilde{m_4}|$ which implies that $|\m{3}| > |\m{4}|\; a.e.$ on $\widetilde{S_\rho}+\V{\pi}_3$ or equivalently $|\widetilde{m_3}(\V{\omega}+\V{\pi}_3)| > |\widetilde{m_4}(\V{\omega}+\V{\pi}_3)|\; a.e.$ on $\widetilde{S_\rho}$. Similarly, we have $|\widetilde{m_4}(\V{\omega}+\V{\pi}_5)| > |\widetilde{m_3}(\V{\omega}+\V{\pi}_5)|\; a.e.$ on $\widetilde{S_\rho}$ following the same analysis on $\widetilde{S_\rho}+\V{\pi}_5\subset T_4$.
On the other hand, $rank(\,\mrow{3}(\V{\omega}),\,\mrow{5}(\V{\omega})\,)= 1$ on $\widetilde{S_\rho}$, hence $\widetilde{m_3}(\V{\omega}+\V{\pi}_3)\widetilde{m_4}(\V{\omega}+\V{\pi}_5) = \widetilde{m_3}(\V{\omega}+\V{\pi}_5)\widetilde{m_4}(\V{\omega}+\V{\pi}_3)$, which contradicts the previous two inequalities.\qed

\begin{lemma}\label{lem: concentrate}
If $\m{1} $ {\rm (}respectively, $\m{6}${\rm)} concentrates in $T_1$ {\rm(}respectively, $T_6${\rm)}, then $|\m{6}| > |\m{1}|\,$ a.e. on $T_6\bigcap \text{supp}(\widetilde{m_6})$ {\rm (}respectively, $|\m{1}| > |\m{6}|$ 
a.e. on $T_1\bigcap\text{supp}(\widetilde{m_1})${\rm )}.
\end{lemma}
\noindent{\it Proof.}
Let $B_6=\{\V{\omega}: |\m{6}| \leq |\m{1}|\}\bigcap T_6\bigcap supp(\widetilde{m_1})$ and $B_1$ be the set symmetric to $B_6$ with respect to $\omega_1 = \omega_2$ and suppose $|B_6|>0$, then $\int_{B_6}|\m{6}|\leq \int_{B_6}|\m{1}|$. On the other hand, since $\m{1}$ concentrates in $T_1$, we know $\int_{B_1}|\m{1}| > \int_{B_6}|\m{1}|$. Moreover, due to the symmetry of $\m{1},\m{6}$ and $B_1,B_6$, $\int_{B_1}|\m{1}| = \int_{B_6}|\m{6}|$, hence $\int_{B_6}|\m{1}| \geq\int_{B_6}|\m{6}| = \int_{B_1}|\m{1}| $ which results in contradiction.\qed

\begin{proposition}\label{prop: zero-corner}
If  $\m{0},\,\m{1}$ and $\m{6}$ concentrate in $C_0,\,T_1$ and $T_6$ respectively, then $\m{6}=0$ a.e. on $S_\rho'+\V{\pi}_1$, where $S_\rho' =  S_\rho\bigcap\{\omega_1<\omega_2\}$.
%then $\m{1} = \m{6} = 0,\, a.e. $ on $ S_\rho + \V{\pi}_1$.
\end{proposition}
\noindent{\it Proof.}
By Lemma \ref{lem: concentrate}, the concentration of $\m{1}$ in $T_1$ implies that $|\widetilde{m_6}(\V{\omega}+\V{\pi}_1)| > |\widetilde{m_1}(\V{\omega}+\V{\pi}_1)|$ a.e. on $S_\rho'\cap \{\V{\omega},\,\widetilde{m_6}(\V{\omega}+\V{\pi}_1)\neq 0\}$. Similarly, the concentration of $\m{6}$ in $T_6$ implies that $|\widetilde{m_1}(\V{\omega}+\V{\pi}_7)| > |\widetilde{m_6}(\V{\omega}+\V{\pi}_7)|$ a.e. on $S_\rho'\cap \{\V{\omega},\,\widetilde{m_1}(\V{\omega}+\V{\pi}_7)\neq 0\}$. Therefore, $|\widetilde{m_1}(\V{\omega}+\V{\pi}_7)\widetilde{m_6}(\V{\omega} + \V{\pi}_1)| > |\widetilde{m_1}(\V{\omega}+\V{\pi}_1)\widetilde{m_6}(\V{\omega}+\V{\pi}_7)|$ a.e. on $S_\rho'\cap \{\V{\omega},\,\widetilde{m_6}(\V{\omega}+\V{\pi}_1)\neq 0\}\cap \{\V{\omega},\,\widetilde{m_1}(\V{\omega}+\V{\pi}_7)\neq 0\}$.

On the other hand, Lemma \ref{lem: rank1} implies that for a.e. $\V{\omega}\in S_\rho'$, $rank(\mrow{1}(\V{\omega}),\,\mrow{7}(\V{\omega})) = 1$, hence $\widetilde{m_1}(\V{\omega}+\V{\pi}_7)\widetilde{m_6}(\V{\omega} + \V{\pi}_1) = \widetilde{m_1}(\V{\omega}+\V{\pi}_1)\widetilde{m_6}(\V{\omega}+\V{\pi}_7)$. Together with the previous result, this forces $|S_\rho'\cap \{\V{\omega},\,\widetilde{m_6}(\V{\omega}+\V{\pi}_1)\neq 0\}\cap \{\V{\omega},\,\widetilde{m_1}(\V{\omega}+\V{\pi}_7)\neq 0\}| = 0$.

The concentration of $\m{0},\m{1}$ and $\m{6}$ in $C_0,\,T_1$ and $T_6$ implies that $\widetilde{m_1}(\V{\omega}+\V{\pi}_7)\neq 0$ on $S_\rho'$, since $\V{\omega}+\V{\pi}_7\not \in C_0\cup T_6,\; \forall \V{\omega}\in S_\rho'$ and neither $\widetilde{m_6}$ or $\widetilde{m_0}$ can dominate at $\V{\omega}+\V{\pi}_7$.
% if $\widetilde{m_1}(\V{\omega}+\V{\pi}_7) = 0$ then $\widetilde{m_6}(\V{\omega}+\V{\pi}_7) = 0$ on $S_\rho'$. It follows that $\widetilde{m_0}(\V{\omega}+\V{\pi}_7) > 0 = \widetilde{m_j}(\V{\omega}+\V{\pi}_7),\, \forall j = 1,\cdots,6$, i.e. $\V{\omega}+\V{\pi}_7 \in \Omega_0$, which will contradict the assumption that $\m{0}$ concentrates in $C_0$, or equivalently $\Omega_0\subset C_0$.  
Therefore, $S_\rho'\cap\{\V{\omega},\widetilde{m_1}(\V{\omega}+\V{\pi}_7)\neq 0\} = S_\rho'$ which implies $|S_\rho'\cap\{\V{\omega},\,\widetilde{m_6}(\V{\omega}+\V{\pi}_1)\neq 0\}|=0$, i.e. $\m{6}=0$ a.e. on $S_\rho'+\V{\pi}_1$.
\qed

%$\exists\,\alpha_{\V{\omega}}\in\mathbb{C}, s.t.\,\mrow{1}(\V{\omega}) = \alpha_{\V{\omega}}\,\mrow{7}(\V{\omega}),\,$ $\forall\, \V{\omega}\in S_\rho',$ i.e. $\widetilde{m_1}(\V{\omega} + \V{\pi}_1) = \alpha_{\V{\omega}}\cdot\widetilde{m_1}(\V{\omega} + \V{\pi}_7)$ and $\widetilde{m_6}(\V{\omega} + \V{\pi}_1) = \alpha_{\V{\omega}}\cdot\widetilde{m_6}(\V{\omega} + \V{\pi}_7)$. On the other hand, Lemma \ref{lem: concentrate} implies that $|\widetilde{m_1}(\V{\omega} + \V{\pi}_7)| \geq |\widetilde{m_6}(\V{\omega} + \V{\pi}_7)|$, hence $|\widetilde{m_1}(\V{\omega} + \V{\pi}_1)| \geq |\widetilde{m_6}(\V{\omega} + \V{\pi}_1)|$. Let $\Omega_6'\doteq (S_\rho+\V{\pi}_1)\bigcap T_6$, then $\int_{\Omega_6'}|\m{1}| \geq\int_{\Omega_6'}|\m{6}|$, which will contradict Lemma \ref{lem: concentrate} unless $|\Omega_6'\bigcap\text{supp}(\widetilde{m_6})| = 0$, or equivalently $\alpha_{\V{\omega}}=0$ and so $\m{6} = \m{1} = 0,\,a.e.$ on $\Omega_6'$. By symmetry, $\m{6}=\m{1} = 0,\,a.e. $ on $(S_\rho+\V{\pi}_1)\setminus \Omega_6'$ as well.\qed

\subsection{Design of input $\m{j}$}\label{app: input design}
{\bf Proof of Lemma \ref{lem: phase-ineq}}:\\[.2em]
{\bf Lemma \ref{lem: phase-ineq}.} 
{\it If $\exists\,\V{\omega}\in D_1:=\{\omega_1=\omega_2,\,\omega_1\in(-\frac{\pi}{2},0)\},\,s.t. \,|m_0(\V{\omega})|\neq 0,$ then $(\V{\eta}_1-\V{\eta}_6)^\top (\V{\pi}_6-\V{\pi}_7)\neq 0(\text{mod}\,2\pi)$. 
}\\ [1em]
\noindent {\it Proof.}
As $\m{1}$ and $\m{6}$ concentrate in $T_1$ and $T_6$ respectively, $\widetilde{m_1}(\V{\omega} + \V{\pi}_i) = 0$ and  $\widetilde{m_6}(\V{\omega} + \V{\pi}_i) = 0$, $i = 1,\cdots, 5$. Due to symmetry, $|\widetilde{m_1}(\V{\omega})| = |\widetilde{m_6}(\V{\omega})|$ on $\{\omega_1=\omega_2\}$. Let $A = |\widetilde{m_1}(\V{\omega}+\V{\pi}_7)| = |\widetilde{m_6}(\V{\omega}+\V{\pi}_7)|$ and $B=|\widetilde{m_1}(\V{\omega}+\V{\pi}_6)| = |\widetilde{m_6}(\V{\omega}+\V{\pi}_6)|$, then the first and the last columns of $\Msub$ are
  \begin{align*}
  \Msub[:,1] = 
 \begin{bmatrix}
 0\\
 \vdots\\
 0\\
 Ae^{i\V{\eta}_1^\top(\V{\omega}+\V{\pi}_6)}\\
 Be^{i\V{\eta}_1^\top(\V{\omega}+\V{\pi}_7)}
 \end{bmatrix}
 \quad\text{and}\quad
  \Msub[:,6] = 
 \begin{bmatrix}
 0\\
 \vdots\\
 0\\
 Ae^{i\V{\eta}_6^\top(\V{\omega}+\V{\pi}_6)}\\
 Be^{i\V{\eta}_6^\top(\V{\omega}+\V{\pi}_7)}
 \end{bmatrix} .
\end{align*}   
By \eqref{eq: m0-cramer}, if $m_0(\V{\omega})>0, \,\V{\omega}\in D_1$ then $\Msub(\V{\omega})$ is full rank, hence its columns are linearly independent.
In particular, $\Msub[:,1]$ and $\Msub[:,6]$ are linearly independent, which implies that $e^{i(\V{\eta}_1^\top\V{\pi}_6 + \V{\eta}_6^\top\V{\pi}_7)}\neq e^{i(\V{\eta}_6^\top\V{\pi}_6 + \V{\eta}_1^\top\V{\pi}_7)}$%$e^{i(\V{\eta}_1-\V{\eta}_6)^\top(\V{\omega}+\V{\pi}_6)}\neq e^{i(\V{\eta}_1-\V{\eta}_6)^\top(\V{\omega}+\V{\pi}_7)}$ 
 or equivalently $(\V{\eta}_1-\V{\eta}_6)^\top(\V{\pi}_6-\V{\pi}_7)\neq 0(\text{mod}2\pi)$. \qed\\

\noindent{\bf Proof of Proposition \ref{prop: origin-det}}\\[.2em]
\noindent{\bf Proposition \ref{prop: origin-det}.} 
{\it If $\widetilde{m_0}(\V{0})\neq 0,$ then $\V{\pi}_1^\top(\V{\eta}_1-\V{\eta}_6)\neq \pi(\text{mod}\,2\pi)$ or $\V{\pi}_3^\top(\V{\eta}_3-\V{\eta}_4)\neq \pi(\text{mod}\,2\pi)$. }\\[1em]
\noindent{\it Proof.}
%$\Msub(\V{0})$ takes the following form
%$$\begin{bmatrix}
%* & 0 & 0 & 0 & 0 & *\\
%0 & * & 0 & 0 & 0 & 0\\
%0 & 0 & * & * & 0 & 0\\
%0 & 0 & 0 & 0 & * & 0\\
%0 & 0 & * & * & 0 & 0\\
%* & 0 & * & * & 0 & *\\
%* & 0 & 0 & 0 & 0 & *
%\end{bmatrix}$$
%The second and the fifth columns of $\Msub$ have single non-zero entry, $\widetilde{m_2}(\V{\pi}_2)$ and $\widetilde{m_5}(\V{\pi}_4)$ respectively, and are orthogonal to all the rest columns, hence the full-rank constraint of $\Msub$ is reduced to the full-rank constraint on its sub-matrix (with permutation of rows and columns)
 Since $\widetilde{m_0}(\V{0})\neq 0$, as shown in Lemma \ref{lem: rank1}, at $\V{\omega} = \V{0}$ $rank(\mrow{1},\mrow{6},\mrow{7},\mrow{3},\mrow{5})= 4$ . This is equivalent to the matrix $\V{A}$ defined in \eqref{eq: matrix-B} to be full rank.
\begin{align}\label{eq: matrix-B}
%\mbox{\V{A}\strut}=
\V{A} = 
\begin{bmatrix}
& & & \\[-1em]
\widetilde{m_1}(\V{\pi}_6) & \widetilde{m_6}(\V{\pi}_6) & \widetilde{m_3}(\V{\pi}_6) & \widetilde{m_4}(\V{\pi}_6) \\
\widetilde{m_1}(\V{\pi}_1) & \widetilde{m_6}(\V{\pi}_1) & 0 & 0\\
\widetilde{m_1}(\V{\pi}_7) & \widetilde{m_6}(\V{\pi}_7) & 0 & 0\\
0 & 0 & \widetilde{m_3}(\V{\pi}_3) & \widetilde{m_4}(\V{\pi}_3)\\
0 & 0 & \widetilde{m_3}(\V{\pi}_5) & \widetilde{m_4}(\V{\pi}_5)\\
\end{bmatrix}
\end{align}
Let $|\widetilde{m_1}(\V{\pi}_1)| = a, \, |\widetilde{m_1}(\V{\pi}_6)|=b$. Due to the symmetry of $\m{j}$,
$|\widetilde{m_1}(\V{\pi}_1)| = |\widetilde{m_1}(\V{\pi}_7)| = |\widetilde{m_6}(\V{\pi}_1)| = |\widetilde{m_6}(\V{\pi}_7)| = |\widetilde{m_3}(\V{\pi}_3)| = |\widetilde{m_3}(\V{\pi}_5)| = |\widetilde{m_4}(\V{\pi}_3)| = |\widetilde{m_4}(\V{\pi}_5)|$ and $|\widetilde{m_1}(\V{\pi}_6)|=| \widetilde{m_6}(\V{\pi}_6)|= | \widetilde{m_3}(\V{\pi}_6)|=| \widetilde{m_4}(\V{\pi}_6)|$. Rewrite $\V{A}$ as follows,
$$\V{A}=
\begin{bmatrix}
b e^{-i\V{\pi}_6^\top\V{\eta}_1} & b e^{-i\V{\pi}_6^\top\V{\eta}_6} & b e^{-i\V{\pi}_6^\top\V{\eta}_3} & b e^{-i\V{\pi}_6^\top\V{\eta}_4}\\
a e^{-i\V{\pi}_1^\top\V{\eta}_1} & a e^{-i\V{\pi}_1^\top\V{\eta}_6} & 0						& 0 \\
a e^{i\V{\pi}_1^\top\V{\eta}_1} & a e^{i\V{\pi}_1^\top\V{\eta}_6} & 0						& 0 \\
0 					& 0 					& a e^{-i\V{\pi}_3^\top\V{\eta}_3} & a e^{-i\V{\pi}_3^\top\V{\eta}_4}\\
0 					& 0 					& a e^{i\V{\pi}_3^\top\V{\eta}_3} & a e^{i\V{\pi}_3^\top\V{\eta}_4}\\
\end{bmatrix}
$$
The product of singular values of $\V{A}$ is 
\begin{align}\label{eq: detB}
\sqrt{\text{det}(\V{A}^* \V{A})} = 4a^3\sqrt{a^2 K_1^2K_2^2 + b^2(Q_1K_2^2 + Q_2K_1^2)},
\end{align}
where $ Q_1 = 1 - \cos(\V{\pi}_6^\top(\V{\eta}_1-\V{\eta}_6))\cos(\V{\pi}_1^\top(\V{\eta}_1-\V{\eta}_6)), Q_2 = 1 - \cos(\V{\pi}_6^\top(\V{\eta}_3-\V{\eta}_4))\cos(\V{\pi}_3^\top(\V{\eta}_3-\V{\eta}_4)), K_1 = \sin(\V{\pi}_1^\top(\V{\eta}_1-\V{\eta}_6)), K_2 = \sin(\V{\pi}_3^\top(\V{\eta}_3-\V{\eta}_4)).$ If $\V{\pi}_1^\top(\V{\eta}_1-\V{\eta}_6) = \V{\pi}_3^\top(\V{\eta}_3-\V{\eta}_4) = \pi (mod\, 2\pi)$, then $K_1 = K_2 = 0$ and $\V{A}$ becomes singular.\qed
\subsection{Solving \eqref{eq: LS-new} and \eqref{eq: identity-cond} for $m_0,\widetilde{m_0}$ and $m_j$}\label{app: solving}
\begin{lemma}\label{lem: null-space}
Let $\V{P}\in\mathbb{C}^{n\times n}$ be a projection matrix of rank $2$ and $\V{a},\V{b},\V{a}',\V{b}'\in\mathbb{C}^n,\, s.t.\, \V{a}^*\V{b} = (\V{a}')^*\V{b}'=1,\, \V{a}'^*\V{b} = \V{a}^*\V{b}' = \V{b}^*\V{b}' = 0.$ If $\V{P}(\V{I}_n - \V{a}\otimes\V{b} - \V{a}'\otimes\V{b}' ) = \V{0}$, then $\V{P}$ is the projection of $span\{\V{b},\V{b}'\}$.
\end{lemma}
\noindent{\it Proof.}
Since $$rank(\V{I}_n) \leq rank(\V{I}_n-\V{a}\otimes\V{b}-\V{a}'\otimes\V{b}') + rank(\V{a}\otimes\V{b}) + rank(\V{a}'\otimes\V{b}'),$$
it follows that $rank(\V{I}_n-\V{a}\otimes\V{b}-\V{a}'\otimes\V{b}')\geq n - 2$. On the other hand, because $rank(\V{P}) = 2$, $\V{P}(\V{I}_n - \V{a}\otimes\V{b} - \V{a}'\otimes\V{b}' ) = \V{0}$ implies that $rank(\V{I}_n-\V{a}\otimes\V{b}-\V{a}'\otimes\V{b}')\leq n - 2$. Hence $rank(\V{I}_n-\V{a}\otimes\V{b}-\V{a}'\otimes\V{b}') = n - 2$ and $\V{P}$ is the projection of $col(\V{I}_n-\V{a}\otimes\V{b}-\V{a}'\otimes\V{b}')^\bot$. On the other hand,
\begin{align*}
\V{b}^*(\V{I}_n-\V{a}\otimes\V{b}-\V{a}'\otimes\V{b}')
&= \V{b}^* - (\V{b}^*\V{a})\V{b}^* - (\V{b}^*\V{a}')(\V{b}')^*\\
&= \V{b}^* - \V{b}^* - 0\cdot (\V{b}')^* = \V{0}^*.
\end{align*}
Therefore, $\V{P}\V{b} = \V{b}$.
Similarly, $(\V{b}')^* (\V{I}_n-\V{a}\otimes\V{b}-\V{a}'\otimes\V{b}')=\V{0}^*$ and $\V{P}\V{b}' = \V{b}'$. Moreover, as $\V{b}^*\V{b}' = 0$ and $rank(\V{P}) = 2$, $\V{P} = \Vert\V{b}\Vert^{-2}\cdot\V{b}\otimes\V{b} + \Vert\V{b}'\Vert^{-2}\cdot\V{b}'\otimes\V{b}'.$\qed

\begin{lemma}
Given $\M[:,\widehat{0}](\V{\omega})$ is full rank $\forall \V{\omega}$, $\M[\widehat{0},:](\V{\omega})$ is singular if \eqref{eq: identity-cond} holds.
\end{lemma}
\noindent{\it Proof. }
If \eqref{eq: identity-cond} holds, then by Lemma \ref{lem: null-space}, $\meven,\,\modd$ are orthogonal to \\$col(\M[:,\widehat{0}])$, therefore $\big[\,\modd,\meven,\M[:,\widehat{0}]\,\big]\in\mathbb{C}^{8\times 8}$ is full rank. Due to \eqref{eq: identity-cond}, $\meven$ and $\mteven$ are not orthogonal to each other, hence $\big[\,\modd,\mteven,\M[:,\widehat{0}]\,\big] = \big[\,\modd, \M\,\big]$ is full rank as well. Because $(\modd)^*\M[:, i] = 0,\,i= 0,\cdots,7$ and $\modd[\widehat{0}]^*\M[\widehat{0}, i] = (\modd)^*\M[:,i]$, $\modd[\widehat{0}]$ is orthogonal to $col(\M[\widehat{0},:])$. Since $\big[ \modd[\widehat{0}], \M[\widehat{0},:]\,\big]\in\mathbb{C}^{7\times 8}$ is full rank, $\M[\widehat{0},:]$ must be singular.\qed\\[1em]

\noindent{\bf Proof of Proposition \ref{prop: m0_formula}}:\\[.2em]
\noindent{\bf Proposition \ref{prop: m0_formula}.} 
{\it Let $\M[odd,\widehat{0}](\V{\omega}),\M[even,\widehat{0}](\V{\omega})\in\mathbb{C}^{4\times6}$ be the submatrices of $\M[:,\widehat{0}](\V{\omega})$ consisting of odd and even indexed rows respectively. $\forall\V{\omega}\in S_0$, suppose {\rm\ref{cond: i}} and \eqref{eq: identity-cond} are satisfied, then {\rm\ref{cond: ii}} holds if and only if $rank(\,\M[odd,\widehat{0}](\V{\omega})\,) = rank(\,\M[even,\widehat{0}](\V{\omega})\,) = 3$ and 
\begin{align}%\label{eq: m0-even-null}
[m_0(\V{\omega}),m_0(\V{\omega}+\V{\pi}_2), m_0(\V{\omega} +\V{\pi}_4), m_0(\V{\omega}+\V{\pi}_6)]\, \M[even,\widehat{0}](\V{\omega}) = \V{0}, \tag{\ref{eq: m0-even-null}}
\end{align}
\begin{align}%\label{eq: m0-odd-null}
[m_0(\V{\omega}+\V{\pi}_1),m_0(\V{\omega}+\V{\pi}_3), m_0(\V{\omega} +\V{\pi}_5), m_0(\V{\omega}+\V{\pi}_7)] \,\M[odd,\widehat{0}](\V{\omega}) = \V{0}. \tag{\ref{eq: m0-odd-null}}
\end{align}
}
\\[.5em]
\noindent{\it Proof.}
Note that $\M[:,\widehat{0}]$ have the same rows at $\V{\omega} + \V{\pi}_i,\,i = 0,\cdots,7$, we define row permutation matrix $\V{P}_i,\; s.t.\,$ $\V{P}_i\big(\M[:,\widehat{0}](\V{\omega} + \V{\pi}_i)\big) = \M[:,\widehat{0}](\V{\omega}). $ Let $\V{P}_{\M}(\V{\omega})$ be the projection matrix of the $col\big(\M[:,\widehat{0}](\V{\omega})\big)^\bot = null(\M[:,\widehat{0}]^*)$, then \ref{cond: ii} is equivalent to $\V{P}_{\M}\V{b}_0'(\V{\omega}) = \V{0}.$ Group this equality at $\V{\omega}+\V{\pi}_i$, we have 
\begin{align}\label{eq: group-proj}
\V{0} & = [\V{P}_i\V{P}_{\M}\V{b}_0'(\V{\omega} + \V{\pi}_i) ]_{i=0,\cdots,7}\notag\\
&= [\V{P}_i\V{P}_{\M}(\V{\omega}+\V{\pi}_i)\V{P}_i^2\V{b}_0'(\V{\omega}+\V{\pi}_i) ]_{i = 0,\cdots,7}\notag\\
&= [\V{P}_{\M}(\V{\omega})\V{P}_i\V{b}_0'(\V{\omega}+\V{\pi}_i)]_{i = 0,\cdots,7}\notag\\
&= \V{P}_{\M}(\V{\omega})[\V{P}_i\V{b}_0'(\V{\omega}+\V{\pi}_i)]_{i = 0,\cdots,7}
\end{align}
Let 
\begin{align*}
\mteven&= [ (1 + i \bmod 2)\cdot\,\sbarmp{0}{i}]_{i=0,\cdots,7}^\top = \M[:,0](\V{\omega}),\\
\mtodd&= [ (i \bmod 2)\cdot\,\sbarmp{0}{i}]_{i=0,\cdots,7}^\top,\\
\meven&= [ (1 + i \bmod 2)\cdot\,m_0(\V{\omega} + \V{\pi}_i)]_{i=0,\cdots,7}^\top,\\
\modd&= [ ( i \bmod 2)\cdot\,m_0(\V{\omega} + \V{\pi}_i)]_{i=0,\cdots,7}^\top.
\end{align*}
The identity constraint \eqref{eq: identity-cond} thus can be written as $(\overlinespace{\meven})^*\,\mteven = 1$ and $(\overlinespace{\modd})^*\,\mtodd = 1$. By definition, 
$$ \V{P}_i\V{b}_0'(\V{\omega}+\V{\pi}_i) = \V{P}_i\big(\V{b}_0 - m_0\M[:,0](\V{\omega}+\V{\pi}_i)\big)  = \V{b}_i - m_0(\V{\omega}+\V{\pi}_i)\V{P}_i\big(\M[:,0](\V{\omega}+\V{\pi}_i)\big)$$
and 
$$\V{P}_i\big(\M[:,0](\V{\omega}+\V{\pi}_i)\big) = 
\begin{cases}
   \M[:,0] = \mteven, &  i \text{ is even}\\[.2em]
    \mtodd,              & i \text{ is odd}
\end{cases}
$$
Substitute the above expression of $ \V{P}_i\V{b}_0'(\V{\omega}+\V{\pi}_i)$ in \eqref{eq: group-proj} and we have
\begin{align}
\V{0}  
%= \V{P}_{\M}(\omega)([\V{b}_i]_{i=0,\cdots,7} - \mteven\otimes\meven - \mtodd\otimes\modd) 
= \V{P}_{\M}(\V{I}_8 -  \mteven\otimes\overline{\meven} - \mtodd\otimes\overline{\modd})
\end{align}
Therefore, by Lemma \ref{lem: null-space}, $\V{P}_{\M}$ is the projection of $span\{\overlinespace{\modd},\overlinespace{\meven}\}$. This is equivalent to \eqref{eq: m0-even-null} and \eqref{eq: m0-odd-null}.
Finally, since $$6 = rank(\M[:,\widehat{0}]) \leq rank(\M[odd,\widehat{0}]) + rank(\M[even,\widehat{0}]) \leq (4-1) + (4-1),$$ $rank(\M[odd,\widehat{0}]) = rank(\M[even,\widehat{0}]) = 3$.\qed

%% file: QCQP.tex
\section{Joint optimization of $c(\V{\omega})$ and $\m{0}$}\label{app: QCQP}

In \ref{alg}, $c(\V{\omega})$ is chosen in step 3. to construct $m_0'(\V{\omega})$, which replaces $m_0(\V{\omega})$ and is used to create the linear constraint in \eqref{eq: opt} in step 4. Since different $c(\V{\omega})$ correspond to different $m_0'(\V{\omega})$, hence different linear constraints \eqref{eq: m0-A} on $\m{0}$; $\m{0}$ obtained in step 4. is optimal with respect to the pre-fixed $c(\V{\omega})$ from step 3., but not necessarily global optimal considering all possible choices of $c(\V{\omega})$. Therefore, we propose an alternative approach that combines step 3. and step 4. in \ref{alg}, where $c(\V{\omega})$ and $\m{0}$ are jointly optimized to obtain $\m{0}$ with the best possible regularity given unregularized $m_0(\V{\omega})$ from step 2.

By the definition in Proposition \ref{prop: mc}, $m_0'(\V{\omega}) = m_0(\V{\omega})c(\V{\omega})$. Furthermore, since $c(\V{\omega})$ is $\pi$-periodic in both $\omega_1,\,\omega_2$, we have $m_0'(\V{\omega}+\V{\pi}_i) = m_0(\V{\omega}+\V{\pi}_i)c(\V{\omega}),\, i = 2,4,6$. Hence the constraint \eqref{eq: identity-cond} on $\m{0}$ with $m_0(\V{\omega})$ replaced by $m_0'(\V{\omega})$ can be reformulated as follows,
\begin{align}
1 & = m_0'\m{0} + m_0'\mp{0}{2} + m_0'\mp{0}{4} + m_0'\mp{0}{6}\notag\\
& = c(\V{\omega})\big(\,m_0\m{0} + m_0\mp{0}{2} + m_0\mp{0}{4} + m_0\mp{0}{6}\,\big).\label{eq: linear-cond_c}
\end{align}
Using the same setup of the optimization \eqref{eq: opt}, we convert \eqref{eq: linear-cond_c} to a constraint on a $2N\times 2N$ grid $\mathcal{G}=\{\V{\omega}_i\}_{i=1}^{4N^2}$ of $[-\pi,\pi)\times[-\pi,\pi)$.
Let $\widetilde{\mathbf{m}_0}\in\mathbb{C}^{4N^2}$ and $\V{A}_0\in\mathbb{C}^{N^2\times 4N^2}$ be the same as in \eqref{eq: m0-A} except that $\V{A}_0$ is constructed by unregularized $m_0$ instead of $m_0'$ for $\V{A}$.
Let $\V{C}\in\mathbb{C}^{N^2\times N^2}$ be a diagonal matrix whose $j$-th diagonal entry is $c(\V{\omega}_j)$, where $\V{\omega}_j\in \mathcal{G}\cap [-\pi,0)\times[-\pi,0)$ in the same order as the rows of $\V{A}_0$. Then \eqref{eq: linear-cond_c} is equivalent to the following constraint on the grid $\mathcal{G}$,
\begin{align}
\V{C}\V{A}_0\,\overlinespace{\widetilde{\mathbf{m}_0}} = \V{1}_{N^2}.
\end{align}
We formulate the joint optimization on $\V{C}$ and $\widetilde{\mathbf{m}_0}$ analogous to \eqref{eq: opt} as follows,
\begin{align}\label{eq: opt-C}
\min_{\xvec\in\mathbb{C}^{4N^2},\;\mathbf{c}\in\mathbb{C}^{N^2}}\; \Vert \V{D}\xvec\Vert^2,\quad 
s.t. \; \V{C}\V{A}_0\,\xvec = \mathbf{1}, \; \V{C} = diag(\mathbf{c}).
\end{align}
Since the objective function does not involve $\mathbf{c}$, $\mathbf{c}$ can be expressed in terms of $\xvec$ as long as $\V{A}_0\,\xvec$ has no zero entry. Therefore, solving \eqref{eq: linear-cond_c} is equivalent to solving the following optimization for $\widetilde{\mathbf{m}_0}$.
\begin{align}\label{eq: opt-ineq}
\min_{\xvec\in\mathbb{C}^{4N^2}}\; \Vert \V{D}\xvec\Vert^2,\quad 
s.t. \; |\V{A}_0\,\xvec| > 0,
\end{align}
where $|\cdot|$ in the constraint is a pointwise operator that computes the absolute value. The constraint $|\V{A}_0\,\xvec| > 0$ can be rewritten as a set of quadratic constraints $\xvec^*\V{Q}_i\xvec > 0,\, i =0,\cdots, N^2-1 $ where $\V{Q}_i = \V{A}_0[i,:]^*\V{A}_0[i,:]$. Therefore, \eqref{eq: opt-ineq} is a quadratically constrained quadratic program. Furthermore, since $\V{Q}_i$ is positive semi-definite, \eqref{eq: opt-ineq} is not convex and is NP-hard in general. One may solve the convex relaxation of \eqref{eq: opt-ineq} using semidefinite programming (SDP). Instead of solving $\xvec$, we solve $\V{X}\doteq \xvec\xvec^*$ and convert \eqref{eq: opt-ineq} into 
\begin{align}\label{eq: opt-trace}
\min_{\V{X}\in\mathbb{C}^{4N^2\times 4N^2}}\;tr(\, \V{D}^*\V{D}\,\V{X}\,),\quad
s.t. \; tr(\,\V{Q}_i\,\V{X}\,) > 0,\, \V{X}\succeq 0,\, rank(\V{X}) = 1,
\end{align}
where $\V{X}\succeq 0$ is the positive semidefinite constraint on $\V{X}$.
By removing the non-convex rank constraint $rank(\V{X}) = 1$, \eqref{eq: opt-trace} becomes a SDP and can be efficiently solved. Yet the solution $\V{X}$ may not be rank 1 and require post processing (e.g. singular value decomposition) to obtain an approximate solution of \eqref{eq: opt-ineq}.

%% file: m0-opt.tex
\pdfoutput=1
\section{Supplementary Numerical Results}\label{app: supp-numerical}
\subsection{Numerical optimization of $\m{0}$ in 1D}\label{subsec: 1D-opt}
To test whether numerical optimization is a practical way to solve \eqref{eq: identity-cond}, we experiment on $m_0$ and $\widetilde{m_0}$ of existing real biorthogonal wavelets. We consider a pair of low frequency filters corresponding to biorthogonal scaling functions $\phi,\, \tilde{\phi}$ with vanishing moments 3 and 5 respectively. 

\begin{figure}
\centering
%\begin{wrapfigure}{r}{.4\textwidth}
\includegraphics[width = .4\textwidth]{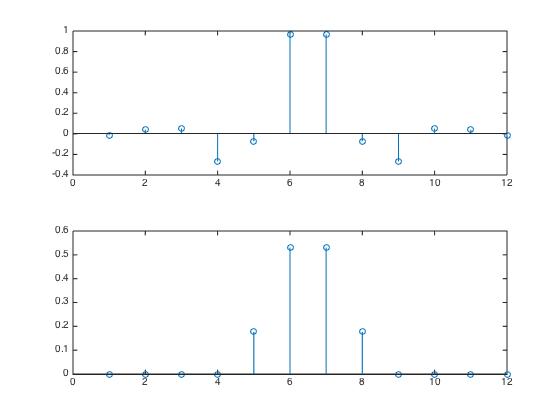}
\caption{1D filters, up: LoD, down: LoR}
\label{fig: filters}
%\end{wrapfigure}
\end{figure}
\begin{figure}%{l}{.4\textwidth}
\begin{minipage}[t]{.45\textwidth}
\includegraphics[width = \linewidth]{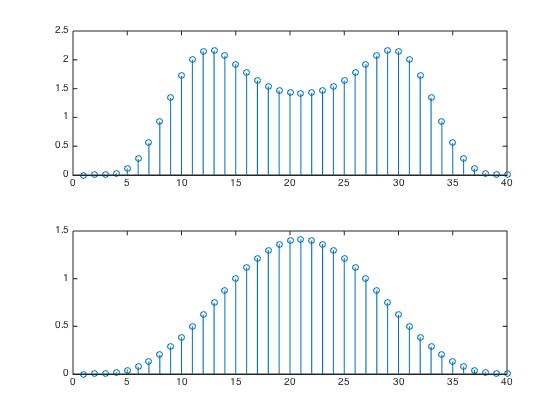}
\caption{$m_0(\omega)$ and $\widetilde{m_0}(\omega)$}
\label{fig: m-funcs}
\end{minipage}
\hfill
\begin{minipage}[t]{.45\textwidth}
\vbox{
\includegraphics[width = \textwidth]{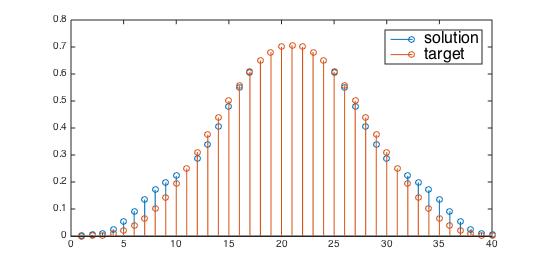}\\
\includegraphics[width = \textwidth]{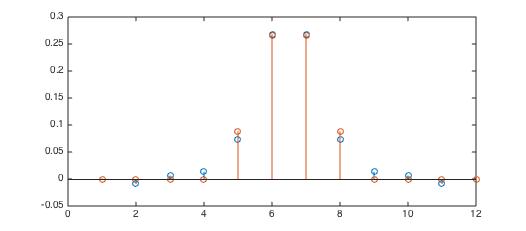}
}
\caption{$\widehat{\widetilde{m_0}}$ vs. $\widetilde{m_0}$, top: frequency domain, bottom: time domain}
\label{fig: 1d-compare}
\end{minipage}
\end{figure}

The 1D filters are shown in Figure \ref{fig: filters}. Suppose we know the decomposition filter, and we want to find the real reconstruction filter, such that it has support as concentrated as possible. 
%The corresponding $m_0$ and $\widetilde{m_0}$ are complex, yet we can shift the phase of $m_0$ such that $m_0$ is real and apply the same phase shift to $\m{0}$. 
%Without loss of generality, \eqref{eq: identity-cond} can be solved assuming that $m_0$ is real.
%It is not necessary that the corredponding $\widetilde{m_0}$ is also real, but in this testing case, $m_0$ and $\widetilde{m_0}$ are both real.
%have the same phase, hence the phase-shifted $\m{0}$ is real as well. 
Figure \ref{fig: m-funcs} shows the ground truth $m_0$ and $\widetilde{m_0}$ considered in this simulation. %and in particular, $|m_0|$ is used as the known coefficients in \eqref{eq: bi-orth-eq}. Hereafter, we use $m_0(\omega)$ and $\m{0}$ to denote the real-valued functions.

Let $\mhat{0}$ be the approximation of $\m{0}$, which is solution of the following optimization problem
\begin{align}
\min_{\xvec}\; \Vert \V{D}\xvec\Vert^2 + \Vert \xvec\Vert^2,\quad s.t. \; \V{A}\xvec = \mathbf{1} \label{eq: opt-1d}
\end{align}
where $\V{A}$ in the constraint is the matrix generated from $m_0\overlinespace{\widetilde{m_0}}(\omega) + m_0\overlinespace{\widetilde{m_0}}(\omega+\pi) = 1$, the 1D version of \eqref{eq: identity-cond}. Since only a single shift of $\pi$ appears in the condition, each row of $\V{A}$ has two non-zero entries. 
%Notice that no symmetry constraint is imposed here, nevertheless, 
Figure \ref{fig: 1d-compare} compares the solution of \eqref{eq: opt-1d} and the ground truth. The support of the solution is slightly more spread out than the ground truth.
%The support of the solution shown in Fig.\ref{fig: 1d-compare} is almost symmetric. On the other hand, its support in the time domain is not as compact as that of $\m{0}$, see the bottom of Fig.\ref{fig: 1d-compare}.

\subsection{Numerical optimization of $\m{0}$ in 2D}
In the 2D case, we use the pair of biorthogonal low-pass filters that are the tensor products of the 1D filters in Section \ref{subsec: 1D-opt} as ground truth. We solve the 2D version of the optimization problem \eqref{eq: opt-1d}. Figure \ref{fig: 2d-compare-1} shows the solution and compares it with the ground truth. 

%{\it 2D version of \eqref{eq: opt-1d}}\\

%The 1D formulation can be easily extended to 2D, where $\V{D} = [\V{D}_x,\V{D}_y]$ consider 1st order derivative in both $x$ and $y$ directions, and $\V{A}$ is generated from \eqref{eq: identity-cond}, each row has four non-zero entries. 
%It is obvious that the solution is not $90^\circ$-rotation invariant. Even worse is the fact that there is much energy in the vertical high-frequency domain.

%{\it weighted L2 norm (Modulation Space$^{[\ref{app: modulation}]}$)}\\
To make the support of $\mhat{0}$ better concentrate within the low frequency domain, we change the squared $\ell_2$-norm penalty in \eqref{eq: opt-1d} to a weighted version (corresponding to Modulation space) as follows,
\begin{align}
\min_{\xvec}\; \Vert\V{ D}\xvec\Vert^2 + \lambda\Vert \wvec\circ\xvec\Vert^2,\quad s.t. \; \V{A}\xvec = \mathbf{1} \label{eq: opt-2d-weight}
\end{align} 
where $\circ$ is Hadamard product and $\wvec$ is a weight vector. In particular, we choose $\forall \V{\omega}, \; \wvec(\V{\omega}) = |\V{\omega}|$. Figure \ref{fig: 2d-compare-2.2} shows the solution of \eqref{eq: opt-2d-weight} with $\lambda=600$. % and $600$ respectively. As $\lambda$ increases, the support of the minimizer concentrates more within the low frequency region. As shown in Fig.\ref{fig: 2d-compare-2}, when $\lambda$ is not huge, the minimizer achieves a certain level of but not full symmetry, whereas Fig.\ref{fig: 2d-compare-2.2} shows that huge $\lambda$ imposes full symmetry.

Compared to \eqref{eq: opt} proposed to solve $\m{0}$, both optimization problems \eqref{eq: opt-1d} and \eqref{eq: opt-2d-weight} in this simulation minimize the squared $\ell_2$-norm of the gradient of $\widetilde{m_0}$ but have an extra (weighted) $\ell_2$ regularization term. Although \eqref{eq: opt-1d} and \eqref{eq: opt-2d-weight} work better than \eqref{eq: opt} for 1D and 2D tensor wavelet construction here, they do not provide solutions with better regularity in the construction of biorthogonal directional wavelets while increasing the computation cost.

\begin{comment}
\begin{minipage}{.9\textwidth}
\centering
\includegraphics[width = .9\textwidth]{2d-m-compare-2-1-eps-converted-to.pdf}\\
\includegraphics[width = .9\textwidth]{2d-filter-compare-2-1-eps-converted-to.pdf}
\captionof{figure}{result of \eqref{eq: opt-2d-weight} $\mhat{0}$ ($\lambda = 60$), target $\m{0}$ and their difference, Top: frequency domain, Bottom: time domain}
\label{fig: 2d-compare-2}
\end{minipage}
\end{comment}

\begin{center}
\begin{minipage}{.8\textwidth}
\centering
\includegraphics[width = \textwidth]{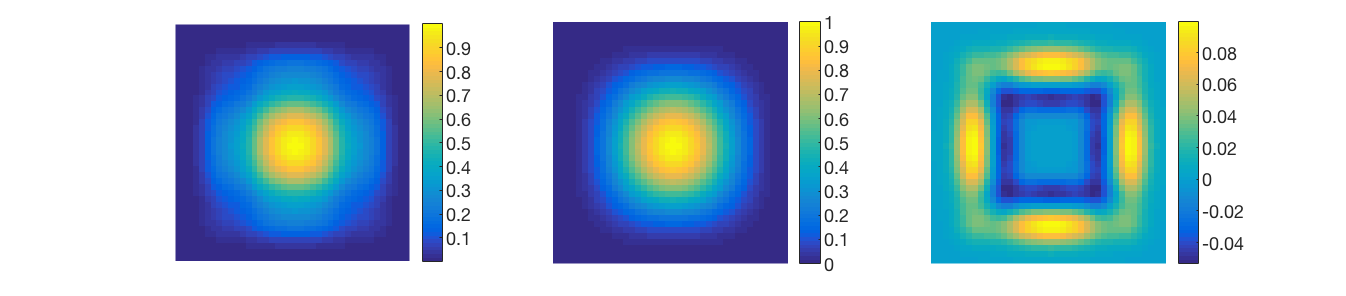}
\captionof{figure}{Left to right: solution of \eqref{eq: opt-1d} in 2D, ground truth and their difference}
\label{fig: 2d-compare-1}
\end{minipage}\\
\vspace*{2em}
\begin{minipage}{.8\textwidth}
\centering
\includegraphics[width = \textwidth]{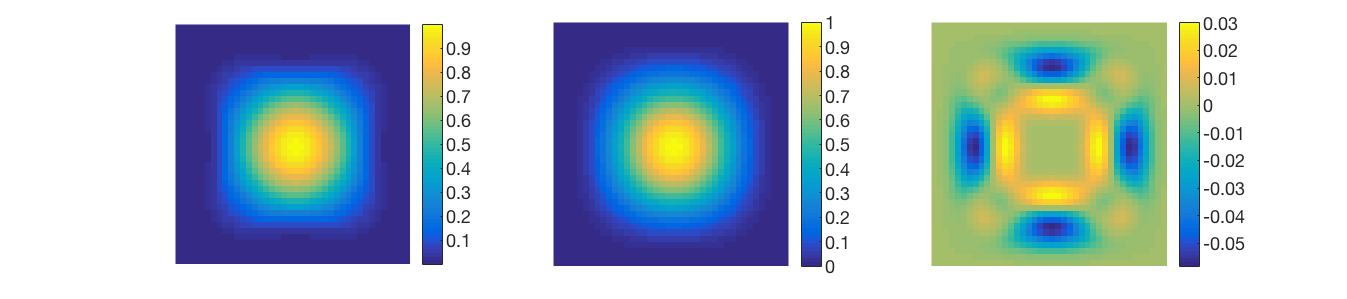}\\
\includegraphics[width = \textwidth]{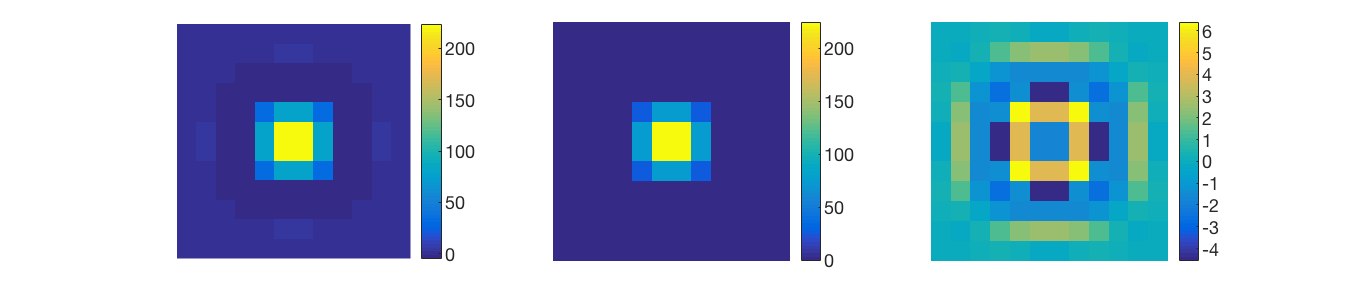}
\captionof{figure}{Left to right: solution of \eqref{eq: opt-2d-weight} ($\lambda = 600$), ground truth and their difference; Top: frequency domain, bottom: time domain.}
\label{fig: 2d-compare-2.2}
\end{minipage}
\end{center}